\documentclass[a4paper,reqno]{amsart}
\usepackage[english]{babel}
\usepackage[utf8]{inputenc}\usepackage{lmodern}
\usepackage[bookmarks,bookmarksopen,bookmarksnumbered,colorlinks,linkcolor=blue,citecolor=red,urlcolor=grey,breaklinks]{hyperref}
\usepackage{graphicx}
\usepackage{tikz}\usetikzlibrary{positioning,decorations.pathmorphing}
\usepackage[orthodox]{nag}
\usepackage{tikz-cd}
\usepackage{mathtools}
\usepackage[mathscr]{euscript}

\theoremstyle{definition}
\newtheorem{Def}{Definition}[section]
\newtheorem*{remark}{Remark}
\newtheorem*{conv}{Convention}

\theoremstyle{plain}
\newtheorem*{Theorem}{Main Theorem}
\newtheorem{Thm}[Def]{Theorem}
\newtheorem{Lem}[Def]{Lemma}
\newtheorem{Cor}[Def]{Corollary}
\newtheorem{Pro}[Def]{Proposition}

\newcommand{\hide}[1]{}

\newcommand{\complex}{\mathbb{C}}
\newcommand{\CC}{\overline{\mathbb{C}}}
\newcommand{\CT}{\mathbb{C}_T}
\newcommand{\R}{\mathbb{R}}
\newcommand{\D}{\mathbb{D}}
\newcommand{\N}{\mathbb{N}}
\newcommand{\Z}{\mathbb{Z}}

\DeclarePairedDelimiter{\abs}{\lvert}{\rvert}
\renewcommand{\Re}{\operatorname{Re}}
\renewcommand{\Im}{\operatorname{Im}}
\newcommand{\M}{\mathscr{M}}

\DeclareMathOperator{\Tr}{Tr}  
\DeclareMathOperator{\It}{It}    
\DeclareMathOperator{\cl}{cl}   
\newcommand{\ad}[1]{\underline{#1}} 
\newcommand{\AD}{\mathcal{S}} 
\newcommand{\itn}[1]{\mathtt{\underline{#1}}} 
\newcommand{\ITN}{\mathcal{S}_\nu}
\newcommand{\Triod}{\mathcal{T}}
\newcommand\restr[2]{\ensuremath{\left.#1\right|_{#2}}} 
\makeatletter
\@namedef{subjclassname@2020}{\textup{2020} Mathematics Subject Classification}
\makeatother

\title[Homotopy Hubbard Trees for exponential maps]{Homotopy Hubbard Trees\\for post-singularly finite exponential maps}

\author{David Pfrang}
\address{prognostica GmbH, Berliner Platz~6, 97080~Würzburg, Germany}
\email{david.pfrang@gmx.de}

\author{Michael Rothgang}
\address{Institut für Mathematik, Humboldt-Universität zu Berlin, Unter den Linden~6, 10099~Berlin, Germany}
\email{rothgami@math.hu-berlin.de}

\author{Dierk Schleicher}
\address{Aix-Marseille Université and CNRS, UMR~7373, Institut de Mathématiques de Marseille, 163~Avenue de Luminy, 13009~Marseille, France}
\email{dierk.schleicher@univ-amu.fr}

\thanks{We gratefully acknowledge that this project was partially supported by a grant from the Deutsche Forschungsgemeinschaft, as well as by the Advanced Grant HOLOGRAM of the European Research Council. Major parts of this work were carried out at Jacobs University Bremen; other parts happened during a visit to Cornell University in the spring of 2018, for which we are grateful as well.}



\begin{document}
\begin{abstract}
We extend the concept of a Hubbard Tree, well established and useful in the theory of polynomial dynamics, to the dynamics of transcendental entire functions. We show that Hubbard Trees in the strict traditional sense, as invariant compact trees embedded in $\mathbb{C}$, do not exist even for post-singularly finite exponential maps; the difficulty lies in the existence of asymptotic values. We therefore introduce the concept of a \emph{Homotopy Hubbard Tree} that takes care of these difficulties.

Specifically for the family of exponential maps, we show that every post-singularly finite map has a Homotopy Hubbard Tree that is unique up to homotopy, and we show that post-singularly finite exponential maps can be classified in terms of Homotopy Hubbard Trees, using a transcendental analogue of Thurston's topological characterization theorem of rational maps.
\end{abstract}

\maketitle

\section{Introduction}

The purpose of this paper is to establish the concept of Hubbard Trees for post-singularly finite exponential maps, i.e., maps of the form $E_\lambda\colon z\mapsto\lambda\exp(z)$ for which the orbit of the singular value $0$ is finite. While Hubbard Trees are known to be an extremely useful concept for the dynamics of polynomials, so far they have not been introduced for the dynamics of entire functions. Part of the problem is that asymptotic values might prevent the existence of Hubbard Trees in the strict traditional sense, as invariant compact trees in the plane, and they do so even for exponential maps. We therefore introduce a modified concept that we call \emph{Homotopy Hubbard Trees} and demonstrate their use for the study of exponential maps. The latter have received a lot of attention over the years, and an understanding of exponential dynamics has often proved to be useful for the study of much more general classes of transcendental entire functions.

This project goes back to the Bachelor's thesis of the second named author \cite{R}.
Previously, Hubbard Trees were defined only for post-critically finite polynomials in \cite{DH1}. The Hubbard Tree is the unique smallest tree embedded in the filled-in Julia set that contains all critical points of the polynomial and is forward invariant under its dynamics (and in addition normalized on bounded Fatou components). It turns out that this definition does not immediately generalize to exponential maps: because of the existence of an asymptotic value, an exponential map cannot have an exactly forward invariant tree containing the post-singular set. We solve this by only requiring the tree to be invariant up to homotopy relative to the post-singular set. Because of the relaxed invariance condition, we call the resulting tree a \emph{Homotopy Hubbard Tree}. In analogy to polynomials, the underlying abstract graph of a Homotopy Hubbard Tree, together with the dynamics of its self-map on the vertices, and a finite amount of extra information, is a useful combinatorial object. Specifically for exponential maps we call this an \emph{abstract exponential Hubbard Tree}. The main results of our paper can be phrased as follows.
\begin{Theorem}
  Every post-singularly finite exponential map has a Homotopy Hubbard Tree, and this tree is unique up to homotopy relative to the post-singular set. Moreover, for every abstract exponential Hubbard Tree, there is a unique post-singularly finite exponential map realizing it.
\end{Theorem}

Hubbard Trees are a convenient tool to read off dynamical properties of the map under consideration (an example is given below). Also, abstract Hubbard Trees provide a way to define the notion of core entropy (introduced by Thurston \cite{ThurstonCoreEntropy} for post-critically finite polynomials) for post-singularly finite parameters and to study properties of the core entropy function on the parameter space. This has been done for quadratic polynomials in \cite{DS} and independently in \cite{TiozzoCoreEntropyQuadratic}.
While entropy of transcendental mappings is always infinite \cite{TranscendentalEntropy}, core entropy of post-singularly finite exponential mappings is always bounded by $\log{2}$ \cite{Ha}. This provides relevant information to the dynamics, in a similar sense that entropy of degree $d$ polynomials always equals $\log{d}$, while core entropy allows to distinguish different polynomial dynamical systems.

We want to point out that the construction of Homotopy Hubbard Trees in this paper is quite explicit. Given the external address of a dynamic ray landing at the singular value, we show how to construct the abstract exponential Hubbard Tree of the corresponding map algorithmically.

Post-singularly finite exponential maps have already been classified in \cite{LSV} in terms of the external addresses of the dynamic rays landing at the singular value. Compared to this previous classification, our result has the advantage that it establishes a bijection between the class of maps and the instances of the combinatorial model: while several addresses might correspond to the same exponential map, there is a natural bijection between abstract exponential Hubbard Trees and post-singularly finite exponential maps. Also, some dynamical properties of the map are obtained more easily from the abstract Hubbard Tree than from an external address. For example, the abstract Hubbard Tree contains information about all periodic branch points of the Julia set, while it is computationally intensive to determine them from the external address.

Although we restrict to exponential maps in this article, our results have wider implications. Indeed, the insights gained here, together with recent findings on the structure of the escaping sets of entire functions with bounded post-singular set established in \cite{BR}, have lead to a general theory of Homotopy Hubbard Trees for all post-singularly finite transcendental entire functions. This is the content of the first named author's PhD thesis \cite{Pf}. The first part of this thesis is available in \cite{PPS}.\\

\textit{Background and relevance.} The dynamics of a holomorphic map is controlled to a large extent by the orbits of its \emph{singular values} (see \cite{M} for a general introduction to holomorphic dynamical systems and \cite{S2} specifically to transcendental entire functions). Singular values are points in the range of the function that do not have a neighborhood on which all branches of the inverse function are well-defined and biholomorphic. For a polynomial $f\colon\complex\to\complex$, every singular value is a \emph{critical value}, i.e., a point $w\in\complex$ such that $f'(z)=0$ for some $z\in f^{-1}(w)$. For a transcendental entire function $f\colon\complex\to\complex$, a singular value can also be an asymptotic value or a limit point of asymptotic and critical values. A point $w\in\complex$ is called an \emph{asymptotic value} if there exists a curve $\gamma\colon[0,+\infty)\to\complex$ satisfying $\lim_{t\to+\infty}\gamma(t)=\infty$ and $\lim_{t\to+\infty}f(\gamma(t))=w$. In any parameter space of holomorphic functions, the easiest maps to understand are those for which all singular values have finite orbits. These functions are called \emph{post-singularly finite} (or \emph{post-critically finite} in the case of polynomials). Not only are they the dynamically simplest maps, but also in many cases the most important ones for the structure of the parameter space.
Figuratively, one reason for studying them is because ``it is easier to search for your lost keys under a lamp post, but there are lamp posts at every important street intersection so they are very helpful to find your way around''.
We explain this by way of analogy to the simplest and best-studied space of polynomial maps,
the space of quadratic polynomials $p_c\colon z\mapsto z^2+c$ for $c\in\complex$. Its connectedness locus is the \emph{Mandelbrot set}
\[
\M:=\{c\in\complex\colon \text{the Julia set}~\mathcal{J}(p_c)~\text{is connected}\}.
\]
All branch points of the Mandelbrot set (in a sense made precise in \cite[Theorem~3.1]{S})  are post-critically finite parameters, and under the assumption that the Mandelbrot set is locally connected (the famous MLC conjecture) these branch points completely describe its topology. In the space of exponential maps, post-singularly finite parameters play a comparable role.

Naturally, a lot of work has gone into investigating the dynamics of post-singularly finite holomorphic functions. For rational maps, a deep characterization theorem by Thurston \cite{DH2} about branched self-covers of the sphere that arise from rational maps has made strong classification results possible. Post-critically finite polynomials \cite{P}, and quite recently post-critically finite Newton maps \cite{LMS1,LMS2}, have been completely classified in terms of custom-tailored combinatorial models using Thurston's theorem.

The fundamental idea of classifying post-critically finite polynomials in terms of their Hubbard Trees originated from \cite{DH1}. Their program has been carried through in greater generality in \cite{BFH} and in full generality in \cite{P}. The following finite amount of combinatorial information is sufficient to completely describe the dynamics of a post-critically finite polynomial:
\begin{itemize}
  \item the graph structure of its Hubbard Tree (without an embedding into the complex plane),
  \item the dynamics of the polynomial on the finite set of vertices of the tree,
  \item the degrees of the critical points of the polynomial,
  \item and for each vertex certain (combinatorial) angles between the edges incident to this vertex.
\end{itemize}
This combinatorial data distinguishes all post-critically finite polynomials. Conversely, if we start with an \emph{expansive} (in a sense defined in Definition~\ref{def:expansive}) dynamical tree with consistent degree and angle information (this finite combinatorial object is known as an \emph{abstract Hubbard Tree}), there exists a post-critically finite polynomial realizing this tree, and this map is unique up to affine conjugation. Therefore, there exists a natural bijection between post-critically finite polynomials (up to affine equivalence) and abstract Hubbard Trees. In the spirit of this classification, we establish a natural bijection between post-singularly finite exponential maps and abstract exponential Hubbard Trees in Section~\ref{sec:classification}.
\\

\textit{Structure of the article.} In Section~\ref{sec:background}, we give a combinatorial description of the escaping set of exponential maps and its dynamics in terms of external addresses. Path-connected components of the escaping set are called \emph{dynamic rays}; topologically, they are arcs terminating at $\infty$ (in $\CC$). The preimages of a dynamic ray landing at the singular value form the boundaries of the \emph{dynamic partition}. Itineraries with respect to this partition distinguish (pre-)periodic points and determine which dynamic rays have a common landing point.

In Section~\ref{sec:definitionHHTs}, we motivate the concept of Homotopy Hubbard Trees and give a precise definition of Homotopy Hubbard Trees for exponential maps. Some technical results on homotopies of embedded trees needed in this paper are discussed.

In Section~\ref{sec:triodAlgorithm}, we show that a Homotopy Hubbard Tree can be chosen to not intersect dynamic rays landing at post-singular points. Furthermore, the itinerary of the singular value with respect to a dynamical partition determines the graph structure and the dynamics of the tree. Together, these two facts imply uniqueness of Homotopy Hubbard Trees.

Section~\ref{sec:separatingRays} deals with the construction of Homotopy Hubbard Trees. The combinatorial description of dynamic rays from Section~\ref{sec:background} is used to show that every triple of post-singular points is separated by dynamic rays landing at a common (pre-)periodic point. Connecting the post-singular set without intersecting these separating rays yields a Homotopy Hubbard Tree.

Finally, in Section~\ref{sec:classification} we give a combinatorial classification of post-singularly finite exponential maps in terms of abstract Hubbard Trees using the transcendental analogue of Thurston's characterization theorem for rational maps established in \cite{HSS}.
In Section~\ref{sec:outlook} we discuss possible extensions and generalizations.\\

\textit{Acknowledgements.}
We would like to thank Mikhail Hlushchanka, Russell Lodge, Dzmitry Dudko, Bernhard Reinke, John Hubbard, and Lasse Rempe for helpful discussions and comments.
We would like to thank the anonymous referee for some helpful comments regarding the presentation.\\

\textit{Notation and terminology.}
The complex plane is denoted by $\complex$ and the Riemann sphere by $\CC$.
We denote the trace of a curve $g$ by $\Tr(g)$. An \emph{arc} is a simple (injective) curve and a \emph{Jordan curve} is a simple closed curve.
A \emph{(pre-)periodic point} is a point that is either periodic or preperiodic under iteration of the map under consideration. We use the term `preperiodic' in the sense of `strictly preperiodic', i.e., excluding the periodic case.

\section{Background on the dynamics of exponential maps}
\label{sec:background}

In the following, let $f\colon\complex\to\complex$ denote an entire function. The \emph{singular set} $S(f)$ of $f$ is defined to be the closure of the set of critical and asymptotic values of $f$. It is the smallest subset of $\complex$ such that the restriction of $f$ to $\complex\setminus f^{-1}(S(f))$ is a covering map (see e.g.\ \cite[Lemma~1.1]{GK}).
The \emph{post-singular set} $P(f)$ of $f$ is given by
\[
P(f):=\overline{\bigcup_{n\geq 0}f^{\circ n}(S(f))}.
\]
We call the function $f$ \emph{post-singularly finite} if $\vert P(f)\vert<\infty$.

Entire functions have an important forward and backward invariant set, the \emph{set of escaping points}; it is often more important than Fatou and Julia sets because it is never empty and never all of $\complex$ (see \cite{E1}), so it provides a non-trivial dynamically invariant decomposition of $\complex$. Here is the formal definition.

\begin{Def}[Escaping Set]
The \emph{escaping set} of $f$ is given by
\[
I(f):=\{z\in\complex\colon f^{\circ n}(z)\to\infty\}.
\]
\end{Def}

In \cite{RRRS} it was shown that for a large class of transcendental entire functions (finite compositions of functions of bounded type and finite order) the escaping
set is organized in the form of disjoint dynamic rays, which are certain arcs consisting of escaping points that terminate at $\infty$.

\begin{Def}[Ray tails, dynamic rays and landing points]
A \emph{ray tail} of $f$ is an injective curve $g\colon [\tau,\infty)\to I(f)$ ($\tau\in\R$) such that for each $n\in\N$ the restriction $\left.f^{\circ n}\right|_g$ is
injective, $\lim_{t\to\infty}f^{\circ n}(g(t))=\infty$, and $f^{\circ n}(g(t))\to\infty$ as $n\to\infty$ uniformly in $t$.

A \emph{dynamic ray} of $f$ is a maximal (in the sense of inclusions of sets) injective curve $g\colon (\tau_0,\infty)\to I(f)$ such that
$\left.g\right|_{[\tau,\infty)}$ is a ray tail for every $\tau>\tau_0$. We say that the dynamic ray $g$ \emph{lands} at the point $z\in\complex$ if $\lim_{t\to \tau_0}g(t)=z$; in this case, $z$ is called the \emph{landing point} of $g$.

We call the dynamic ray $g$ \emph{periodic} if $f^{\circ n}(\Tr(g))\subseteq \Tr(g)$ for some $n\in\N$; we call it \emph{preperiodic} if some forward iterate $f^{\circ k}(g)$ of the ray $g$ (which by definition is again a dynamic ray) is periodic, but not $g$ itself. \label{def:DynamicRays}
\end{Def}

In this paper we focus on the case of \emph{post-singularly finite (psf) exponential maps}: these are (up to affine conjugation) maps of the form $E_\lambda\colon \complex\to\complex$, $E_\lambda(z)=\lambda\exp(z)$, where the parameter $\lambda\in\complex\setminus\{0\}$ is chosen such that the orbit of the only singular value $0$ is finite and hence strictly preperiodic (the orbit of $0$ cannot be periodic because $0$ is an omitted value).
For these functions, the set of escaping points has been described and classified in \cite{SZ1}, using the combinatorial concept of \emph{external addresses} that distinguish dynamic rays.

\begin{Def}[External addresses and the shift map]
An \emph{external address} $\ad{s}$ is a sequence $\ad{s}=s_1 s_2 s_3\ldots$ over the integers. We denote by $\AD=\Z^\N$ the space of all external addresses and by $\sigma\colon \AD\to\AD$, $\sigma(s_1 s_2 s_3\ldots)=s_2 s_3\ldots$, the \emph{left shift map}.
\end{Def}

The shift space $\AD$ can be totally ordered using lexicographic order (for  $\ad{s},\ad{t}\in\mathcal{S}$ we write $\ad{s}< \ad{t}$ if and only if $\ad{s}=s_1 s_2 s_3\dots$ and $\ad t=t_1 t_2 t_3\dots$ have a $k$ so that $s_1=t_1$, \dots, $s_k=t_k$ and $s_{k+1}<t_{k+1}$).
The lexicographic order defines the order topology on $\AD$. Any total order induces a \emph{cyclic order} on the same set: for distinct $\ad{s},\ad{t},\ad{u}\in\AD$ we write
\[  \ad{s}\prec\ad{t}\prec\ad{u} :\Leftrightarrow
 (\ad{s}<\ad{t}<\ad{u}) \lor (\ad{t}<\ad{u}<\ad{s}) \lor (\ad{u}<\ad{s}<\ad{t}). \]

We write $(\ad{s},\ad{t})$ for the open interval between $\ad{s}$ and $\ad{t}$ w.r.t.\ the cyclic order on $\AD$, i.e., we have $\ad{u}\in(\ad{s},\ad{t})$ if and only if $\ad{s}\prec\ad{u}\prec\ad{t}$.

The following theorem is a weak version of the classification result proved in \cite{SZ1}, but it is all we need for this work.

\begin{Pro}[The escaping set of a post-singularly finite exponential map]
Let $E_\lambda$ be a psf exponential map. Then every escaping point either lies on a unique dynamic ray or is the landing point of a unique ray. In particular, distinct dynamic rays are disjoint.

We can assign to each (pre-)periodic external address $\ad{s}$ a dynamic ray $g_{\ad{s}}\colon (0,\infty)\to\complex$ in such a way that
\[ E_\lambda(g_{\ad{s}}(t))=g_{\sigma(\ad{s})} (F(t)) \quad  \text{for all } t>0 \]
where $F(t)=e^t-1$, and so that
\begin{equation}
\lim_{t\to\infty} \Re g_{\ad s}(t) = +\infty
\quad and \quad
\lim_{t\to\infty} \Im g_{\ad s}(t) = -\Im \log\lambda+2\pi s_1 \;.
\label{eq:LimitDynRay}
\end{equation}
Here, we chose the branch of the logarithm for which $\Im \log\lambda\in(-\pi,\pi]$.
\label{pro:escaping_set}\end{Pro}

\begin{remark}
Note that the landing point of a dynamic ray can be an escaping point; in this case there cannot be another dynamic ray with the same landing point. In contrast, a non-escaping point can be the landing point of several dynamic rays. Indeed, (pre-)periodic rays landing together play an important role in this paper.
\end{remark}

We partition the complex plane into horizontal strips of the form
\[
S_k := \left\{ z\in\complex\colon -\Im \log\lambda-\pi+2\pi k < \Im z < -\Im \log\lambda+\pi+2\pi k \right\}
\]
with $k\in\Z$. This partition is called the \emph{static partition} for $E_\lambda$.

Equation \eqref{eq:LimitDynRay} implies that every dynamic ray is, for all sufficiently large potentials $t$, contained in a single sector of the static partition, i.e., there exists a $k\in\Z$ such that $g(t)\in S_k$ for $t$ large enough.

\begin{Def}[External address of a dynamic ray]
Let $g$ be a dynamic ray for a psf exponential map $E_\lambda$. The \emph{external address} $\text{Ad}(g)=s_1 s_2\ldots\in\AD$ \emph{of the dynamic ray} $g$ is the unique external address such that for every $k\in\N$ we have
\[ E_\lambda^{\circ (k-1)}(g(t))\in S_{s_k}  \quad\text{for } t \text{ large enough (depending on $k$)}. \]
\end{Def}

\begin{remark}
We see that $\text{Ad}(g_{\ad{t}})=\ad{t}$ for the dynamic rays $g_{\ad{t}}$ defined in Proposition~\ref{pro:escaping_set}. Different dynamic rays have different external addresses (as follows from the full version of the classification result in \cite{SZ1}), but not every external address is the address of a dynamic ray. A sequence $\ad s\in\AD$ occurs as the external address of a dynamic ray if and only if it is \emph{exponentially bounded} (see \cite[Definition~4.1 and Theorem~4.2]{SZ1}).
\end{remark}

Since dynamic rays are disjoint and converge to $\infty$ in a controlled way as described in \eqref{eq:LimitDynRay}, they have a well-defined vertical order, defined as follows.

\begin{Def}[Vertical order of dynamic rays]
 Let $g$ and $g'$ be two dynamic rays of $E_\lambda$. Then for sufficiently large $\xi>0$, the ray $g$ disconnects the right half plane $\{z\in\complex\colon \Re(z)>\xi \}$ into exactly two unbounded parts (plus possibly some bounded ones), and the curve $g'(t)$ must be contained in a single one of them for all sufficiently large $t$. We say that \emph{$g$ lies above $g'$} if $g'$ (for large potentials) is contained in the lower of these two unbounded components, and write $g>g'$; otherwise we say that \emph{$g$ lies below $g'$} and write $g<g'$.
\end{Def}

It follows easily from the mapping properties of the exponential map that the vertical order of dynamic rays coincides with the lexicographical order of their external addresses (see \cite[Lemma~3.9]{FS}).

\begin{Lem}[Order of rays and external addresses]
For any two dynamic rays $g_{\ad{t}}$ and $g_{\ad{t}'}$ of $E_\lambda$, the ray $g_{\ad{t}}$ lies above $g_{\ad{t}'}$ if and only if $\ad{t}> \ad{t}'$ in lexicographic ordering.\qed
\label{lem:verticalLexicographicOrder}\end{Lem}

As described above for the space of external addresses, the linear order induces a cyclic order on the set of dynamic rays; we write $g\prec g'\prec g''$ if $g'$ lies between $g$ and $g''$ in this cyclic order. The map assigning to each ray its external address preserves the cyclic order (as it preserves the linear order inducing it).

For the purposes of this paper, we are most interested in (pre-)periodic dynamic rays and their landing behavior (compare Definition~\ref{def:DynamicRays}). It was shown in \cite[Theorem~3.2]{SZ2} that every (pre-)periodic dynamic ray of a psf exponential map that is not eventually mapped onto a ray landing at the singular value $0$ lands at a (pre-)periodic point. Conversely, every (pre-)periodic point is the landing point of at least one (pre-)periodic dynamic ray by \cite[Theorem~5.3]{SZ2}. This result is in analogy to the Douady-Hubbard landing theorem for polynomials.

\begin{Thm}[A landing theorem, \cite{SZ2}]
	For every post-singularly finite exponential map, every periodic or perperiodic point is the landing point of at least one and at most finitely many periodic, respectively preperiodic, dynamic rays. Rays landing at the same point have the same preperiod and period.\qed
\label{thm:LandingTheorem}\end{Thm}

We construct a partition of the plane that allows us to describe in combinatorial terms which dynamic rays are landing together: Choose a preperiodic dynamic ray $g_{\ad{s}}$ landing at $0$.
The preimage $E_\lambda^{-1}(\Tr(g_{\ad{s}}))$ consists of countably many disjoint simple curves which are translates of each other by integer multiples of $2\pi i$. Furthermore, if $g$ is any lift of $g_{\ad{s}}$ by $E_\lambda$, then we have
\[ \lim_{t\to 0}\Re (g(t))=-\infty \quad \text{and} \quad \lim_{t\to \infty}\Re (g(t))=+\infty \]
by the mapping properties of exponential maps.

We define a \emph{sector} to be a connected component of $\complex\setminus E_\lambda^{-1}(\Tr(g_{\ad{s}}))$. These sectors partition $\complex$ so that the sector boundaries are exactly the preimages of the ray $g_{\ad{s}}$.

Since distinct dynamic rays are disjoint, there is a unique sector $D_0$ containing the ray $g_{\ad{s}}$ and its landing point $0$. For $k\in\Z$, the sector obtained by translating $D_0$ by $2k\pi i$ is called $D_k$. Observe that for all $k\in\Z$, the restriction $E_\lambda\colon D_k\to \complex\setminus(\Tr(g_{\ad{s}})\cup\{0\})$ is biholomorphic.
We call $\mathcal{D}:=\bigcup_{k\in\Z}\{D_k\}$ the \emph{dynamical partition} for $E_\lambda$ w.r.t.\ $g_{\ad{s}}$ and we call $\partial\mathcal{D}:=E_\lambda^{-1}(\Tr(g_{\ad{s}}))$ the \emph{boundary} of the partition $\mathcal{D}$. Note that $\mathcal{D}$ is a collection of subsets of $\complex$, and $\partial\mathcal{D}\subset\complex$.

As dynamic rays are parametrized by external addresses, the dynamical partition can also be constructed on the level of external addresses. Consider again the shift space $\AD$ endowed with the lexicographic order. In the following, terms like $t_0\ad{t}\in\AD$ with $t_0\in\Z$ and $\ad{t}\in\AD$ will denote concatenation of an integer and an external address.

We start with the external address $\ad{s}=s_1 s_2\ldots$ of the dynamic ray $g_{\ad{s}}$. This address is not constant, so we either have $\ad{s}\in((s_1-1)\ad{s}, s_1\ad{s})$, or $\ad{s}\in(s_1\ad{s}, (s_1+1)\ad{s})$. Denote the interval containing $\ad{s}$ by $I_0$.
For all $k\in\Z$ we define the interval
\[ I_{k}:=\{ t_1 t_2 t_3\ldots\in\AD\colon (t_1-k)t_2 t_3\ldots\in I_0 \}.\]
Observe that $\mathcal{I}:=\bigcup_{k\in\Z} \{I_{k}\}$ is a partition of the shift space $\AD$, which we call the \emph{dynamical partition} of $\AD$ w.r.t.\ $\ad{s}$. We denote by $\partial\mathcal{I}:=\sigma^{-1}(\ad{s})\subset\AD$ the \emph{boundary} of the partition $\mathcal{I}$.

It follows from Lemma~\ref{lem:verticalLexicographicOrder} that a dynamic ray $g$ is contained in the sector $D_k$ of the dynamical partition $\mathcal{D}$ of the complex plane if and only if its address $\text{Ad}(g)\in\AD$ is contained in the sector $I_k$ of the dynamical partition of the shift space.
Recording the sectors in $\mathcal{D}$ to which a point in the plane is mapped under iteration of $E_\lambda$ yields so-called itineraries. The same applies to the partition $\mathcal{I}$ of the shift space.

\begin{Def}[Itineraries of external addresses and the kneading sequence]~\\
 Let $\star$ be a formal symbol not contained in $\Z$ and let $\ad{t}\in\AD$ be an external address. Then the \emph{itinerary} $\It(\ad{t}~\vert~\ad{s})$ of $\ad{t}$ \emph{with respect to $\ad{s}$} is the unique sequence $\It(\ad{t}~\vert~\ad{s})=\itn{t}=\mathtt{t}_1\mathtt{t}_2\ldots\in (\Z\cup\{\star\})^\N$ such that
 \[
 \mathtt{t}_k:=
 \begin{cases}
 i  & \text{if } \sigma^{\circ(k-1)}(\ad{t})\in I_i,\\
 \star & \text{if } \sigma^{\circ(k-1)}(\ad{t})\in \partial\mathcal{I}.
 \end{cases}
 \]
 We call the itinerary $\nu:=\It(\ad{s}~\vert~\ad{s})$ the \emph{kneading sequence} of the external address $\ad{s}$. \label{def:Itineraries}
\end{Def}

\begin{Def}[Itineraries of points and dynamic rays]
 Let $\star$ be a formal symbol not contained in $\Z$. For a point $z\in\complex$, we define the \emph{itinerary} of $z$ w.r.t. $\mathcal{D}$ to be the sequence $\It(z~\vert~\mathcal{D})=\itn{t}=\mathtt{t}_1\mathtt{t}_2\ldots\in (\Z\cup\{\star\})^\N$ such that
  \[
 \mathtt{t}_k:=
 \begin{cases}
 i      & \text{if } E_\lambda^{\circ (k-1)}(z)\in D_i,\\
 \star & \text{if } E_\lambda^{\circ (k-1)}(z)\in \partial\mathcal{D}.
 \end{cases}
 \]

\noindent Let $g$ be a dynamic ray of $E_\lambda$. As distinct dynamic rays do not intersect, and the partition boundary $\partial\mathcal{D}$ consists of rays, every iterate of $g$ is either fully contained in a single sector of the dynamical partition or part of the partition boundary. Therefore, all points on $g$ have equal itineraries and we set $\It(g~\vert~\mathcal{D}):=\It(w~\vert~\mathcal{D})$ for an arbitrary point $w\in\Tr(g)$.
\end{Def}

Itineraries are written in a different font to distinguish them from external addresses. Observe that the itinerary of a (pre-)periodic point is itself (pre-)periodic. Furthermore, it does not contain the symbol $\star$ because $\partial\mathcal{D}$ consists of escaping points.

The following proposition describes the landing behavior of (pre-)periodic dynamic rays in terms of itineraries. It is crucial for the construction of Homotopy Hubbard Trees for exponential maps. A proof can be found in \cite[Lemma~5.2 and Theorem~5.3]{SZ2}.

\begin{Pro}[Landing behavior of (pre-)periodic rays, \cite{SZ2}]
Let $E_\lambda$ be psf, let $\ad{s}$ be the external address of a dynamic ray landing at the singular value, and let $\mathcal{D}$ the corresponding dynamical partiton; denote the kneading sequence of $\ad{s}$ by $\nu$.

The map $z\mapsto\It(z~\vert~\mathcal{D})$ is a bijection between (pre-)periodic points and (pre-)periodic itineraries $\itn{t}\in\Z^\N$ satisfying $\sigma^{\circ n}(\itn{t})\neq\nu$ for all $n\in\N$. The (pre-)periodic ray $g_{\ad{t}}$ lands at the (pre-)periodic point $z$ if and only if $\It(z~\vert~\mathcal{D})=\It(\ad{t}~\vert~\ad{s})$.\qed
\label{pro:LandingBehavior}\end{Pro}

One of the main goals of this paper is to classify post-singularly finite exponential maps in terms of their abstract Hubbard Trees. A different classification of psf exponential maps has already been obtained in \cite[Theorems 2.6 and 2.7]{LSV}. 
\begin{Thm}[Classification of post-singularly finite exponential maps, \cite{LSV}]
For every preperiodic external address $\ad{s}$ starting with the entry $0$, there is a unique post-singularly finite exponential map such that the dynamic ray at external address $\ad{s}$ lands at the singular value.

Every post-singularly finite exponential map is associated in this way to a positive finite number of preperiodic external addresses starting with $0$. Two such external addresses $\ad{s}$ and $\ad{s}'$ are associated to the same exponential map if and only if $\It(\ad{s}'~\vert~\ad{s})=\It(\ad{s}~\vert~\ad{s})$.\qed
\label{thm:SpiderClassification}\end{Thm}

We do not use this result, neither for the construction of Homotopy Hubbard Trees for exponential maps nor in the proof of our own classification result. Still, it is a nice insight to have in mind during the upcoming constructions.

\begin{conv}
Let $E_\lambda$ be a post-singularly finite exponential map. For the rest of the paper, $\ad{s}$ will always denote the external address of a dynamic ray landing at the singular value and $\nu$ will denote its kneading sequence.
\end{conv}

\section{Homotopy Hubbard Trees}
\label{sec:definitionHHTs}

Hubbard Trees have been defined for polynomials over thirty years ago in \cite{DH1}. We give here a more conceptual, but equivalent definition of Hubbard Trees which has better chances to be generalized to the case of exponential maps than the original one given in \cite{DH1} (which uses the concept of filled-in Julia sets). In order to avoid difficulties arising from unnecessary generality, we don't give the most general definition, but a definition valid for \emph{unicritical polynomials} (polynomials with only one critical point in the complex plane) with preperiodic critical value. This makes sense from a conceptual viewpoint as exponential maps are a dynamical limit of unicritical polynomials \cite{DGH} and for post-singularly finite exponential maps the unique singular value is preperiodic.

Let us make the concept of an embedded tree precise. Our definition differs from the standard definition (based on a topological quotient) in the case of infinite trees: infinite trees occur naturally as part of our construction, and the topology of these trees at vertices of infinite degree deviates from the usual quotient topology.
\begin{Def}[Embedded graphs and trees]
A topological space $G$ is called a \emph{topological graph} if it is homeomorphic to a space
\[ X=\left(\dot\bigsqcup_{i\in I} [0^{(i)}, 1^{(i)}] \right)/_\sim \]
of disjoint copies of the unit interval $[0, 1]$, where $I$ is some index set (for our purposes, we may restrict to countable index sets) and $\sim$ is an equivalence relation identifying some of the endpoints of the intervals $[0^{(i)}, 1^{(i)}]$.

Away from equivalence classes of infinite cardinality, we define the topology on $X$ to be the usual quotient topology. In particular, if $I$ is finite, we equip $X$ with the usual quotient topology.

If $x\in X$ is a vertex of infinite degree, i.e., an equivalence class of infinite cardinality, there are index sets $I_0, I_1\subset I$, at least one of which is infinite, such that $0^{(i)}\in x$ if and only if $i\in I_0$ and $1^{(i)}\in x$ if and only if $i\in I_1$.
For all $0 < \epsilon < 1$, we set
\[ U_\epsilon := \left(\left(\bigcup_{i\in I} [0^{(i)}, \epsilon^{(i)})\right)\cup \left(\bigcup_{i\in I} (1-\epsilon^{(i)}, 1^{(i)}] \right)\right)/_\sim,\]
and define $\{U_\epsilon \}_{0<\epsilon<1}$ to be a neighborhood basis of $x$ in $X$.

We call $G$ \emph{finite} if the index set $I$ may be chosen to be finite. The topological graph $G$ is called a \emph{topological tree} if it is connected and has trivial fundamental group.
\end{Def}

\begin{Def}[Branch and endpoints in embedded trees]
  A \emph{branch} of the topological tree $H$ at the point $p\in H$ is the closure of a connected component of $H\setminus\{p\}$. We denote the number of branches of $H$ at $p$ by $\text{deg}_H(p)$ and call $p$ a \emph{branch point} of $H$ if $\text{deg}_H(p)\geq 3$. If $\text{deg}_H(p)=1$, we call $p$ an \emph{endpoint} of $H$.
\end{Def}

\begin{Def}[Hubbard Trees for preperiodic unicritical polynomials]

Let $p$ be a post-critically finite unicritical polynomial with preperiodic critical value. The \emph{Hubbard Tree $H\subset\complex$} of $p$ is the unique finite embedded tree such that:
\begin{itemize}
\item $P(p)\subset H$, i.e., the forward orbit of every critical point of $p$ is contained in the tree.
\item All endpoints of $H$ are contained in $P(p)$.
\item $p(H)\subset H$, i.e., $H$ is forward invariant under the dynamics of $p$.
\end{itemize}
\label{def:HubbardPoly}\end{Def}

The naive way to generalize this definition to the case of exponential maps fails because of the existence of an asymptotic value.

\begin{Thm}[Hubbard Trees of exponential maps must contain escaping points]
 Let $E_\lambda$ be a post-singularly finite exponential map. There does not exist a finite embedded tree $H\subset\complex$ that is forward invariant under the dynamics of $E_\lambda$ and contains $P(E_\lambda)$.
\label{thm:denseEscapingPoints}\end{Thm}
\begin{proof}
 Let $\ad{s}$ be the external address of a dynamic ray landing at the singular value $0$ and let $\mathcal{D}$ be the dynamical partition for $E_\lambda$ w.r.t.\ $g_{\ad{s}}$. Furthermore, let $p,q\in H$, $p\neq q$, be arbitrary points of well-defined and different itineraries w.r.t.\ $\mathcal{D}$. For example, one can easily verify that there are two points on the forward orbit of $0$ which differ by a multiple of $2\pi i$ and hence lie in different sectors of the dynamical partition $\mathcal{D}$, so one could take $p$ and $q$ to be those two points. Since $H$ is a tree, there is a unique arc $\gamma\colon[0,1]\to H$ connecting $p$ and $q$. Since $H$ is forward-invariant, for every $n\geq 0$ the points $E_\lambda^{\circ n}(p)$ and $E_\lambda^{\circ n}(q)$ are connected within $H$ by $E_\lambda^{\circ n}(\gamma)$.
 By hypothesis, for some $k\in\N$ the points $E_\lambda^{\circ k}(p)$ and $E_\lambda^{\circ k}(q)$ lie in different sectors of the dynamical partition, thus by continuity $E_\lambda^{\circ k}(\gamma)$ crosses the partition boundary and contains an escaping point, contradicting the forward invariance of $H$.
\end{proof}

The idea that leads us to a meaningful definition of Hubbard Trees for exponential maps is to relax the invariance condition: we do not require a Hubbard Tree to be exactly forward invariant, but only invariant up to homotopy relative to the post-singular set. Note that this relaxation is not only necessary (because there is no exactly invariant tree), but also natural from the point of view of Thurston theory: an exponential Hubbard Tree should determine a psf exponential map up to Thurston equivalence, and the Thurston equivalence class of a map is invariant under homotopies relative to the post-singular set in the domain and co-domain of the map (see Definition~\ref{def:ThurstonEquiv}).

Since homotopies rel $P(E_\lambda)$ cannot be pushed forward by $E_\lambda$ (because of the existence of non post-singular preimages of post-singular points), the right way to formulate the forward invariance condition is via backwards iteration. We want to say that the preimage $E_\lambda^{-1}(H)$ of a Hubbard Tree $H$ contains $H$ up to homotopy rel $P(E_\lambda)$. This statement bears a problem. The preimage $E_\lambda^{-1}(H)$ is disconnected because $H$ contains the singular value $0$. Different connected components of $E_\lambda^{-1}(H)$ contain post-singular points which are by definition contained in $H$, and are not allowed to move during the homotopy. Hence, by connectedness, $H$ cannot be homotoped into its preimage $E_\lambda^{-1}(H)$. By adding a preimage $-\infty$ of $0$ to the complex plane, the preimage of $H$ becomes connected and in fact becomes an (infinite) embedded tree, so it makes sense to require that $H$ can be homotoped into $E_{\lambda}^{-1}(H)\cup\{-\infty\}$ rel $P(E_\lambda)$ in the extended plane.

Let us make this idea precise. Let $\CT$ (as a set) be defined as the disjoint union $\CT:=\complex\mathbin{\dot{\cup}}\{-\infty\}$, where for now the point $-\infty$ is just an abstract point not contained in the complex plane. We turn $\CT$ into a topological space by choosing a neighborhood basis $(V_n)$ of $0$ and declaring the sets $U_n:=E_\lambda^{-1}(V_n)\cup\{-\infty\}$ to be a neighborhood basis of $\{-\infty\}$. The extension $\widehat{E_\lambda}\colon \CT\to\complex$ of $E_\lambda$ defined by $\widehat{E_\lambda}(-\infty):=0$ is continuous by definition. The completion $\CT$ of $\complex$ is a special case of a far more general construction of defining a completion of the domain of a holomorphic function by adding all transcendental singularities of its inverse function. See \cite{E2} for further information.
The extended map $\widehat{E_\lambda}$ is not a covering map any more, but we can still lift homeomorphisms and homotopies of the complex plane that fix $0$ under $\widehat{E_\lambda}$.

\begin{Lem}[Lifting homeomorphisms]
Let $A\subset\complex$ be a set containing $0$ and $\varphi\colon\complex\to\complex$ be a homeomorphism which is homotopic to the identity relative to $A$. There exists a unique homeomorphism $\Phi\colon\CT\to\CT$ which is homotopic to the identity relative to $\widehat{E_\lambda}^{-1}(A)$ such that the diagram
\[ \begin{tikzcd}
    \CT \arrow{r}{\Phi} \arrow[swap]{d}{\widehat{E_\lambda}}
                                       & \CT  \arrow{d}{\widehat{E_\lambda}} \\
    \complex \arrow{r}{\varphi} & \complex
\end{tikzcd} \]
commutes. We call $\Phi$ the \emph{preferred lift} of $\varphi$. Every homotopy between $\varphi$ and $\text{id}$ rel $A$ lifts to a homotopy between $\Phi$ and $\text{id}$ rel $\widehat{E_\lambda}^{-1}(A)$.
\label{lem:Lifting}\end{Lem}
\begin{proof}
It follows from the Homotopy Lifting Principle that every homotopy between $\text{id}$ and $\varphi$ rel $A$ lifts to a homotopy between $\text{id}$ and a homeomorphism $\Phi\colon \complex\to\complex$ relative to $E_\lambda^{-1}(A)$.
 Both the map $\Phi$ and the homotopy between $\Phi$ and the identity extend continuously to $\CT$ and fix $-\infty$ since $\varphi$ as well as the homotopy between $\varphi$ and the identity fix $0$.
\end{proof}

In order to state our definition of Homotopy Hubbard Trees, we need a bit more vocabulary.
\begin{Def}[Spanned subtrees]
	For an embedded tree $H\subset X$ and a finite subset $W\subset H$ we write $[W]_H$ (or just $[W]$) for the smallest subtree of $H$ containing $W$. We say that $H$ \emph{is spanned by $W$} if $H=[W]$. Usually, we write $[p_1,p_2,\ldots, p_n]_H$ for the smallest subtree containing the points $p_i\in H$ (omitting the curly brackets).
\end{Def}
Assume that $H\subset\complex$ is a finite embedded tree spanned by $P(E_\lambda)$. By the mapping properties of exponential maps, $\widehat{E_\lambda}^{-1}(H)\subset\CT$ is an infinite embedded tree where the only point of infinite degree is $-\infty$ and every branch of $\widehat{E_\lambda}^{-1}(H)$ at $-\infty$ is homeomorphic to $H$. (This step uses our tailored definition of infinite tree: the topology near the vertex $-\infty$ is compatible with the extension $\CT$.)

As $P(E_\lambda)$ is forward invariant, we have $P(E_\lambda)\subset\widehat{E_\lambda}^{-1}(H)$, and therefore it makes sense to talk about the subtree $H':=[P(E_\lambda)]_{\widehat{E_\lambda}^{-1}(H)}$ of the preimage tree spanned by $P(E_\lambda)$. We say that $H$ is \emph{invariant up to homotopy rel $P(E_\lambda)$} if $H'$ is homotopic to $H$ in the extended plane $\CT$ relative to $P(E_\lambda)$. Let us make the last statement precise.

\begin{Def}[Relative homotopies of embedded trees]
Given a subset $A\subset X$ of the space $X$ and two embedded trees $H,H'\subset X$, we say that $H$ \emph{is homotopic to} $H'$ \emph{relative to (rel) $A$} if there exists a continuous map $I\colon H\times[0,1]\to X$ with the following properties:
	\begin{itemize}
		\item $I^0:=I(\cdot,0)=\text{id}$ and $I^1:=I(\cdot,1)\colon H\to H'$ is a homeomorphism.
		\item For all $p\in H\setminus A$, we have $I(p,t)\in X\setminus A$ for all $t\in[0,1]$.
		\item For all $p\in H\cap A$, the homotopy $I$ is constant on $\{p\}\times[0,1]$.
	\end{itemize}
\label{def:RelativeHomotopies}\end{Def}

As one might expect, this defines an equivalence relation:
if two embedded trees $H$ and $H'$ are homotopic rel $A$ through through $I\colon H\times [0,1]\to X$, also $H'$ and $H$ are homotopic rel $A$ through the \emph{reversed homotopy} $\overline{I}\colon H'\times[0,1]\to X$, $\overline{I}(p,t):=I((I^1)^{-1}(p),1-t)$.
Given a second homotopy $I'\colon H'\times[0,1]\to X$ between $H'$ and an embedded tree $H''\subset X$ rel $A$, we obtain a homotopy between $H$ and $H''$ rel $A$ by forming the \emph{concatenation} $I\cdot I':H\times[0,1]\to X$ via
\[
I\cdot I'(p,t):=
\begin{cases}
I(p,2t) & \text{for}~t\leq\tfrac{1}{2},\\
I'(I^1(p),2(t-\tfrac{1}{2})) & \text{for}~t\geq\tfrac{1}{2}.
\end{cases}
\]
The concatenation satisfies $(I\cdot I')^1=(I')^1\circ I^1$.

Assume from now on that $H$ is invariant up to homotopy rel $P(E_\lambda)$. While $H$ is not forward invariant as a set, we still obtain a self-map of $H$ which is well-defined up to a certain equivalence relation:
the homotopy $I\colon H\times[0,1]\to\CT$ from Definition~\ref{def:RelativeHomotopies} yields an identification of $H$ and $H'$ via the homeomorphism $\psi:=I_1\colon H\to H'$ and the composition $f:=\widehat{E_\lambda}\circ\psi\colon H\to H$ is a self-map of the tree $H$. We call $f$ the \emph{induced self-map} of $H$. There is a distinguished point $v_T:=\psi^{-1}(-\infty)\in H$, the \emph{singular point} of $(H,f)$, with the property that $f$ is not injective at $v_T$ while it is a local homeomorphism elsewhere, and its image $f(v_T)=0$ is the singular value. We call $V_f:=\{f^{\circ n}(v_T)~|~n\in\N_0\}\cup\{p\in H~|~\text{deg}_H(p)\geq 3\}$ the set of \emph{marked points} of $(H,f)$. It is the union of the forward orbit of $v_T$ under $f$ and the set of branch points of $H$, and as such contains $P(E_\lambda)$.

Different homotopies between $H$ and $H'$ yield different self-maps of $H$ and we want to investigate this ambiguity. Let $\varphi\colon H\to H'$ be another identification of $H$ with $H'$ obtained by a homotopy between $H$ and $H'$ rel $P(E_\lambda)$, let $f':=\widehat{E_\lambda}\circ\varphi$ be the corresponding self-map of $H$ and set $v_T':=\varphi^{-1}(-\infty)$. As both maps $\psi$ and $\varphi$ are induced by homotopies relative to $P(E_\lambda)$, we have $\psi(p)=\varphi(p)=p$ for all $p\in P(E_\lambda)$. Therefore the ``change of identification'' $\theta:=\varphi^{-1}\circ\psi\colon  H\to H$ restricts to the identity on $P(E_\lambda)$ and in particular on the set of endpoints of $H$. A (graph-theoretic) isomorphism between finite trees is uniquely determined by its values on endpoints; this implies $\theta(b)=b$ for all branch points $b\in H$. It might happen that $v_T\neq v_T'$, but we still have $\theta(v_T)=v_T'$ by definition of the singular point. We claim that $\theta(f(v))=f(v)$ for all $v\in V_f$. Indeed, if $v\neq v_T$ is a branch point of $H$, then $f$ is a local homeomorphism at $v$, so $f(v)$ is also a branch point of $H$, and we have $\theta(f(v))=f(v)$. If $v\in P(E_\lambda)$, then $f(v)\in P(E_\lambda)$ since $P(E_\lambda)$ is forward invariant under the dynamics of $f$, so again we have $\theta(f(v))=f(v)$. Finally, if $v=v_T$, then $f(v_T)=0\in P(E_\lambda)$, so we have $\theta(f(v))=f(v)$. It follows that $\theta$ restricts to a conjugation between $f$ and $f'$ on $V_f$ because $\theta(f(v))=f(v)=f'(\theta(v))$ for all $v\in V_f$.
\begin{figure}[ht]
	\includegraphics[width=\textwidth]{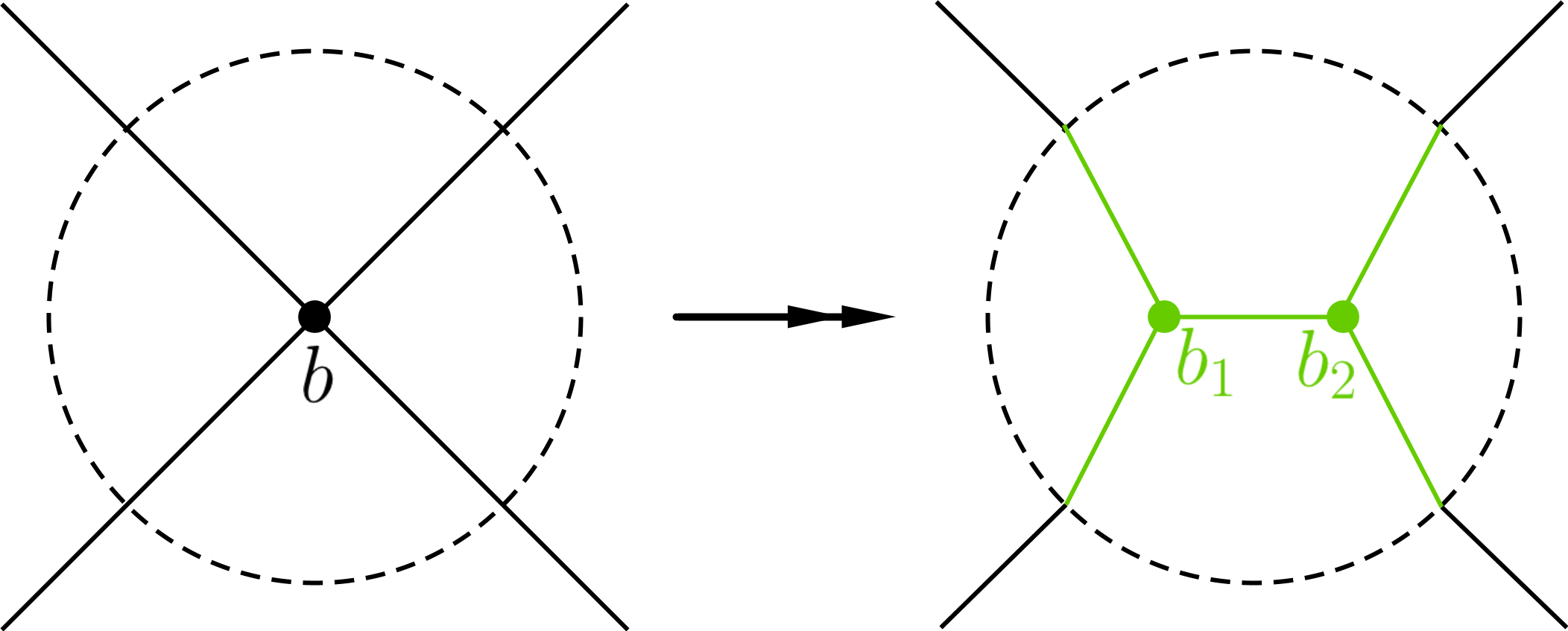}
	\caption{Replacing one degree four branch point $b$ by two degree three branch points $b_1$ and $b_2$.}
\label{fig:SplitBranchpoint}\end{figure}

It turns out that if we just replace the forward invariance condition of Hubbard Trees by the weaker condition of being invariant up to homotopy rel $P(E_\lambda)$, the graph structure of the resulting tree is not uniquely determined. For example, if there exists a periodic branch point $b\in H$ of degree four, one can split it into two degree three branch points $b_1$ and $b_2$, changing the tree only in an arbitrary small neighborhood of $b$ as indicated in Figure~\ref{fig:SplitBranchpoint}. If we perform this change consistently along the (forward and backward) orbit of $b$ under the self-map $f$, we again obtain an invariant tree. Note, however, that in the new tree we have $f^{\circ{n}}([b_1,b_2])=[b_1,b_2]$, where $n$ is the period of $b$ (and hence also of $b_1$ and $b_2$) under $f$, i.e., the in-tree connection of $b_1$ and $b_2$ is forward invariant.

To obtain uniqueness of the graph structure, we want $f$ to be expansive in the following sense.
\begin{Def}[Expansive self-map]
  The self-map $f\colon H\to H$ is called \emph{expansive} if for every pair of marked points $p,q\in V_f$ with $p\neq q$ there exists an $n\geq 0$ such that $v_T\in f^{\circ n}([p,q])$.
\label{def:expansive}\end{Def}

Since different self-maps of $H$ are conjugate on their sets of marked points, the definition of expansivity does not depend on the choice of homotopy between $H$ and $H'$. Again, the motivation for this definition comes from the polynomial case: the map obtained by restricting a unicritical polynomial with preperiodic critical value to its Hubbard Tree is always expansive in the sense that some iterated image of the in-tree connection of two marked points contains the unique critical point.

Summing up, we obtain the following definition of Homotopy Hubbard Trees for exponential maps:

\begin{Def}[Homotopy Hubbard Trees for exponential maps]
Let $E_\lambda$ be a post-singularly finite exponential map. A \emph{Homotopy Hubbard Tree} $H\subset\complex$ for $E_\lambda$ is a finite embedded tree $H$ such that:
 \begin{itemize}
   \item $H$ is spanned by $P(E_\lambda)$.
   \item $H$ is invariant up to homotopy:
   \[ H':=[P(E_\lambda)]_{\widehat{E_\lambda}^{-1}(H)}\subset\CT \text{ is homotopic to } H \text{ relative to } P(E_\lambda).\]
   \item The induced self-map of $H$ is expansive.
 \end{itemize}
\label{def:HomotopyHubbardTree}\end{Def}

Homotopy Hubbard Trees for exponential maps are only required to be invariant up to homotopy. In order to prove a meaningful uniqueness statement, we have to define a suitable equivalence relation on Homotopy Hubbard Trees that deals with the increased flexibility compared to the polynomial case. Naturally, being a Homotopy Hubbard Tree is a property of homotopy classes of embedded trees relative to the post-singular set.

\begin{Lem}[Equivalent Homotopy Hubbard Trees]
Let $H,\tilde{H}\subset\complex$ be finite embedded trees spanned by $P(E_\lambda)$ and assume that $H$ and $\tilde{H}$ are homotopic rel $P(E_\lambda)$. Then $H$ is a Homotopy Hubbard Tree for $E_\lambda$ if and only if $\tilde{H}$ is a Homotopy Hubbard Tree for $E_\lambda$.\label{lem:Transformation}\end{Lem}
\begin{proof}
By symmetry of Definition~\ref{def:RelativeHomotopies}, both directions are equivalent, and it suffices to show one of them.
Assume that $H$ is a Homotopy Hubbard Tree for $E_\lambda$, and let $I_0$ be a homotopy between $H$ and $\tilde{H}$ rel $P(E_\lambda)$. In analogy to Lemma~\ref{lem:Lifting}, one shows that $I_0$ lifts to a homotopy $\widehat{I}_1\colon\widehat{E_\lambda}^{-1}(H)\times[0,1]\to\CT$ between the preimages of $H$ and $\tilde{H}$ rel $\widehat{E_\lambda}^{-1}(P(E_\lambda))$. As $P(E_\lambda)\subset\widehat{E_\lambda}^{-1}(P(E_\lambda))$, the homotopy $\widehat{I}_1$ fixes $P(E_\lambda)$ pointwise, so the restriction
$I_1:=\restr{\widehat{I}_1}{H'\times[0,1]}$
is a homotopy between $H'$ and $\tilde{H'}$ rel $P(E_\lambda)$. By the invariance of $H$, there exists a homotopy $I\colon H\times[0,1]\to\CT$ between $H$ and $H'$. We see that $\tilde{H}$ is homotopic to $\tilde{H}'$ rel $P(E_\lambda)$ via the concatenation $\overline{I_0}\cdot I\cdot I_1$. The self-map $\tilde{f}:=\widehat{E_\lambda}\circ I_1^1\circ I^1 \circ\overline{I_0}^1$ of $\tilde{H}$ obtained via this homotopy is conjugate to the self-map $f:=\widehat{E_\lambda}\circ I^1$ of $H$. This can be seen using the fact that $\widehat{E_\lambda}\circ I_1^1=I_0^1\circ \widehat{E_\lambda}$ because $I_1^1$ is a lift of $I_0^1$. Indeed,
\[
\tilde{f}=\widehat{E_\lambda}\circ I_1^1\circ I^1 \circ\overline{I_0}^1 = I_0^1\circ\widehat{E_\lambda}\circ I^1\circ\overline{I_0}^1
=I_0^1\circ f\circ (I_0^1)^{-1}.
\]
Therefore, the expansivity of $\tilde{f}$ follows from the expansivity of $f$.
\end{proof}

Let us look more closely at the homotopy involved in the invariance condition of Homotopy Hubbard Trees. We want to see that after a small modification this homotopy can be replaced by a stronger kind of homotopy called an ambient isotopy.

\begin{Def}[Ambient isotopies]
	Let $X$ be a topological space and $A\subset X$ a subspace.
  Two embedded trees $H,H'\subset X$ are called \emph{ambient isotopic} relative to $A$ if there exists a continuous map $I\colon X\times[0,1]\to X$ such that the following conditions are satisfied:
    \begin{itemize}
    	\item $I(\cdot,0)=\text{id}$ and $I(H\times\{1\})=H'$.
    	\item For each $a\in A$, the homotopy $I$ is constant on $\{a\}\times[0,1]$.
    	\item For each $t\in[0,1]$, the time-$t$ map $I(\cdot,t)\colon X\to X$ is a homeomorphism.
    \end{itemize}
    Two homeomorphisms $\varphi_0,\varphi_1\colon X\to X$ are called \emph{isotopic} relative to $A$ if there exists a homotopy between $\varphi_0$ and $\varphi_1$ which is constant on $A$ and restricts to a homeomorphism for each fixed $t\in[0,1]$.\label{def:Isotopies}
\end{Def}

The distinction between relative homotopy and isotopy is known to be subtle; in general, $\phi$ being homotopic to the identity $\text{id}$ rel $A$ does not imply $\phi$ being isotopic to $\text{id}$ rel $A$. In our setting, however, these notions are equivalent; we will use this in Section~\ref{sec:classification}.

\begin{Pro}[{e.g.\ \cite[Theorem~1.12]{FM}}]
  Let $S$ be a closed oriented surface (i.e., compact without boundary) and $A\subset S$ a finite set of marked points. An orientation-preserving homeomorphism $\phi\colon S\to S$ is homotopic to $\text{id}$ rel $A$ iff it is isotopic to $\text{id}$ rel $A$.\qed
\label{pro:homotopyImpliesIsotopy}\end{Pro}
For later use (in the proof of Proposition~\ref{pro:AccessesContainRays}), we want to prove the following result regarding the homotopy involved in Definition~\ref{def:HomotopyHubbardTree}:

\begin{Lem}[On the invariance condition]
	Let $H$ be a Homotopy Hubbard Tree for $E_\lambda$. For every neighborhood $U\subset\CT$ of $-\infty$, there exists an embedded tree $H''\subset\complex$ such that $H''$ is homotopic to $H'$ rel $\CT\setminus U$, and $H''$ is ambient isotopic to $H$ in $\complex$ rel $P(E_\lambda)$.
  \label{lem:ExtendedIsotopy}
\end{Lem}

One might hope to obtain an ambient isotopy between $H$ and $H'$ rel $P(E_\lambda)$. Unfortunately, such an isotopy does not exist. The tree $H'$ contains $-\infty$, while $H$ does not, so the ambient isotopy would have to send $-\infty$ to some point in the complex plane. However, the space $\CT$ is not locally compact at $-\infty$, while it is locally compact at every point in the complex plane, so there does not exist a homeomorphism of $\CT$ sending $-\infty$ to a point in the complex plane.

Therefore, we have to take an intermediate step, homotoping \emph{only the tree} $H'$ to a tree $H''\subset\complex$. Then, we use a classical result on homotopies between graphs on surfaces to find an ambient isotopy between the modified tree $H''$ and $H$.

\begin{Lem}[Homotoping $-\infty$ into the plane]
	Let $H\subset\complex$ be a finite embedded tree spanned by $P(E_\lambda)$, and assume that $0$ is an endpoint of $H$. Let $H'=[P(E_\lambda)]_{\widehat{E_\lambda}^{-1}(H)}\subset\CT$ be the subset of the preimage tree spanned by $P(E_\lambda)$. For every neighborhood $U$ of $-\infty$ in $\CT$, there is a homotopy between $H'$ and an embedded tree $H''\subset\complex$ relative to $\CT\setminus U$.
\label{lem:HomotopingMinusInfty}\end{Lem}
\begin{proof}
	The image $\tilde{V}:=\widehat{E_\lambda}(U)$ is a neighborhood of $0$. As $0$ is an endpoint of $H$, we can find a Jordan domain $V\subset\tilde{V}$ such that $H\cap V=\Tr(\gamma)$ for some arc $\gamma\colon[0,1]\to\overline{V}$ connecting $0$ to some point $v\in\partial V$ and satisfying $\gamma([0,1))\subset V$.
  The preimage $\tilde{U}:=\widehat{E_\lambda}^{-1}(V)\subset U$ is a neighborhood of $-\infty$ and $\partial U\subset\complex$ is a $2\pi i$-periodic arc.
  The preimage curves of $\gamma$ have a natural vertical order, and only finitely many of them are contained in $H'$. Let $\gamma_u$ be the uppermost and $\gamma_l$ be the lowermost preimage curve of $\gamma$ contained in $H'$.
  The arcs $\gamma_u$ and $\gamma_l$ together with the part of $\partial\tilde{U}$ between their endpoints on $\partial\tilde{U}$ bound a simply connected domain $D\subset\tilde{U}$ which is contained in $\complex$.
  The map $\Phi\colon\CT\to\CC$ defined by $\Phi(z)=z$ for $z\in\complex$, and $\Phi(-\infty)=\infty$, is continuous and injective, but its inverse is not continuous. The domain $\Phi(D)$ is a Jordan domain in $\CC$, and it is easy to see that $\tilde{H}':=\Phi(H')$ is homotopic to a tree $\tilde{H}''\subset\complex$ rel $\CC\setminus\cl_{\CC}(D)$.
  The closure $\text{cl}_{\CT}(D)$ in the extended plane $\CT$ is homeomorphic to its closure $\text{cl}_{\CC}(D)$ on the Riemann sphere because every sequence $(z_n)_n\subset D$ with $\lim_{n\to\infty}\abs{z_n}=\infty$ satisfies $\lim_{n\to\infty}\Re(z_n)=-\infty$. (Note that this does not hold if we just require $(z_n)_n\in U$ because the sequence could escape in vertical direction.) Therefore, the homotopy between $\tilde{H}'$ and $\tilde{H}''$ in $\CC$ is pushed forward to a homotopy between $H'$ and some $H''\subset\complex$ relative to $\CT\setminus\text{cl}_{\CT}(D)$.
\end{proof}

It remains to show that $H''$ is ambient isotopic to $H$ rel $P(E_\lambda)$ in $\complex$. The question under which conditions homotopic embedded graphs on a surface are ambient isotopic has already been studied extensively. One rather general result can be found in \cite[Lemma~2.9]{FM}. We state a weaker version of this result that is sufficient for our purposes.

\begin{Lem}[Isotopies of curve systems on marked spheres]
	Let $A\subset\CC$ be a finite set. Let $\gamma_1$, \ldots, $\gamma_n$ be a collection of pairwise non-homotopic proper arcs in $\CC$ that do not intersect each other except possibly at their endpoints. If $\gamma_1'$, \ldots, $\gamma_n'$ is another such collection so that $\gamma_i'$ is homotopic to $\gamma_i$ relative to $A$ for each $i$, then there exists an ambient isotopy $I\colon\CC\times[0,1]\to\CC$ rel $A$ satisfying $I^1(\Tr(\gamma_i))=\Tr(\gamma_i')$ for all $i$ simultaneously.
\label{lem:CurveSystem}\end{Lem}

Here, a \emph{proper arc} is an arc $\gamma\colon[0,1]\to\CC$ satisfying $\gamma^{-1}(A)=\{0,1\}$, and two arcs are called \emph{homotopic rel $A$} if they are homotopic rel $A$ in the sense of Definition~\ref{def:RelativeHomotopies}. We cannot apply Lemma~\ref{lem:CurveSystem} directly to our setting because not post-singular branch points of $H$ are allowed to move during the homotopy, so there is no decomposition of $H$ into proper arcs. The following lemma allows us to deal with this problem.

\begin{Lem}
Let $D\subset\CC$ be a Jordan domain, and let $H,\tilde{H}\subset\overline{D}$ be finite embedded trees such that $H\cap\partial D=\tilde{H}\cap\partial D=\{p_1,\ldots, p_n\}$ where the $p_i$ are indexed according to their cyclic order on $\partial D$. Assume that there exists a homeomorphism $\varphi\colon H\to\tilde{H}$ such that $\varphi(p_i)=p_i$ for all $i\in\{1,\ldots, n\}$, and that the set of endpoints of $H$ (and therefore also of $\tilde{H}$) is a subset of $\{p_1,\ldots, p_n\}$. Then $H$ and $\tilde{H}$ are ambient isotopic rel $\partial D$.\label{lem:ClosedDiskHomotopy}
\end{Lem}
\begin{proof} We prove a more general statement with $H$ and $\tilde{H}$ replaced by finite unions of pairwise disjoint embedded trees. The points $p_i$ and the branch points of $H$ are the \emph{vertices} of $H$, and the subarcs of $H$ joining adjacent vertices are the \emph{edges} of $H$. The proof works via induction on the number of edges.

There exists a component $T$ of $H$ and an index $i$ such that $p_i,p_{i+1}\in T$. Let $D'$ be the subdomain of $D$ bounded by $[p_i,p_{i+1}]_H$ and $[p_{i+1},p_i]_{\partial D}$. By hypothesis, we have $H\subset\overline{D'}$. Let $\tilde{T}$ be the component of $\tilde{H}$ with the same endpoints. There exists an isotopy of $\overline{D}$ rel $\partial{D}$ mapping $[p_i,p_{i+1}]_{\tilde{T}}$ to $[p_i,p_{i+1}]_T$. Both $H'=\overline{H\cap D'}$ and $\tilde{H}'=\overline{\tilde{H}\cap D'}$ consist of finitely many embedded trees with endpoints $p_1',\ldots, p_m'\in\partial D$, but the total number of edges of these trees has been reduced.
\end{proof}


\begin{Lem}\label{lem:HomotopyImpliesAmbIsotopy}
Let $H$ be a Homotopy Hubbard Tree, and let $H''\subset\complex$ be homotopic to $H'=[P(E_\lambda)]_{\widehat{E_\lambda}^{-1}(H)}$ rel $P(E_\lambda)$ in $\CT$. Then $H''$ is ambient isotopic to $H$ rel $P(E_\lambda)$ in $\complex$.
\end{Lem}
\begin{proof}
  It is easy to see that $H''$ is homotopic to $H$ rel $P(E_\lambda)$ in $\complex$. Hence, it remains to show that $H''$ is ambient isotopic to $H$ relative to $P(E_\lambda)$ in $\complex$. Choose a positively oriented simple closed curve $\Gamma\colon [0,1]\to\complex$ with $\Gamma(0)=\Gamma(1)=0$ containing $P(E_\lambda)$ and satisfying the following properties:
\begin{itemize}
	\item We can index the post-singular set as $P(E_\lambda)=\{0=p_0,p_1,\ldots, p_n=p_0\}$ and subdivide $\Gamma$ into arcs $\gamma_i\colon [t_i,t_{i+1}]\to\complex$ ($i\in\{0,\ldots, n-1\}$) with $\gamma_i(t_i)=p_i$, $\gamma_i(t_{i+1})=p_{i+1}$.
	\item We have $\gamma_i((t_i,t_{i+1}))\cap H=\emptyset$ for all $i\in\{0,\ldots, n-1\}$.
	\item $\gamma_i$ is homotopic to the arc $[p_i,p_{i+1}]_H$ relative to $P(E_\lambda)$.
\end{itemize}
Hence, $\Gamma$ traverses the boundary of $H$, touching it only once at every post-singular point. In the same way, we choose a simple closed curve $\Gamma''$ that traverses the boundary of $H''$. As $H$ and $H''$ are homotopic relative to $P(E_\lambda)$, $\gamma_i$ and $\gamma_i''$ are homotopic rel $P(E_\lambda)$ by the third property on the list. By Lemma~\ref{lem:CurveSystem}, there is an ambient isotopy $I\colon \CC\times[0,1]\to\CC$ between $\cup\gamma_i$ and $\cup\gamma_i''$ relative to $P(E_\lambda)\cup\{\infty\}$. The curve $\Gamma''$ traverses the boundaries of $H^1:=I^1(H)$ and $H''$ simultaneously. By Lemma~\ref{lem:ClosedDiskHomotopy}, the trees $H^1$ and $H''$ are ambient isotopic rel $\Tr(\Gamma'')\supset P(E_\lambda)$.
\end{proof}

\begin{proof}[Proof of Lemma~\ref{lem:ExtendedIsotopy}]
By Lemma~\ref{lem:BasicProperties}, the singular value $0$ is an endpoint of $H$. Hence, Lemma~\ref{lem:ExtendedIsotopy} is a direct consequence of Lemma~\ref{lem:HomotopingMinusInfty} and Lemma~\ref{lem:HomotopyImpliesAmbIsotopy}.
\end{proof}

\section{The triod algorithm: determining the graph structure}
\label{sec:triodAlgorithm}

In this section, we show that for a fixed post-singularly finite exponential map the structure of a Homotopy Hubbard Tree as a dynamical tree (see Definition~\ref{def:ExpDynamicalTree}) is uniquely determined by the kneading sequence (see Definition~\ref{def:Itineraries}) of the external address of a ray landing at the singular value. Together with the fact that a Homotopy Hubbard Tree does not intersect dynamic rays landing at post-singular points (up to homotopy), this implies uniqueness of Homotopy Hubbard Trees.

If several dynamic rays land at the singular value, the kneading sequences of their external addresses agree by the following lemma, so we speak of \emph{the} kneading sequence of a post-singularly finite exponential map. A proof of the result can be found in \cite[Proof of Theorem~2.7, $(1)\Rightarrow(5)$]{LSV}.

\begin{Lem}[Different dynamical partitions]
	Let $E_\lambda$ be a psf exponential map for which the dynamic rays $g_{\ad{s}}$ and $g_{\ad{s'}}$ both land at $0$. Then the kneading sequences $\It(\ad{s}'~\vert~\ad{s'})$ and $\It(\ad{s}~\vert~\ad{s})$ agree. Therefore, a post-singular point $p\in P(E_\lambda)$ is contained in the $i$-th sector $D_i$ of $\mathcal{D}$ if and only if it is contained in the $i$-th sector $D_i'$ of $\mathcal{D'}$.
\label{lem:DifferentPartitions}\end{Lem}

Let $H$ be a Homotopy Hubbard Tree, and $f$ be an induced self-map of $H$. We will see that two post-singular points are contained in the same branch of $H$ at the singular point $v_T$ if and only if their itineraries start with the same integer $k_i$, and we index the different branches $B_{k_i}$ by the corresponding integers. If we accept this fact for the moment, the triple $(H,f,B_{k_i})$ satisfies the following definition.
\begin{Def}[Exponential dynamical tree]
	An \emph{exponential dynamical tree} is a triple $(H,f,B_{k_i})$ where $H$ is a finite topological tree and $f\colon H\to H$ is a self-map of $H$ satisfying the following conditions:
	\begin{itemize}
		\item There exists a distinguished point $v_T\in H$, the \emph{singular point}, such that $f$ is not injective at $v_T$ whereas $f$ is a local homeomorphism at each point $p\in H\setminus\{v_T\}$.
		\item All endpoints of $H$ are on the singular orbit.
		\item The \emph{singular value} $f(v_T)$ is strictly preperiodic.
		\item If $p,q\in H$ with $p\neq q$ are branch points or points on the singular orbit, then there is an $n\geq 0$ such that $f^{\circ n}([p,q])$ contains the singular point $v_T$ (expansivity).
	\end{itemize}
In addition, the connected components $B_{k_0}, B_{k_1},\ldots, B_{k_n}$ of $H\setminus\{v_T\}$ are indexed by distinct integers $k_i\in \Z$ such that $k_0=0$ and $f(v_T)\in B_0$.

The set $V_f$ of \emph{marked points} consists of the forward orbit of $v_T$ under $f$ (including $v_T$) and the set of branch points of $H$.

We call two exponential dynamical trees $(H,f,B_{k_i})$ and $(\tilde{H},\tilde{f},\tilde{B}_{k_i})$ \emph{equivalent} if there exists a homeomorphism $\varphi\colon H\to\tilde{H}$ that restricts to a conjugation between $f$ and $\tilde{f}$ on the set of marked points $V_f$ and maps each branch $B_{k_i}$ to the branch $\tilde{B}_{k_i}$ of the same index. \label{def:ExpDynamicalTree}
\end{Def}

Let us state some simple properties of exponential dynamical trees first. Their proof is not very hard and works in complete analogy to the proof of \cite[Lemma~2.3]{BKS}, so we omit it here.

\begin{Lem}[Basic properties of exponential dynamical trees]
	Let $(H,f,B_{k_i})$ be an exponential dynamical tree. The singular value $f(v_T)$ is an endpoint of $H$. Each branch point is periodic or preperiodic. The restriction of $f$ to any branch of $H$ at $v_T$ is injective.\qed
\label{lem:BasicProperties}\end{Lem}

We have seen in the previous section that different self-maps of a Homotopy Hubbard Tree $H$ yield equivalent exponential dynamical trees (see the paragraphs after Definition~\ref{def:RelativeHomotopies}). Furthermore, if $\tilde{H}$ and $H$ are equivalent Homotopy Hubbard Trees, we can choose their induced self-maps to be conjugate to each other (see the proof of Lemma~\ref{lem:Transformation}). Therefore, we obtain a well-defined map from the set of equivalence classes of Homotopy Hubbard Trees to the set of equivalence classes of exponential dynamical trees. The main goal of this section is to give a constructive proof of the following statement.

\begin{Thm}
	Let $H$ and $\tilde{H}$ be two Homotopy Hubbard Trees for a post-singularly finite exponential map $E_\lambda$. Then $H$ and $\tilde{H}$ yield equivalent exponential dynamical trees. \label{thm:UniqueGraphStructure}
\end{Thm}

Note that we do not require $H$ and $\tilde{H}$ to be equivalent. The proof of this theorem works in two steps. To explain the outline of the proof, we have to define kneading sequences for exponential dynamical trees also.
\begin{Def}[Itineraries and kneading sequence]
Let $(H,f,B_{k_i})$ be an exponential dynamical tree. The \emph{itinerary} of a point $p\in H$ is the infinite sequence $\It(p)=\itn{t}=\mathtt{t_1 t_2\ldots}$ where
	  \[
	 \mathtt{t}_i:=
	 \begin{cases}
	 i      & \text{if } f^{\circ (k-1)}(p)\in B_{k_i},\\
	 \star & \text{if } f^{\circ (k-1)}(p)=v_T.
	 \end{cases}
	 \]
The itinerary $\It(f(v_T))=:\nu=\nu_1\nu_2\ldots$ of the singular value is called the \emph{kneading sequence} of $(H,f,B_{k_i})$. \label{def:ItinerariesAbstract}
\end{Def}
The proof of Theorem~\ref{thm:UniqueGraphStructure} works in two steps. First we show that different Homotopy Hubbard Trees for the same psf exponential map yield exponential dynamical trees with equal kneading sequences. The second step is to show that an exponential dynamical tree is already determined (up to equivalence) by its kneading sequence.

In order to prove the first statement, we will show that for every Homotopy Hubbard Tree there exists a dynamic ray landing at $0$ such that up to homotopy the tree does not intersect the ray. Therefore, the preimage of the tree does not intersect the boundary of the corresponding dynamical partition (up to homotopy). This implies the statement above since itineraries are independent of the choice of partition by Lemma~\ref{lem:DifferentPartitions} and every partition sector contains at most one branch of $H'=[P(E_\lambda)]_{\widehat{E_\lambda}^{-1}(H)}$ at $-\infty$.

For the proof we need a topological characterization of dynamic rays landing at periodic post-singular points.
Let $p\in P(E_\lambda)$ be periodic. An arc $\gamma\colon [0,1]\to\CC$ with $\gamma(0)=p$, $\gamma(1)=\infty$, and $\Tr(\gamma)\cap P(E_\lambda)=\{p\}$ is called a \emph{leg at $p$}. We denote its homotopy class relative to $P(E_\lambda)\cup\{\infty\}$ by $[\gamma]$. By definition a leg $\tilde{\gamma}$ at $p$ is in the same homotopy class as $\gamma$ if there is a homotopy between $\gamma$ and $\tilde{\gamma}$ in $\CC\setminus(P(E_\lambda)\setminus\{p\})$ fixing the endpoints.

For every leg $\gamma$ we set $\tilde{\mathcal{L}}(\gamma)=\Gamma$, where $\Gamma\colon [0,1]\to\CC$ with $\Gamma(0)=p'$ and $\Gamma(1)=\infty$ is the lift of $\gamma$ starting at the unique periodic preimage $p'$ of $p$. Lifting is always possible as $0\notin\Tr(\gamma)$ and $E_\lambda\colon \complex\to\complex\setminus\{0\}$ is a covering map. As $P(E_\lambda)$ is forward invariant, the arc $\Gamma$ satisfies $\Tr(\Gamma)\cap P(E_\lambda)=\{p'\}$, so $\Gamma$ is a leg at $p'$.

The map $\tilde{\mathcal{L}}$ defined in this way descends to a well-defined map $\mathcal{L}([\gamma])=[\Gamma]$ because homotopies lift under covering maps. We call $\mathcal{L}$ the \emph{leg map}, using the terminology from \cite[Definition~4.3]{B}.

\begin{Pro}[Topological Characterization of dynamic rays]
Let $p\in P(E_\lambda)$ be periodic, and let $\gamma$ be a leg at $p$. There exists a dynamic ray $g\in[\gamma]$ if and only if $[\gamma]$ is periodic under the leg map.

If $[\gamma]$ is periodic under the leg map, the ray $g\in[\gamma]$ is unique. Stated differently, distinct dynamic rays are not homotopic rel $P(E_\lambda)\cup\infty$.
\label{pro:TopRays}\end{Pro}
\begin{proof}
It is clear that $[\gamma]$ is periodic under iteration of $\mathcal{L}$ if there is a dynamic ray $g\in[\gamma]$: every dynamic ray landing at a periodic point is periodic as a set, so $[g]=[\gamma]$ is periodic under the leg map. The proof of the other direction is essentially the same as the proof of \cite[Theorem~4.11]{B}. Every iterate $f:=E_\lambda^{\circ n}$ has the same dynamic rays as $E_\lambda$, so by passing to a suitable iterate $f$ we can assume that $[\gamma]$ is fixed under iteration of the leg map for $f$. Applying the leg map to $[\gamma]$ at least twice, we can further assume that $\gamma$ is eventually contained in a single fundamental domain $F_0$ for $f$ for some static partition of $f$. One can then prove just as in \cite[Theorem~4.11]{B} that $\gamma$ is homotopic to the dynamic ray $g_{\ad{t}}$ of $f$ with external address $\ad{t}=F_0 F_0\ldots$. If there was a second dynamic ray $g'\in[\gamma]$, it would also have external address $\ad{t}$, contradicting the fact that every external address is the address of at most one ray. See \cite{B} for the definition of the terminology and the details of the proof.
\end{proof}

\begin{remark}
We are slightly imprecise in the formulation of Proposition~\ref{pro:TopRays}: a dynamic ray $g\colon (0,\infty)\to\CC$ landing at $p$ neither contains $p$ nor $\infty$ by Definition~\ref{def:DynamicRays}. When talking about the homotopy class of $g$, we treat it as an arc in $\CC$ containing its endpoints $p$ and $\infty$, so that $g$ becomes a leg at $p$.
\end{remark}

We prove that for every access to a Homotopy Hubbard Tree at a post-singular point there exists a dynamic ray that approaches $H$ through this access and does not intersect $H$ up to homotopy. Before stating the theorem, let us formally define what we mean by an access.

\begin{Def}[Accesses]
Let $H$ be a Homotopy Hubbard Tree for $E_\lambda$. An \emph{access} to $H$ at $p\in P(E_\lambda)$ is a prime end $A=[\{N_j\}]$ for the pair $(\CC\setminus H, H)$ with impression $\{p\}$ (see \cite[Section~17]{M}).

We say that a leg $\gamma\colon [0,1]\to\CC$ \emph{approaches $p$ through $A$} if $\gamma(0)=p$, $\gamma((0,1])\cap H=\emptyset$, and for every $j\in\N$ we have $\gamma(t)\in N_j$ for $t>0$ small enough.

\end{Def}

\begin{Pro}[Accesses contain dynamic rays]
Let $H$ be a Homotopy Hubbard Tree for $E_\lambda$. There exists an equivalent Homotopy Hubbard Tree $\tilde{H}\cong H$ such that for each post-singular point $p\in P(E_\lambda)$ and each access $A$ to $p$ in $\tilde{H}$ there exists a dynamic ray $g$ approaching $p$ through $A$.
\label{pro:AccessesContainRays}\end{Pro}

\begin{proof}
Let $\{0=p_0,p_1=E_\lambda(p_0),\ldots,p_n=E_\lambda^{\circ n}(p_0), p_{n+1}=p_{l+1}\}$ be the post-singular orbit, where $p_{l+1}$ is the first periodic point on the forward orbit of $0$. Denote the accesses to $H$ at $p_m$ by $A_m^{(1)},\ldots, A_m^{(d_m)}$, where $d_m\in\N$ is the number of branches of $H$ at $p_m$ and the accesses $A_m^{(i)}$ are numbered according to their cyclic order at $p_m$. For every access $A_m^{(i)}$, choose a leg $\gamma_m^{(i)}\colon [0,1]\to\CC$ approaching $p_m$ through $A_m^{(i)}$. Note that if $\Gamma_m^{(i)}$ is another leg approaching $p_m$ through $A_m^{(i)}$, then $\gamma_m^{(i)}$ is homotopic to $\Gamma_m^{(i)}$ rel $P(E_\lambda)$ (which is easy to see using the theory of prime ends).
We first consider the periodic part of the forward orbit of $0$. For $m\in\{l+1,\ldots,n\}$, the number of branches of $H$ at $p_m$ is independent of $m$ (say equal to $d\in\N$) because the induced self-map $f\colon H\to H$ is a local homeomorphism at $p_m$ and $p_m$ is periodic. Choose a neighborhood $V$ of $0$ such that none of the legs $\gamma_m^{(i)}$ for $m\in\{l+1,\ldots, n\}$ intersect $V$. We set $U:=\widehat{E_\lambda}^{-1}(V)$ and apply Lemma~\ref{lem:ExtendedIsotopy} to obtain an embedded tree $H''\subset\complex$ homotopic to $H'$ rel $\CT\setminus U$. The lifts $\tilde{\mathcal{L}}(\gamma_m^{(i)})$ of the $\gamma_m^{(i)}$ along the periodic post-singular orbit do not intersect $U$ and therefore are disjoint from $H''$ (as $H'$ and $H''$ are equal on $\complex\setminus U$) except for their endpoints. Now we see why we need Lemma~\ref{lem:ExtendedIsotopy}: there exists an ambient isotopy $I\colon \complex\times[0,1]\to\complex$ between $H''$ and $H$ rel $P(E_\lambda)$, and it extends to an ambient isotopy $I\colon\CC\times[0,1]\to\CC$ fixing $\infty$. Hence, $\Gamma_m^{(i)}:=\left. I\right|_{\CC\times\{1\}}(\tilde{\mathcal{L}}(\gamma_m^{(i)}))$ approaches $p_{m-1}$ (or $p_n$ for $m=l+1$) through an access $A_{m-1}^{(\sigma_m(i))}$. By the above, $\Gamma_m^{(i)}$ is homotopic to $\gamma_{m-1}^{(\sigma_m(i))}$, so we have
\[
\mathcal{L}([\gamma_m^{(i)}])=[\gamma_{m-1}^{(\sigma_m(i))}]~~~\text{for all}~m\in\{l+1,\ldots,n\}~\text{and all}~i\in\{1,\ldots,d\}.
\]
As $E_\lambda$ is locally an orientation-preserving homeomorphism and as ambient isotopies preserve the cyclic order of curves landing at a common point, we see that each $\sigma_m\colon \{1,\ldots,d\}\to\{1,\ldots,d\}$ is a power of the $d$-cycle $(1~2~\ldots~d)$. Therefore, there exists an $N\in\N$ such that $\mathcal{L}^{\circ N}([\gamma_m^{(i)}])=[\gamma_m^{(i)}]$ for all $m$ and $i$. By Proposition~\ref{pro:TopRays}, each $\gamma_m^{(i)}$ is homotopic to a dynamic ray $g_m^{(i)}$ rel $P(E_\lambda)\cup\{\infty\}$.

Since the $\gamma_m^{(i)}$ as well as the $g_m^{(i)}$ are pairwise non-homotopic rel $P(E_\lambda)\cup\{\infty\}$ and intersect each other at most at their endpoints, the two curve systems fulfill the hypotheses of Lemma~\ref{lem:CurveSystem}. Thus, there exists an ambient isotopy $I\colon \CC\times[0,1]\to\CC$ rel $P(E_\lambda)\cup\{\infty\}$ satisfying $I^1(\gamma_m^{(i)})=g_m^{(i)}$ for all $m\in\{l+1,\ldots,n\}$ and all $i\in\{1,\ldots,d\}$.
By Lemma~\ref{lem:Transformation}, the tree $H_{0}:=I^1(H)$ is again a Homotopy Hubbard Tree and it is equivalent to $H$.

The last part of the proof works by induction. Assume that $H_j\cong H$ is equivalent to our initial Homotopy Hubbard Tree and for each post-singular point $p\in P(E_\lambda)$ for which $E_\lambda^{\circ j}(p)$ is periodic and every access $A$ to $p$ in $H_j$, there exists a dynamic ray $g$ approaching $p$ through $A$.
Choose a small neighborhood $V_j$ of $0$ such that $V_j$ does not intersect the dynamic rays landing at the points $p_{l+1-j},\ldots,p_n$. The preimage $U_l:=\widehat{E_\lambda}^{-1}(V_l)$ does not intersect any of the dynamic rays landing at the points $p_{l-j},\ldots,p_n$. By Lemma~\ref{lem:ExtendedIsotopy}, there exists a tree $H_{j}''\subset\complex$ which is isotopic to $H_{j}'$ rel $\CT\setminus U_j$ and isotopic to $H_{j}$ rel $P(E_\lambda)$ in $\complex$.
By Lemma~\ref{lem:Transformation}, the tree $H_{j+1}:=H_{j}''$ is itself a Homotopy Hubbard Tree, and we have $H_{j+1}\cong H_{j}\cong H$. By construction, for every $p\in\{p_{l-j},\ldots,p_n\}$ and every access to $H_{j+1}$ at $p$, there exists a dynamic ray approaching $p$ through this access. Iterating this procedure $l$ times yields a tree $\tilde{H}:=H_l$ satisfying the conditions of the theorem.
\end{proof}

\begin{Cor}[Equal kneading sequences]\label{cor:TreeRespectsBoundaries}
Let $H$ be a Homotopy Hubbard Tree for the post-singularly finite exponential map $E_\lambda$, let $f\colon H\to H$ be an induced self-map of $H$, and let $v_T$ be the singular point of $f$. Two post-singular points $p,q\in P(E_\lambda)$ are contained in the same branch of $H$ at $v_T$ if and only if the first entries of the itineraries of $p$ and $q$ w.r.t.\ some, and hence any, dynamical partition (compare Lemma~\ref{lem:DifferentPartitions}) are equal.

As this property is independent of the choice of $H$, different Homotopy Hubbard Trees yield exponential dynamical trees with equal kneading sequences.
\end{Cor}
\begin{proof}
If $\tilde{H}\cong H$ and $\tilde{f}$ is an induced self-map of $\tilde{H}$, then there exists a homeomorphism $\theta\colon H\to\tilde{H}$ that restricts to a conjugation between $f$ and $\tilde{f}$ on the set $V_f$ of marked points. This follows from the elaborations on different self-maps of the same tree after Definition~\ref{def:RelativeHomotopies} together with the proof of Lemma~\ref{lem:Transformation}. As $\theta$ maps branches of $H$ at $v_T$ to branches of $\tilde{H}$ at $\tilde{v_T}$, it is enough to prove the Corollary~\ref{cor:TreeRespectsBoundaries} for an equivalent Homotopy Hubbard Tree $\tilde{H}$.

By Proposition~\ref{pro:AccessesContainRays}, there exist a dynamic ray $g$ landing at $0$ and a Homotopy Hubbard Tree $\tilde{H}\cong H$ such that $\Tr(g)\cap\tilde{H}=\emptyset$. Let $\mathcal{D}$ be the dynamical partition w.r.t.\ $g$. Then we have $\tilde{H}'\cap\partial\mathcal{D}=\emptyset$, where $\tilde{H}'$ is the subtree of the preimage tree of $\tilde{H}$ spanned by $P(E_\lambda)$. It follows that post-singular points from different sectors can't be contained in the same branch of  $\tilde{H}'$ at $-\infty$. Furthermore each partition sector contains at most one branch of  $\tilde{H}'$ at $-\infty$, as the singular value is an endpoint of $\tilde{H}$ by Lemma~\ref{lem:BasicProperties}. Hence, $p$ and $q$ are contained in the same branch of $\tilde{H}'$ at $-\infty$ if and only if they are contained in the same sector of $\mathcal{D}$, i.e., if their itineraries start with the same integer. Branches of $\tilde{H}'$ at $-\infty$ get identified with branches of $H$ at $\tilde{v}_T$, so the corollary follows.
\end{proof}

We will now show that an exponential dynamical tree is already determined (up to equivalence) by its kneading sequence. The following ideas, leading to the proof of this fact, are inspired by and in many parts analogous to results of \cite{BKS}. Let us start with a simple observation.

\begin{Pro}[Distinct itineraries]\label{pro:DistinctItineraries}
Let $(H,f,B_{k_i})$ be an exponential dynamical tree, and let $p,q\in V_f$ be distinct marked points. Then $\It(p)\neq\It(q)$.
\end{Pro}

\begin{proof}
We write $\It(p)=\mathtt{p_1 p_2\ldots}$ and $\It(q)=\mathtt{q_1 q_2\ldots}$. If $\mathtt{p_1}\neq\mathtt{q_1}$, we are done. Else, we have $\mathtt{p_1}=\mathtt{q_1}=k_1$, and hence $p,q\in B_{k_1}$ for the branch $B_{k_1}$ of $H$ at $v_T$. By Lemma~\ref{lem:BasicProperties}, the restriction $\restr{f}{B_{k_1}}$ is injective, so we have $f([p,q])=[f(p),f(q)]$. By expansivity of $(H,f)$, there exists a smallest $n\geq 0$, such that $v_T\in f^{\circ n}([p,q])=[f^{\circ n}(p),f^{\circ n}(q)]$. Therefore, the itineraries $\It(f^{\circ n}(p))$ and $\It(f^{\circ n}(q))$ have different initial entries.
\end{proof}

Let $(H,f,B_{k_i})$ be an exponential dynamical tree and let $p,q,r\in V_f$ be distinct marked points. The subtree $[p,q,r]$ spanned by $p$, $q$ and $r$ is called a \emph{triod}. It must look in one of two ways:
if $[p,q,r]$ is homeomorphic to the letter {\sffamily Y}, we call it \emph{branched}. Else, it is homeomorphic to the letter {\sffamily I} and we call it \emph{linear}. It should now be clear how to define the \emph{middle point} $b[p,q,r]\in H$ \emph{of the triod $[p,q,r]$}. Note that the middle point $b[p,q,r]\in V_f$ is also a marked point.

The shape of the triod $[p,q,r]$ as well as the itinerary of the middle point $b[p,q,r]$ can be determined algorithmically from the itineraries of the marked points $p$, $q$, and $r$. This algorithm is known as the \emph{triod algorithm}.
The triod algorithm is purely combinatorical and could be applied to any triple of preperiodic sequences (not just the itineraries of $p$, $q$ and $r$). We apply it only to itineraries of vertices of an abstract exponential tree $(H, f, B_{k_i})$; let $\nu$ be its kneading sequence.
\begin{Def}[Formal triods and the formal triod map]~\\
Let $\star$ be a formal symbol not contained in $\Z$. The \emph{space of formal (pre-)periodic points} $\ITN\subset (\Z\cup\{\star\})^\N$ consists of all (pre-)periodic sequences $\itn{t}\in(\Z\cup\{\star\})^\N$ such that either $\itn{t}\in\Z^\N$ and $\sigma^{\circ n}(\mathtt{\itn{t}})\neq\nu$ for all $n\in\N$ or $\itn{t}$ is contained in the backwards orbit
\[
\Omega^-(\star\nu):= \{k_1\ldots k_n\star\nu~\vert ~n\in\N_0, k_i\in\Z \}\subset (\Z\cup\{\star\})^\N
\]
of finite integer sequences followed by $\star\nu$. We call a point $\itn{t}\in\ITN$ a \emph{pre-singular} point if $\itn{t}\in\Omega^-(\star\nu)$.

We have the \emph{shift map} $\sigma\colon (\Z\cup\{\star\})^\N\to(\Z\cup\{\star\})^\N$ acting on the full space of sequences. Note that we are using the same symbol for the shift-map on $\mathcal{S}$. It will always be clear from the context which map we are considering.

Any triple of distinct sequences $\itn{t},\itn{u},\itn{v}\in\ITN$ is called a \emph{formal triod} $\mathtt{[\itn{t},\itn{u},\itn{v}]}$. Given a formal triod we define the \emph{formal triod map} $\Triod$ as follows:
\[
\Triod[\itn{t},\itn{u},\itn{v}]:=
\begin{cases}
[\sigma(\itn{t}), \sigma(\itn{u}), \sigma(\itn{v})]&\text{if}~\mathtt{t}_1=\mathtt{u}_1=\mathtt{v}_1,\\
[\sigma(\itn{t}), \sigma(\itn{u}), \nu]&\text{if}~\mathtt{t}_1=\mathtt{u}_1\neq \mathtt{v}_1,\\
[\sigma(\itn{t}), \nu, \sigma(\itn{v})]&\text{if}~\mathtt{t}_1=\mathtt{v}_1\neq \mathtt{u}_1,\\
[\nu, \sigma(\itn{u}), \sigma(\itn{v})]&\text{if}~\mathtt{t}_1\neq \mathtt{u}_1=\mathtt{v}_1,\\
\textbf{stop}&\text{if}~ \mathtt{t}_1, \mathtt{u}_1, \mathtt{v}_1~\text{distinct}.
\end{cases}
\]\label{def:FormalTriods}
\end{Def}
In all cases other than the \textbf{stop} case, $\Triod[\itn{t},\itn{u},\itn{v}]$ is again a formal triod: since $\ITN$ is forward invariant under $\sigma$ and $\nu\in\ITN$, the three image sequences are contained in $\ITN$, and they are distinct, as we have $\{\itn{t},\itn{u},\itn{v}\}\cap\{k\nu~\vert~k\in\Z\}=\emptyset$. By construction, the only sequence that starts with $\star$ is $\star\nu$, so in all cases other than the \textbf{stop} case at least two of the first entries of the involved sequences are equal integers. If exactly two of the three first entries are equal, we say that the sequence whose first entry differs from the other two gets \emph{chopped off} under iteration of $\Triod$.\\
The formal triod map can be iterated as long as the \textbf{stop} case is not reached. We write $\Triod^{\circ n}[\itn{t},\itn{u},\itn{v}]$ for the resulting formal triod after $n$ iterations of $\Triod$ (if iteration is possible). Note that, if the triod can be iterated indefinitely, at least two sequences each get chopped off under iteration of $\Triod$ infinitely many times: otherwise, there would exist an $n\in\N$ such that two of the three sequences of the triod $\Triod^{\circ n}[\itn{t},\itn{u},\itn{v}]$ never get chopped off. However, this implies that these sequences are equal, contradicting the fact that all three sequences stay distinct under iteration of $\Triod$.

\begin{Def}[Majority vote and middle point of a formal triod]~\\
Let $[\itn{t},\itn{u},\itn{v}]$ be a formal triod.
If $\Triod[\itn{t},\itn{u},\itn{v}]\neq\textbf{stop}$, then as noted above, at least two of the three sequences start with the same integer. We denote this integer by $\mathtt{m}(\itn{t},\itn{u},\itn{v})$ and call it the \emph{majority vote} of the triod $[\itn{t},\itn{u},\itn{v}]$.
Let $i_0\in\N_0$ be chosen such that $\Triod^{\circ i_0}[\itn{t},\itn{u},\itn{v}]=\textbf{stop}$, if the triod eventually reaches the stop case, and let $i_0=\infty$ otherwise. We define a sequence $\itn{b}(\itn{t},\itn{u},\itn{v})\in\Z^\N\cup\Omega^-(\star\nu)$ by setting
\[
(\itn{b}(\itn{t},\itn{u},\itn{v}))_i:=
\begin{cases}
\mathtt{m}(\Triod^{\circ (i-1)}[\itn{t},\itn{u},\itn{v}])&\text{if}~i<i_0,\\
\star&\text{if}~i=i_0,\\
\nu_{i-i_0}&\text{if}~i>i_0,
\end{cases}
\]
and we call $\itn{b}(\itn{t},\itn{u},\itn{v})$ the \emph{middle point} of the triod.
If we have $\itn{b}(\itn{t},\itn{u},\itn{v})\notin\{\itn{t},\itn{u},\itn{v}\}$, we call $[\itn{t},\itn{u},\itn{v}]$ \emph{branched}, and otherwise we call it \emph{linear}. Sometimes, we want to be more precise and call a triod \emph{pre-singularly branched} or \emph{pre-singularly linear} if $i_0<\infty$, i.e., if it eventually reaches the \textbf{stop} case under iteration of the triod map.\label{def:BranchPoints}
\end{Def}

We now prove that the triod algorithm really determines the itineraries of the branch points of an exponential dynamical tree. The proof is an adaption of \cite[Proposition~3.5]{BKS}.

\begin{Lem}[Correctness of the triod algorithm]\label{lem:Correctness}
Let $(H,f,B_{k_i})$ be an exponential dynamical tree and let $[p_1,p_2,p_3]\subset H$ be a triod. Then $[p_1,p_2,p_3]$ is branched if and only if the formal triod $[\It(p_1),\It(p_2),\It(p_3)]$ is branched. Furthermore, we have
\[
\It(b[p_1,p_2,p_3])=\itn{b}(\It(p_1),\It(p_2),\It(p_3)).
\]
\end{Lem}
\begin{proof}
If $p_1$, $p_2$, and $p_3$ are contained in distinct branches of $H$ at $v_T$, then $b[p_1,p_2,p_3]=v_T$, and this is the branch point determined (on the level of itineraries) by the triod algorithm.  If one of the $p_i$ equals $v_T$ and the other two points are contained in different branches of $H$, then $b[p_1,p_2,p_3]=p_i$, and again this is the output of the triod algorithm.

If both of these cases do not occur, there exist a branch $B_{k_1}$ of $H$ at $v_T$ and distinct indices $i,j\in\{1,2,3\}$, such that $p_i,p_j\in B_{k_1}\setminus\{v_t\}$. Hence, we also have $b[p_1,p_2,p_3]\in B_{k_1}\setminus\{v_T\}$. We see, that the first entry of $\It(b[p_1,p_2,p_3])$ is calculated correctly by the triod algorithm.

If all of the $p_l$ are contained in $B_{k_1}$, the spanned subtree $[p_1,p_2,p_3]$ is also entirely contained in $B_{k_1}$. By Lemma~\ref{lem:BasicProperties}, the restriction $\restr{f}{[p_1,p_2,p_3]}$ is injective, so $[f(p_1),f(p_2),f(p_3)]$ is also a triod and we have $b[f(p_1),f(p_2),f(p_3)]=f(b[p_1,p_2,p_3])$. On the level of formal triods, we have $\Triod[\It(p_1),\It(p_2),\It(p_3)]=[\It(f(p_1)),\It(f(p_2)),\It(f(p_3))]$.
If instead $p_l\notin B_{k_1}$, while $p_i,p_j\in B_{k_1}$, then the chopped off triod $[v_T, p_i, p_j]$ still gets mapped forward injectively. Hence, we have $b[0,f(p_i),f(p_j)]=f(b[p_1,p_2,p_3])$. On the level of formal triods, we have $\Triod[\It(p_1),\It(p_2),\It(p_3)]=[\It(0),\It(f(p_i)),\It(f(p_j))]$.

If we are in one of the two preceding cases (we haven't reached the \textbf{stop} case in the first iteration step), we apply the same reasoning as before to the image triod. Hence, the triod algorithm correctly calculates $\It(b[p_1,p_2,p_3])$. It remains to show, that $[p_1,p_2,p_3]$ is branched if and only if the formal triod $[\It(p_1),\It(p_2),\It(p_3)]$ is branched. If $[p_1,p_2,p_3]$ is branched, then
\[
\It(b[p_1,p_2,p_3])=\itn{b}(\It(p_1),\It(p_2),\It(p_3))\neq\It(p_i)~~~\text{for all}~i\in\{1,2,3\}
\]
by Proposition~\ref{pro:DistinctItineraries}, so the formal triod $[\It(p_1),\It(p_2),\It(p_3)]$ is also branched. If $[p_1,p_2,p_3]$ is linear, then
\[
\It(b[p_1,p_2,p_3])=\itn{b}(\It(p_1),\It(p_2),\It(p_3))=\It(p_i)~~~\text{for some}~i\in\{1,2,3\},
\]
so the formal triod $[\It(p_1),\It(p_2),\It(p_3)]$ is also linear.
\end{proof}

\begin{proof}[Proof of Theorem~\ref{thm:UniqueGraphStructure}]
  It remains to show that two exponential dynamical trees $(H,f,B_{k_i})$ and $(\tilde{H},\tilde{f},\tilde{B}_{k_i})$ with the same kneading sequence are equivalent. Since every endpoint of $H$ is a post-singular point, every branch point of $H$ is also the branch point of a post-singular triod. Thus, the kneading sequence fully determines the itineraries of all marked points of $H$ and for every triod $[p,q,r]$ of marked points it determines their incidence relation by Lemma~\ref{lem:Correctness}. As $(H,f,B_{k_i})$ and $(\tilde{H},\tilde{f},\tilde{B}_{k_i})$ have the same kneading sequence, their marked points have the same itineraries, and we can define a homeomorphism $\varphi\colon H\to \tilde{H}$ by sending the marked points of $(H,f,B_{k_i})$ to the marked points of $(\tilde{H},\tilde{f},\tilde{B}_{k_i})$ with the same itinerary and extending this map to the edges between the marked points. The exponential dynamical trees $(H,f,B_{k_i})$ and $(\tilde{H},\tilde{f},\tilde{B}_{k_i})$ are equivalent via $\varphi$.
\end{proof}

Theorem~\ref{thm:UniqueGraphStructure} is an important step in proving uniqueness of Homotopy Hubbard Trees, but we still have to see that the embedding of the tree into the plane is also unique. Let $H$ be a Homotopy Hubbard Tree and let $[p,q,r]$ be a triod of post-singular points. If $[p,q,r]$ is linear with middle point $p$, Proposition~\ref{pro:AccessesContainRays} implies that there are two dynamic rays $g_{\ad{p}},g_{\ad{p}'}$ landing at $p$ and separating $q$ from $r$ in the sense that $q$ and $r$ are contained in different connected components of $\complex\setminus(\Tr(g_{\ad{p}})\cup\Tr(g_{\ad{p}'})\cup\{p\}$. Conversely, if such separating rays exist, the triod $[p,q,r]$ is linear with middle point $p$:
\begin{Lem}[Separating rays determine triod type]
Let $H$ be a Homotopy Hubbard Tree and let $[p,q,r]$ be a triod of post-singular points. If there are two dynamic rays $g_{\ad{p}},g_{\ad{p}'}$ landing at $p$, such that $q$ and $r$ are contained in different connected components of $\complex\setminus(\Tr(g_{\ad{p}})\cup\Tr(g_{\ad{p}'})\cup\{p\})$, the triod $[p,q,r]$ is linear with middle point $p$.\label{lem:TriodType}
\end{Lem}
\begin{proof}
Let $\It(p)=\itn{p}=\mathtt{p_1 p_2\ldots}$, $\It(q)=\itn{q}=\mathtt{q_1 q_2\ldots}$, and $\It(r)=\itn{r}=\mathtt{r_1 r_2\ldots}$ be the itineraries of the points $p$, $q$ and $r$. By correctness of the triod algorithm (see Lemma~\ref{lem:Correctness}), it is enough to show that $\itn{b}(\itn{p},\itn{q},\itn{r})=\itn{p}$.

By Proposition~\ref{pro:LandingBehavior} and Theorem~\ref{thm:LandingTheorem}, there exist addresses $\ad{q}$ and $\ad{r}$ satisfying $\It(\ad{q}~\vert~\ad{s})=\itn{q}$ and $\It(\ad{r}~\vert~\ad{s})=\itn{r}$. The triod $[\ad{p},\ad{q},\ad{r}]$ is called a triod of external addresses associated to $[\itn{p},\itn{q},\itn{r}]$ (this terminology is introduced rigorously in Section~\ref{sec:separatingRays}). It will be shown in Section~\ref{sec:separatingRays} that $[\itn{p},\itn{q},\itn{r}]$ has the same shape as the triod $[\ad{p},\ad{q},\ad{r}]$ of associated external addresses.
Finally, Lemma~\ref{lem:CombSeparatingRays} implies that $\itn{b}(\itn{p},\itn{q},\itn{r})=\itn{p}$.
\end{proof}

We are now in the position to prove uniqueness of Homotopy Hubbard Trees.
\begin{Thm}[Uniqueness of exponential Hubbard Trees]~\\
	Let $H$ and $\tilde{H}$ be Homotopy Hubbard Trees for the post-singularly finite exponential map $E_\lambda$. Then $H$ and $\tilde{H}$ are homotopic relative to the post-singular set $P(E_\lambda)$.
\end{Thm}
\begin{proof}
	Let $P(E_\lambda)=\{p_0=0, p_1=E_\lambda(p_0), \ldots, p_n=E_\lambda^{\circ n}(0), p_{n+1}=p_{l+1}\}$ denote the post-singular orbit of $E_\lambda$.
	At every post-singular point $p_i$, the trees $H$ and $\tilde{H}$ have the same number of branches by Theorem~\ref{thm:UniqueGraphStructure}; denote this number by $d_i$.
	By Proposition~\ref{pro:AccessesContainRays}, we can assume w.l.o.g.\ that there are dynamic rays $g_1^{(i)},\ldots, g_{d_i}^{(i)}$ landing at $p_i$ indexed according to their cyclic order at $p_i$ and satisfying $H\cap\Tr(g_l^{(i)})=\emptyset$ such that $[p_i,p_j,p_k]_H$ is linear with middle point $p_i$ if and only if $p_j$ and $p_k$ are contained in different connected components of $\complex\setminus((\cup_l\Tr(g_l^{(i)}))\cup\{p_i\})$. In the same way, we choose rays $\gamma_l^{(i)}$ for $\tilde{H}$, and we claim that (possibly after re-indexing) $g_l^{(i)}$ is homotopic to $\gamma_l^{(i)}$ rel $P(E_\lambda)$.

	Otherwise, there are indices $i_0$ and $l_0$ such that $\gamma_{l_0}^{(i_0)}$ is not homotopic to any of the rays $g_l^{(i_0)}$. There is an index $l'$ such that $ g_{l'}^{(i_0)}\prec \gamma_{l_0}^{(i_0)}\prec g_{l'+1}^{(i_0)}$ holds. Furthermore, there are rays $g_{l_j}^{(j)}$ and $g_{l_k}^{(k)}$ landing at distinct post-singular points $p_j$ and $p_k$ different from $p_{i_0}$ such that $ g_{l'}^{(i_0)}\prec g_{l_j}^{(j)}\prec\gamma_{l_0}^{(i_0)}$ and $ \gamma_{l_0}^{(i_0)}\prec g_{l_k}^{(k)}\prec g_{l'+1}^{(i_0)}$ hold: otherwise,  $\gamma_{l_0}^{(i_0)}$ would either be homotopic to $g_{l'}^{(i_0)}$ or to $g_{l'+1}^{(i_0)}$. But then $p_j$ and $p_k$ are separated by $\gamma_{l_0}^{(i_0)}$ and $g_{l'}^{(i_0)}$, so $[p_{i_0},p_j,p_k]$ is linear with middle point $p_{i_0}$ by Lemma~\ref{lem:TriodType}, yielding a contradiction.

	By Proposition~\ref{pro:TopRays}, homotopic dynamic rays landing together at a \emph{periodic} post-singular point are equal, so we actually have $\gamma_l^{(i)}=g_l^{(i)}$ for all $i\in\{l+1,\ldots, n\}$. The argument in the preceding paragraph also shows that the number of rays landing at $p_i$ equals $d_i$.
  Analyzing the proof of Proposition~\ref{pro:AccessesContainRays}, we conclude that we can choose $H$ and $\tilde{H}$ as to not intersect any dynamic ray landing at any post-singular point (periodic or not). In particular, we can find a ray $g$ landing at $0$ such that $H\cap E_\lambda^{\circ i}(\Tr(g))=\tilde{H}\cap E_\lambda^{\circ i}(\Tr(g))=\emptyset$ for all $i\in\{0,\ldots, n\}$.
  Let $\varphi\colon \D\to\complex\setminus\cup_i(E_\lambda^{\circ i}(\Tr(g))\cup\{p_i\})$ be a conformal map.
  By Carath\'eodory's Theorem, $\varphi$ extends continuously to $\partial\D$ and by Lemma~\ref{lem:ClosedDiskHomotopy} the trees $\varphi^{-1}(H)$ and $\varphi^{-1}(\tilde{H})$ are homotopic rel $\partial\D$.
  Pushing this homotopy forward via $\varphi$, we see that $H$ and $\tilde{H}$ are homotopic rel $P(E_\lambda)$, hence they are equivalent.
\end{proof}

\section{Separating dynamic rays: embedding the tree into the plane}
\label{sec:separatingRays}

In this section, we show the existence of Homotopy Hubbard Trees via an explicit construction. We use the triod algorithm from the preceding section to determine the middle points of post-singular triods, and thereby the set of marked points of the Homotopy Hubbard Tree, on the level of itineraries. Using Proposition~\ref{pro:LandingBehavior}, we pick (pre-)periodic points in the complex plane realizing the itineraries in the formal set of marked points (for pre-singular itineraries, the construction is a bit more involved because there is no actual (pre-)periodic point in the plane realizing this itinerary).

It remains to find the right way to embed the edges of the tree into the complex plane. Our main idea is to find an embedded tree which does not intersect any dynamic rays landing at a marked point. (This statement is only approximately true; some rays landing at pre-singular marked points need to be intersected, but we do so in a controlled way.)

This is partially motivated by Proposition~\ref{pro:AccessesContainRays}: if there exists a Homotopy Hubbard Tree, then up to homotopy rel $P(E_\lambda)$, every access to the tree at every post-singular point contains a dynamic ray. In particular, this ray does not intersect the Homotopy Hubbard Tree. Another reason is the analogy to polynomials. A polynomial Hubbard Tree does not intersect any dynamic rays. If a polynomial does not have bounded Fatou components, every access to every point on its Hubbard Tree contains a dynamic ray.

The key observation is that the union of the rays landing at marked points partitions the plane in a meaningful way: there exists an embedded tree spanned by $P(E_\lambda)$ that does not intersect these rays (except for some rays landing at pre-singular points, as noted above), and this tree is unique up to homotopy rel $P(E_\lambda)$. As the rays landing at marked points form a forward invariant set, the preimage tree also does not intersect them, so the embedded tree is invariant up to homotopy. Expansivity of the induced self-map follows because different marked points have different itineraries, so the embedded tree is a Homotopy Hubbard Tree for $E_\lambda$.

Let us begin with the construction of Homotopy Hubbard Trees. At first, we determine the set of marked points on the level of itineraries, using the triod algorithm. Throughout this section, $\nu$ denotes the kneading sequence of the exponential map $E_\lambda$ (recall Lemma~\ref{lem:DifferentPartitions}).
Let $\Omega^+(\star\nu):=\{\sigma^{\circ k}(\star\nu)~\colon~k\geq 0\}$
denote the forward orbit of $\star\nu$ under the shift map $\sigma$. This is a finite set because $\nu$ is preperiodic. We set
\[
V_\nu:=\Omega^+(\star\nu)\cup\bigcup_{[\itn{t},\itn{u},\itn{v}]}\{\itn{b}(\itn{t},\itn{u},\itn{v})\},
\]
where the union runs over all formal triods consisting of sequences in $\Omega^+(\star\nu)$: we are adding all branch points of triods formed by post-singular points and $\star\nu$ to the set of post-singular points. We call $V_\nu$ the \emph{formal vertex set}. The following properties of $V_\nu$ will be proved later in this section.

\begin{Lem}[The formal vertex set]
	The formal vertex set $V_\nu$ has the following properties:
	\begin{enumerate}
		\item It is forward invariant under the shift map, i.e., $\sigma(V_\nu)\subset V_\nu$.
		\item It consists entirely of formal (pre-)periodic points, i.e., $V_\nu\subset\ITN$.
		\item It is closed under taking triods: for every triod $[\itn{t},\itn{u},\itn{v}]$ of formal vertices $\itn{t},\itn{u},\itn{v}\in V_\nu$ we have $\itn{b}(\itn{t},\itn{u},\itn{v})\in V_\nu$.
	\end{enumerate}\label{lem:VertexSet}
\end{Lem}

Next, we embed the formal vertex set into the complex plane. The resulting set will become the set of marked points of the yet to be constructed Homotopy Hubbard Tree. First, assume that $\itn{t}\in V_\nu$ is not a pre-singular point, i.e., the itinerary $\itn{t}\in\Z^\N$ is (pre-)periodic and $\sigma^{\circ n}(\itn{t})\neq\nu$ for all $n\in\N$. By Proposition~\ref{pro:LandingBehavior}, there is exactly one (pre-)periodic point $v_\itn{t}\in\complex$ with $\It(v_\itn{t}~\vert~\ad{s})=\itn{t}$. This is the point we assign to our formal vertex $\itn{t}$.

For a pre-singular itinerary $\itn{t}\in V_\nu$, there is no (pre-)periodic point in $\complex$ of itinerary $\itn{t}$. This issue could be addressed by adding further points at infinity corresponding to iterated preimages of $-\infty$ to the plane. For brevity, we take a more hands-on approach: we choose surrogate points in the plane which are sufficiently close to these iterated preimages. In our terms, this means that the itinerary of the vertex $v_{\itn{t}}\in\complex$ shares sufficiently many entries with $\itn{t}$. To make this precise, let
\[
N:=\max\{n\in\N_0~\vert~\exists~\itn{t}\in V_\nu\colon\itn{t}=k_1\ldots k_n\star\nu\}
\]
and pick a closed neighborhood $U$ of the singular value $0$ with the following properties:
\begin{enumerate}
	\item $E_\lambda^{-n}(U)\cap E_\lambda^{-m}(U)=\emptyset$ for all distinct $n,m\in\{0,\ldots, N+1\}.$
	\item For every $n\in\{1,\ldots, N+1\}$, the preimage $E_\lambda^{-n}(U)$ does not intersect the union of the set of non pre-singular vertices and the dynamic rays landing at these points.
	\item $U$ is bounded by a Jordan curve such that the dynamic ray $g_{\ad{s}}$ landing at the singular value intersects the boundary of $U$ exactly once, i.e., there is a unique potential $t_0$ such that $g_{\ad{s}}(t_0)\in\partial U$, whereas $g_{\ad{s}}(t)\in U$ for $t<t_0$ and $g_{\ad{s}}(t)\in\complex\setminus\overline{U}$ for $t>t_0$. This is only to simplify topological considerations involving $U$.
\end{enumerate}
The first condition is equivalent to $U\cap E_\lambda^{\circ n}(U)=\emptyset$ for all $n\in\{1,\ldots, N+1\}$ and this is obviously fulfilled for $U$ small enough because $0$ is not periodic. The second condition is equivalent to $U$ not intersecting the union of the set of non-presingular vertices and the dynamic rays landing at those except for the singular value $0$. Again, this is true for $U$ small enough. Finally, condition (3) can always be ensured by shrinking a neighborhood that fulfills (1) and (2).

Let $\log_{\ad{s},k}\colon \complex\setminus (\Tr(g_{\ad{s}})\cup\{0\})\to D_k$ be the branch of the inverse of $E_\lambda$ with the stated domain and co-domain. For every itinerary $\itn{t}\in\mathcal{S}_\nu$ of the form $\itn{t}=k_1\ldots k_n\star\nu$ with $n\in\{0,\ldots, N\}$, we define
\[
U_{\itn{t}}:=\log_{\ad{s},k_1}\ldots \log_{\ad{s},k_n}(E_\lambda^{-1}(U)),
\]
i.e., $U_{\itn{t}}$ is the iterated preimage under $E_\lambda$ of the domain $U$ constructed above by the branches of the logarithm prescribed by the entries of $\itn{t}$. Note that by property (1) of $U$, distinct $U_{\itn{t}}$ and $U_{\itn{u}}$ are disjoint.
We assign to $\star\nu\in V_\nu$ an arbitrary point $v_{\star\nu}\in U_{\star\nu}$ which will later become the singular point (see the paragraph after Definition~\ref{def:RelativeHomotopies}) of the Homotopy Hubbard Tree. If $\itn{t}=k_1\ldots k_n\star\nu$ is pre-singular, we associate to $\itn{t}$ the point
\[
v_\itn{t}:=\log_{\ad{s},k_1}\ldots \log_{\ad{s},k_n}(v_{\star\nu})\in U_{\itn{t}}.
\]
To every formal vertex $\itn{t}\in V_\nu$ we have thus associated a point $v_{\itn{t}}\in\complex$.

\begin{Def}[Vertex set and triods of vertices]
	We define the \emph{vertex set} $V\subset\complex$ to be the set of all $v_{\itn{t}}$ for $\itn{t}\in V_\nu$. A \emph{triod of vertices} $[v_\itn{t}, v_\itn{u}, v_\itn{v}]$ is a triple of distinct vertices $v_\itn{t}, v_\itn{u}, v_\itn{v}\in V$.
\label{def:PointTriods}\end{Def}

As triods of vertices are in natural bijection to triods of formal vertices, we use the terminology introduced for formal triods for vertex triods, too. In particular, we call a triod of vertices (pre-singularly) branched or (pre-singularly) linear if the corresponding formal triod is of this type.

We now turn our attention to constructing the edges of a Homotopy Hubbard Tree. The partitioning property mentioned at the beginning of this section rests on a result about dynamic rays separating the vertex set $V$: if $p, q, r\in V$ are distinct vertices and the corresponding formal triod $[p, q, r]$ is branched, then there exist (pre-)periodic dynamic rays $g_1,g_2,g_3$ of itinerary $\itn{b}[p, q, r]$, such that $p$, $q$, and $r$ are separated by the rays $g_i$, i.e., such that these three points are contained in different connected components of $\complex\setminus\bigcup\Tr(g_i)$. If $[p,q,r]$ is linear, there are two (pre-)periodic rays landing at the middle point and separating the other two points from each other. To prove the existence of such separating dynamic rays, we use the language of external addresses. This requires a variant of the triod algorithm operating on the level of external addresses.

\begin{Def}[Formal triods and the formal triod map]~\\
  Recall our notation $\mathcal{I}=\{I_k\}_{k\in\Z}$ for the dynamical partition of the shift space $\mathcal{S}$ w.r.t.\ $\ad{s}$ (see Section~\ref{sec:background}). In particular, $\ad{s}$ is the external address of a dynamic ray $g_{\ad{s}}$ landing at the singular value $0$.\\
	A \emph{formal triod of external addresses} is a triple $[\ad{t}, \ad{u}, \ad{v}]$ of external addresses $\ad{t},\ad{u}, \ad{v}\in\mathcal{S}$ such that $\ad{t}\prec\ad{u}\prec\ad{v}$, and $\It(\ad{t}~\vert~\ad{s}),\It(\ad{u}~\vert~\ad{s}),\It(\ad{v}~\vert~\ad{s})\in\ITN$ are distinct itineraries. We define the \emph{formal triod map} $\Triod_\mathcal{S}$ \emph{on the level of external addresses} as follows:
	\[
	\Triod_\mathcal{S}[\ad{t}, \ad{u}, \ad{v}]:=
	\begin{cases}
	[\sigma(\ad{t}), \sigma(\ad{u}), \sigma(\ad{v})]&\text{if}~\ad{t},\ad{u},\ad{v}\in I_k~\text{for some}~k\in\Z,\\
	[\sigma(\ad{t}), \sigma(\ad{u}), \ad{s}]&\text{if}~\ad{t},\ad{u}\in I_k~\text{for some}~k\in\Z,~\ad{v}\notin I_k,\\
	[\sigma(\ad{t}), \ad{s}, \sigma(\ad{t})]&\text{if}~\ad{t},\ad{v}\in I_k~\text{for some}~k\in\Z,~\ad{u}\notin I_k,\\
	[\ad{s}, \sigma(\ad{u}), \sigma(\ad{v})]&\text{if}~\ad{u},\ad{v}\in I_k~\text{for some}~k\in\Z,~\ad{t}\notin I_k,\\
	\textbf{stop}&\text{otherwise}.
	\end{cases}
	\]
	\label{def:AddressTriods}
\end{Def}
We claim that, if the \textbf{stop} case is not reached, the image is again a formal triod. The following result will help us to show this fact.

\begin{Lem}[Order-preserving restrictions]
Let $I=(k\ad{s},(k+1)\ad{s})$ be a partition sector and let $I^-:=(k\ad{s},(k+1)\ad{s}]$. The restriction
\[
\restr{\sigma}{I^-}\colon I^-\to\AD
\]
is a bijection and preserves the cyclic order.
\label{lem:OrderPreserving}\end{Lem}

\begin{figure}
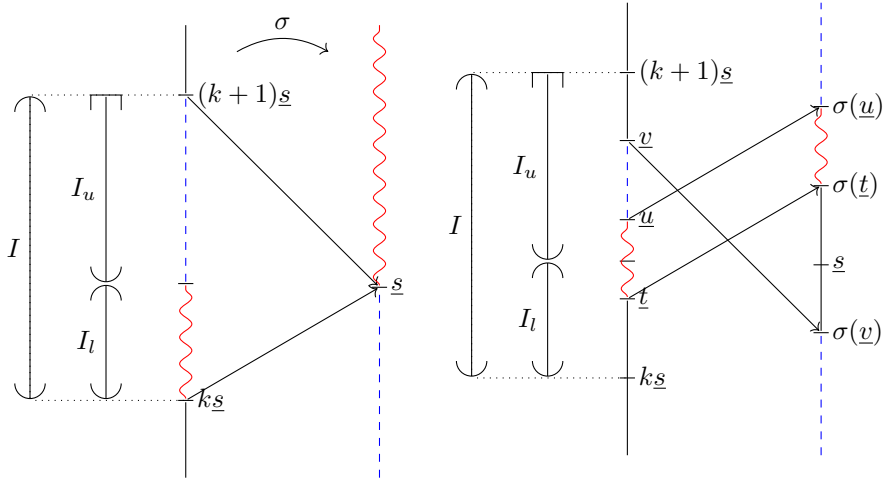

	\begin{minipage}{0.45\textwidth}
		\input{tikzdefs.tex}
\begin{tikzpicture}
 \draw pic {left_axis};
 \draw pic {partition_intervals};
 \draw pic {right_axis};
 \draw (base_left) to (ks) (k+s) to (top_left);
 \draw[interval1] (ks) -- (kinf);
 \draw[interval2] (kinf) to (k+s);
 \draw[->] (ks) to (s); \draw[->] (k+s) to (s);
 \draw[interval2] (base_right) to (s);
 \draw[interval1] (s) to (top_right);

 \node (sigmaArrowLeft)  [below right=0.3 and 0.4 of top_left] {};
 \node (sigmaArrowRight) [below left=0.3 and 0.4 of top_right] {};
 \draw [->] (sigmaArrowLeft) to [bend left] node[above] {$\sigma$} (sigmaArrowRight);
\end{tikzpicture}
	\end{minipage}
	\begin{minipage}{0.45\textwidth}
		\input{tikzdefs.tex}
\begin{tikzpicture}
 \draw pic {left_axis};
 \draw pic {partition_intervals};
 \draw pic {right_axis};

 \node[empty] (t)   [above=of ks,label=right:$\ad{t}$]       {};
 \node[empty] (u)   [above=of t,label=right:$\ad{u}$]        {};
 \node[empty] (v)   [above=of u,label=right:$\ad{v}$]        {};
 \draw (base_left) to (t);
 \draw[interval1] (t) to (u);
 \draw[interval2] (u) to (v);
 \draw (v) to (k+s) to (top_left);

 \node[empty] (sigma_v) [above=1.6 of base_right,label=right:$\sigma(\ad{v})$] {};
 \node[empty] (sigma_t) [above=of s,label=right:$\sigma(\ad{t})$] {};
 \node[empty] (sigma_u) [above=of sigma_t,label=right:$\sigma(\ad{u})$] {};
 \draw (sigma_v) to (sigma_t);
 \draw[interval1] (sigma_t) to (sigma_u);
 \draw[interval2] (base_right) to (sigma_v) (sigma_u) to (top_right);

 \foreach \n in {t, u, v, sigma_t, sigma_u, sigma_v} {
  \draw (\n) ++(-0.1,0) -- ++(0.2,0);
 }

 \draw[->] (t) to (sigma_t);
 \draw[->] (u) to (sigma_u);
 \draw[->] (v) to (sigma_v);
\end{tikzpicture}
	\end{minipage}
	\caption{Sketch illustrating the mapping behavior of the shift map. The left sketch shows how the intervals $I_l$ and $I_u$ are exchanged globally, while the restriction of the shift to each of them is strictly monotonically increasing. The right figure shows that the linear order of the image addresses might change, but the cyclic order stays the same. It also illustrates that an interval splits if it contains the address $\ad{s}$.}
	\label{fig:preimageLemmaSketches}
\end{figure}

\begin{proof}
The shift map $\sigma$ is strictly monotonically increasing w.r.t.\ the linear order on $\AD$ on both of the sets $I_l:=\{\ad{t}\in I\colon\ad{t}=k\ldots\}$ and $I_u:=\{\ad{t}\in I\colon\ad{t}=(k+1)\ldots\}\cup\{(k+1)\ad{s}\}$. It maps $I_l$ bijectively onto $\{\ad{t}\in\AD\colon\ad{t}>\ad{s}\}$ and $I_u$ bijectively onto $\{\ad{t}\in\AD\colon\ad{t}\leq\ad{s}\}$. In particular, it swaps the two sets globally, i.e.,
$\sigma(\ad{t})>\sigma(\ad{u})~\text{for all}~\ad{t}\in I_l~\text{and all}~\ad{u}\in I_u$. See Figure~\ref{fig:preimageLemmaSketches} for a sketch of the mapping behavior of $\sigma$. It follows that $\restr{\sigma}{I^-}$ is a bijection. One can see that $\restr{\sigma}{I^-}$ preserves the cyclic order by checking all possible configurations of the three external addresses w.r.t.\ $I_l$ and $I_u$. See Figure~\ref{fig:preimageLemmaSketches} for a configuration under which the cyclic order, but not the linear order, of the addresses $\ad{t}$, $\ad{u}$, and $\ad{v}$ is preserved.
\end{proof}

Let us write $[\ad{\tilde{t}},\ad{\tilde{u}},\ad{\tilde{v}}]:=\Triod_\mathcal{S}[\ad{t}, \ad{u}, \ad{v}]$. Lemma~\ref{lem:OrderPreserving} shows that $\ad{\tilde{t}}\prec\ad{\tilde{u}}\prec\ad{\tilde{v}}$. It remains to show that $\itn{\tilde{x}}:=\It(\ad{\tilde{x}}~\vert~\ad{s})$ are distinct for $x\in\{t,u,v\}$. The map $[\ad{t},\ad{u},\ad{v}]\mapsto[\itn{t},\itn{u},\itn{v}]$ (where $\itn{x}:=\It(\ad{{x}}~\vert~\ad{s})$ as before) is a semi-conjugation between $\Triod_\AD$ and $\Triod$. The sequences $\itn{\tilde{x}}$ are distinct since $[\itn{\tilde{t}},\itn{\tilde{u}},\itn{\tilde{v}}]=\Triod[\itn{t},\itn{u},\itn{v}]$, and triods of formal (pre-)periodic points are mapped to triods of formal (pre-)periodic points under iteration of $\Triod$ by the considerations after Definition~\ref{def:FormalTriods}. We call the triods $[\ad{t}, \ad{u}, \ad{v}]$ and $[\itn{t},\itn{u},\itn{v}]$ \emph{associated to each other}.
We define the \emph{middle point} $\itn{b}[\ad{t},\ad{u},\ad{v}]$ to be the middle point of $[\itn{t},\itn{u},\itn{v}]$, and the majority vote $\itn{m}[\ad{t}, \ad{u}, \ad{v}]$ to be the majority vote of $[\itn{t},\itn{u},\itn{v}]$. We call $[\ad{v},\ad{t},\ad{u}]$ \emph{branched} if $[\itn{t},\itn{u},\itn{v}]$ is branched, and linear otherwise. If $\Triod^{\circ i_0}[\itn{t},\itn{u},\itn{v}]=\textbf{stop}$, then $[\ad{t}, \ad{u}, \ad{v}]$ reaches the \textbf{stop} case at the same iteration step and vice versa.
We have seen in Section~\ref{sec:triodAlgorithm} that if $[\itn{t},\itn{u},\itn{v}]$ can be iterated indefinitely, all of the three sequences will eventually be contained in $\ITN\setminus\Omega^-(\star\nu)$ and $[\itn{t},\itn{u},\itn{v}]$ is (pre-)periodic under iteration of $\Triod$. By Theorem~\ref{thm:LandingTheorem} and Proposition~\ref{pro:LandingBehavior}, there are only finitely many external addresses associated to a formal (pre-)periodic point $\itn{t}\in\ITN\setminus\Omega^-(\star\nu)$. Therefore, $[\ad{t}, \ad{u}, \ad{v}]$ is eventually periodic under iteration of $\Triod_\mathcal{S}$.

For the proof of the existence of separating rays, which works on the combinatorial level, the following results will be useful.
\begin{Lem}[Pullbacks of Intervals]
	Let $[\ad{t}^1, \ad{t}^2, \ad{t}^3]$ be a formal triod of external addresses which can be iterated at least once before reaching the stop case, let $[\ad{u}^1, \ad{u}^2, \ad{u}^3]:=\Triod_\mathcal{S}[\ad{t}^1, \ad{t}^2, \ad{t}^3]$ be its image triod and let $k:=\mathtt{m}(\ad{t}^1, \ad{t}^2, \ad{t}^3)$ be its majority vote. Write $J_n:=(\ad{t}^n,\ad{t}^{n+1})$ and $J_n':=(\ad{u}^n, \ad{u}^{n+1})$ for $n\in\{1,2,3\}$ (with indices labeled mod 3) for the intervals of the partition of the shift space by the initial triod and the image triod respectively. Then we have
	\[
	\sigma^{-1}(J_n')\cap I_k\subseteq J_n.
	\]
\label{lem:IntervalPullbacks}\end{Lem}
\begin{proof}
	Since not all $\ad{t}^i$ are contained in distinct partition sectors by hypothesis, we have $\ad{t}^j\notin I_k$ for at most one $j$. We define a triod $[\ad{\tilde{t}}^1, \ad{\tilde{t}}^2, \ad{\tilde{t}}^3]$ by replacing such a $\ad{t}^j$ (if any) by the unique preimage of $\ad{s}$ in $I_k^-$. Then $\Triod_\AD[\ad{t}^1, \ad{t}^2, \ad{t}^3]=\Triod_\AD[\ad{\tilde{t}}^1, \ad{\tilde{t}}^2, \ad{\tilde{t}}^3]=[\sigma(\ad{\tilde{t}}^1), \sigma(\ad{\tilde{t}}^2), \sigma(\ad{\tilde{t}}^3)]$. Furthermore, we have $\tilde{J}_n\cap I_k= J_n\cap I_k$, where $\tilde{J}_n:=(\ad{\tilde{t}}^n,\ad{\tilde{t}}^{n+1})$. Let $\ad{v}\in\sigma^{-1}(J_n')\cap I_k$ be an external address. We have $\sigma(\ad{v})\in J_n'$, i.e., $\sigma(\ad{\tilde{t}}^n)\prec\sigma(\ad{v})\prec\sigma(\ad{\tilde{t}}^{n+1})$. It follows from Lemma~\ref{lem:OrderPreserving} that $\ad{\tilde{t}}^n\prec\ad{v}\prec\ad{\tilde{t}}^{n+1}$, and therefore $\ad{v}\in\tilde{J}_n\cap I_k\subset J_n$.
\end{proof}

\begin{Lem}[Splitting of intervals]
Let $I,I'\in\mathcal{I}$ be partition sectors, and let $J\subset I$ be an interval. If $\ad{s}\notin J$, then $J':=\sigma^{-1}(J)\cap I'$ is an interval and it is of the form $J'=\{k\ad{t}\colon\ad{t}\in J\}$ for some $k\in\Z$.
\label{lem:IntervalSplittling}\end{Lem}
\begin{proof}
  We have $I'=(k\ad{s},(k+1)\ad{s})$ for some $k\in\Z$. As $J$ is entirely contained in some partition sector and $\ad{s}\notin J$, we either have $\ad{t}>\ad{s}$ for all $\ad{t}\in J$ or $\ad{t}<\ad{s}$ for all $\ad{t}\in J$. Assume that the first case is true (the second case works analogously). Then we have $k\ad{t}\in I'$ for all $\ad{t}\in J$. As $\restr{\sigma}{I'}$ is injective, we have $J'=\{k\ad{t}\colon\ad{t}\in J\}$. The right part of Figure~\ref{fig:preimageLemmaSketches} illustrates why an interval containing $\ad{s}$ splits.
\end{proof}

\begin{Lem}[Unlinked addresses]
	Let $\itn{t},\itn{u}\in\Z^\N\cup\Omega^-(\star\nu)$ be distinct (pre-)periodic itineraries, let $T:=(\ad{t}^i)_{i\in I}$ be the set of external addresses with $\It(\ad{t}^i~\vert~\ad{s})=\itn{t}$ and let $U:=(\ad{u}^j)_{j\in J}$ be the set of external addresses with $\It(\ad{u}^j~\vert~\ad{s})=\itn{u}$. Then $T$ and $U$ are \emph{unlinked} in the sense that there are no addresses $\ad{t},\ad{t}'\in T$ and $\ad{u},\ad{u}'\in U$ such that $\ad{t}\prec\ad{u}\prec\ad{t}'\prec\ad{u}'$.
\label{lem:Unlinked}\end{Lem}
\begin{proof}
	Assume to the contrary that there exist four external addresses as above. First consider the case that $\mathtt{t}_1=\star$. Then $\ad{t}$ and $\ad{t'}$ are contained in the boundary $\partial\mathcal{I}$ of the dynamical partition, so $\ad{u}$ and $\ad{u}'$ are contained in different partition sectors, contradicting the fact that they have the same itinerary. Hence, we can assume $\mathtt{t}_1,\mathtt{u}_1\neq\star$. But then, $\ad{t}\prec\ad{u}\prec\ad{t}'\prec\ad{u}'$ implies $\mathtt{t}_1=\mathtt{u}_1$. As $\restr{\sigma}{I_{\mathtt{t}_1}}$ preserves the cyclic order, we have $\sigma(\ad{t})\prec\sigma(\ad{u})\prec\sigma(\ad{t}')\prec\sigma(\ad{u}')$. Repeating this argument inductively, we obtain $\itn{u}=\itn{t}$, contradicting our assumptions.
\end{proof}

Let us now prove the main combinatorial lemma of this section. Major ideas for the proof are taken from \cite[Lemma~5.2]{SZ2}.
\begin{Lem}[Combinatorial version of separating dynamic rays]
	Let $[\ad{t}^1,\ad{t}^2,\ad{t}^3]$ be a triod of external addresses and let $\itn{b}\in\Z^\N\cup\Omega^-(\star\nu)$ be (pre-)periodic.

	The triod $[\ad{t}^1,\ad{t}^2,\ad{t}^3]$ is branched with $\itn{b}=\itn{b}[\ad{t}^1,\ad{t}^2,\ad{t}^3]$, if and only if there are three distinct (pre-)periodic external addresses $\ad{s}^1, \ad{s}^2, \ad{s}^3\in\mathcal{S}$ such that $\It(\ad{s}^i~\vert~\ad{s})=\itn{b}$ and $\ad{s}^i\in(\ad{t}^i,\ad{t}^{i+1})$ (where again indices are labeled modulo 3).

	The triod $[\ad{t}^1,\ad{t}^2,\ad{t}^3]$ is linear with $\itn{b}[\ad{t}^1,\ad{t}^2,\ad{t}^3]=\It(\ad{t}^j~\vert~\ad{s})=:\itn{b}$ if and only if there are two distinct (pre-)periodic external addresses $\ad{s}^j, \ad{s}^{j+1}\in\mathcal{S}$ with $\ad{s}^j=\ad{t}^j$ and $\ad{s}^{j+1}\in(\ad{t}^{j+1},\ad{t}^{j+2})$, such that $\It(\ad{s}^i~\vert~\ad{s})=\itn{b}$ holds.
\label{lem:CombSeparatingRays}\end{Lem}
\begin{proof}
	We start by proving the ``only if'' direction.
	\vspace{2mm}

	\textbf{Claim 1.} If $\Triod_\AD[\ad{t}^1,\ad{t}^2,\ad{t}^3]=[\ad{u}^1,\ad{u}^2,\ad{u}^3]$, and the result is true for $[\ad{u}^1,\ad{u}^2,\ad{u}^3]$, then it also holds for $[\ad{t}^1,\ad{t}^2,\ad{t}^3]$.
	\vspace{2mm}

	Assume first, that $[\ad{t}^1,\ad{t}^2,\ad{t}^3]$, and hence also $[\ad{u}^1,\ad{u}^2,\ad{u}^3]$, is branched, and let $\ad{\tilde{s}}^j$ be separating addresses for the image triod. Then, setting $k:=\mathtt{m}(\ad{t}^1, \ad{t}^2, \ad{t}^3)$ and $\ad{s}^j:=\restr{\sigma}{I_k}^{-1}(\ad{\tilde{s}}^j)$, we have $\ad{s}^j\in(\ad{t}^j,\ad{t}^{j+1})$ by Lemma~\ref{lem:IntervalPullbacks}. Furthermore, we have $\It(\ad{s}^j~\vert~\ad{s})=\itn{b}[\ad{t}^1,\ad{t}^2,\ad{t}^3]$, so the addresses $\ad{s}^j$ are separating addresses for the triod $[\ad{t}^1,\ad{t}^2,\ad{t}^3]$. The linear case works analogously, and is left to the reader.\hfill$\triangle$
	\vspace{2mm}

	\textbf{Claim 2.} The result is true if $\itn{b}[\ad{t}^1,\ad{t}^2,\ad{t}^3]\in\Omega^-(\star\nu)$.
	\vspace{2mm}

	There exists an $n\geq 0$, such that $\Triod_\AD[\ad{t}^1,\ad{t}^2,\ad{t}^3]=\textbf{stop}$. Assume first, that $[\ad{t}^1,\ad{t}^2,\ad{t}^3]$ is branched. The image triod $[\ad{u}^1,\ad{u}^2,\ad{u}^3]:=\Triod_\mathcal{S}^{\circ (n-1)}[\ad{t}^1,\ad{t}^2,\ad{t}^3]$ consists of external addresses lying in distinct sectors of the dynamical partition $\mathcal{I}$. Therefore, there exists addresses $\ad{\tilde{s}}^i:=k_i\ad{s}$ for suitable $k_i\in\Z$, such that $\ad{\tilde{s}}^i\in(\ad{u}^i,\ad{u}^{i+1})$, and we have $\It(\ad{\tilde{s}}^i~\vert~\ad{s})=\star\nu=\itn{b}[\ad{u}^1,\ad{u}^2,\ad{u}^3]$ for every $i\in\{1,2,3\}$. Hence, we have proved the lemma for the image triod $[\ad{u}^1,\ad{u}^2,\ad{u}^3]$. By Claim~1, it also holds for $[\ad{t}^1,\ad{t}^2,\ad{t}^3]$. The linear case works analogously and is left to the reader. \hfill$\triangle$
	\vspace{2mm}

	Now assume that $\itn{b}[\ad{t}^1,\ad{t}^2,\ad{t}^3]\notin \Omega^-(\star\nu)$, so the triod $[\ad{t}^1,\ad{t}^2,\ad{t}^3]$ can be iterated indefinitely under $\Triod_\AD$. Every such triod is eventually periodic, so by Claim~1, we can assume w.l.o.g.\ that $[\ad{t}^1,\ad{t}^2,\ad{t}^3]$ is periodic, say of period $N$. Hence, $\itn{b}[\ad{t}^1,\ad{t}^2,\ad{t}^3]$ is also periodic (under iteration of $\Triod$), possibly of smaller period. Let us write $\itn{b}:=\itn{b}[\ad{t}^1,\ad{t}^2,\ad{t}^3]=:\overline{\mathtt{b}_0\mathtt{b}_1\ldots \mathtt{b}_{n-1}}$, where the indices of $\itn{b}$ are labeled modulo $n$.
	\vspace{2mm}

	\textbf{Claim 3.} For every $i\in\{1,2,3\}$, there exists an address $\ad{s}^i\in\cl(\ad{t}^i,\ad{t}^{i+1})$ satisfying $\It(\ad{s}^i~\vert~\ad{s})=\itn{b}$.
	\vspace{2mm}

	Our strategy for the proof is to pull back a suitable interval along the inverse branches of $\sigma$ prescribed by $\itn{b}$, to obtain a nested sequence of intervals whose intersection is an external address with the desired properties. The forward orbit $T:=\Omega^+(\ad{s})$ is finite, so $I_{\mathtt{b}_0}\setminus T$ consists of a finite number of disjoint open intervals.
	Write $J_{1,i}^{(0)}, \ldots, J_{n_i, i}^{(0)}$ for the subintervals of $(\ad{t}^i, \ad{t}^{i+1})$ contained in $I_{\mathtt{b}_0}\setminus T$. Inductively, we define
	\[
	J_{j,i}^{(m)}:=\sigma^{-1}(J_{j,i}^{(m-1)})\cap I_{\mathtt{b}_{-m}}.
	\]
	As $\AD\setminus T$ is backward invariant, we have $J_{j.i}^{(m)}\cap T=\emptyset$, so every set $J_{j,i}^{(m)}$ is an interval by Lemma~\ref{lem:IntervalSplittling}. Writing $[\ad{t}^{1, (j)},\ad{t}^{2, (j)},\ad{t}^{3, (j)}]:=\Triod_\mathcal{S}^{\circ j}[\ad{t}^1,\ad{t}^2,\ad{t}^3]$, we have $J_{j.i}^{(m)}\subseteq (\ad{t}^{i, (-m)},\ad{t}^{i+1, (-m)})$ by Lemma~\ref{lem:IntervalPullbacks}.

	After $N$ pullbacks, we arrive at subintervals of the initial partition sector $I_{\mathtt{b}_0}$ satisfying $J_{j,i}^{(N)}\subseteq (\ad{t}^i,\ad{t}^{i+1})$. Hence, each $J_{j,i}^{(N)}$ must be a subinterval of one of our initial intervals $J_{j,i}^{(0)}$. More precisely, we get three self-maps $\rho_i\colon \{1,\ldots, n_i\}\to\{1,\ldots,n_i\}$ such that
	\[
	J_{j,i}^{(N)}\subseteq J_{\rho_i(j),i}^{(0)}~~~\text{for each}~i\in\{1,2,3\}~\text{and}~j\in\{1,\ldots, n_i\}.
	\]

	For each $i\in\{1,2,3\}$ we can find a $j_i\in\{1,\ldots, n_i\}$ and an $l_i\in\N$ such that $\rho_i^{\circ l_i}(j_i)=j_i$ and setting $N_i:=l_i\cdot N$ we have $J_{j_i,i}^{(N_i)}\subset J_{j_i,i}^{(0)}$. Furthermore, the first $N_i$ entries are the same for all $\ad{t}\in J_{j_i,i}^{(N_i)}$, say equal to $s_0^{(i)}s_1^{(i)}\ldots s_{N_i-1}^{(i)}$. Setting
	$J_i^{(m)}:=J_{j_i,i}^{(mN_i)}$, the $J_i^{(m)}$ form a nested sequence of (open) intervals and every address $\ad{t}\in J_i^{(m)}$ begins with $m$ times the sequence $s_0^{(i)}s_1^{(i)}\ldots s_{N_i-1}^{(i)}$. We define $\ad{s}^i:=\overline{s_0^{(i)}s_1^{(i)}\ldots s_{N_i-1}^{(i)}}$.
	Recall that for an address $\ad{t}\in\AD$, the cylinder sets $\{\ad{u}\in\AD\colon u_0\ldots u_i=t_0\ldots t_i\}$ form an open neighborhood basis of $\ad{t}$. As every address $\ad{t}\in J_i^{(m)}$ begins with $m$ times the sequence $s_0^{(i)}s_1^{(i)}\ldots s_{N_i-1}^{(i)}$, we see that $U\cap J_i^{(m)}\neq\emptyset$ for every open neighborhood $U$ of $\ad{s}^i$ and every $m\geq 0$. Hence, we have $\ad{s}^i\in\bigcap_{m\in\N}\text{cl}(J_i^{(m)})$. This implies $\sigma^{\circ m}(\ad{s}^i)\in\overline{I_{\mathtt{b}_{m}}}$ for all $m\in\N_0$. But we have $\sigma^{\circ m}(\ad{s}^i)\notin\partial I_{\mathtt{b}_{m}}$ because otherwise $\sigma^{\circ m}(\ad{s}^i)$ would be strictly preperiodic, contradicting the fact that $\ad{s}^i$ is periodic. Hence, we have $\It(\ad{s}^i~\vert~\ad{s})=\itn{b}$.\hfill$\triangle$
	\vspace{2mm}

	If $[\ad{t}^1,\ad{t}^2,\ad{t}^3]$ is branched, we actually have $\ad{s}^i\in (\ad{t}^i,\ad{t}^{i+1})$ because $\It(\ad{t}^i~\vert~\ad{s})\neq\itn{b}$ for all $i$.
 If instead $[\ad{t}^1,\ad{t}^2,\ad{t}^3]$ is linear with $\ad{t}^j$ in the middle,
 we have $\ad{s}^{j+1}\in(\ad{t}^{j+1},\ad{t}^{j+2})$ since $\It(\ad{s}^{j+1}~\vert~\ad{s})=\itn{b}$ is distinct from both $\It(\ad{t}^{j+1}~\vert~\ad{s})$ and $\It(\ad{t}^{j+2}~\vert~\ad{s})$.
   Hence, $\ad{t}^j$ and $\ad{s}^{j+1}$ form separating external addresses in this case. This finishes the proof of the ``only if'' direction.

	Let us now prove the other direction. Assume first, that there exists a (pre-)periodic $\itn{b}\in\Z^\N\cup\Omega^-(\star\nu)$ and three distinct (pre-)periodic external addresses $\ad{s}^1, \ad{s}^2, \ad{s}^3\in\mathcal{S}$ such that $\It(\ad{s}^i~\vert~\ad{s})=\itn{b}$ and $\ad{s}^i\in(\ad{t}^i,\ad{t}^{i+1})$. We want to see that $[\ad{t}^1,\ad{t}^2,\ad{t}^3]$ is branched with $\itn{b}[\ad{t}^1,\ad{t}^2,\ad{t}^3]=\itn{b}$. Assume to the contrary that $[\ad{t}^1,\ad{t}^2,\ad{t}^3]$ is branched with $\itn{b}[\ad{t}^1,\ad{t}^2,\ad{t}^3]\neq\itn{b}$. By the above, there are three distinct (pre-)periodic external addresses $\ad{\tilde{s}}^1, \ad{\tilde{s}}^2, \ad{\tilde{s}}^3\in\mathcal{S}$ such that $\It(\ad{\tilde{s}}^i~\vert~\ad{s})=\itn{b}[\ad{t}^1,\ad{t}^2,\ad{t}^3]$ and $\ad{\tilde{s}}^i\in(\ad{t}^i,\ad{t}^{i+1})$. This contradicts Lemma~\ref{lem:Unlinked} because the sets $\{\ad{s}^i\}$ and $\{\ad{\tilde{s}}^j\}$ would not be unlinked. The same lemma leads to a contradiction if we assume, that $[\ad{t}^1,\ad{t}^2,\ad{t}^3]$ is linear. The linear case works analogously and is left to the reader.
\end{proof}

Let us translate Lemma~\ref{lem:CombSeparatingRays} into the language of dynamic rays and their landing points. This requires Lemma~\ref{lem:VertexSet} which is why we prove it now.

\begin{proof}[Proof of Lemma~\ref{lem:VertexSet}]
	The formal vertex set $V_\nu$ is forward invariant under $\sigma$ since $\Omega^+(\star\nu)$ is forward invariant and $\sigma(\itn{b}(\itn{t},\itn{u},\itn{v}))=\itn{b}(\Triod[\itn{t},\itn{u},\itn{v}])$ as long as $\Triod[\itn{t},\itn{u},\itn{v}]\neq$ \textbf{stop}. Also, if $\Triod[\itn{t},\itn{u},\itn{v}]=$\textbf{stop}, then $\itn{b}(\itn{t},\itn{u},\itn{v})=\star\nu\in V_\nu$.

	For the proof of $(2)$, assume to the contrary that there exists a triod $[\tilde{\itn{t}},\tilde{\itn{u}},\tilde{\itn{v}}]$ with $\tilde{\itn{t}},\tilde{\itn{u}},\tilde{\itn{v}}\in\Omega^+(\star\nu)$ and $\itn{b}(\tilde{\itn{t}},\tilde{\itn{u}},\tilde{\itn{v}})\notin \ITN$. Since $\tilde{\itn{t}},\tilde{\itn{u}},\tilde{\itn{v}}$ are all (pre-)periodic, so is $\itn{b}(\tilde{\itn{t}},\tilde{\itn{u}},\tilde{\itn{v}})$, and we must have $\itn{b}(\tilde{\itn{t}},\tilde{\itn{u}},\tilde{\itn{v}})=k_1\ldots k_n\nu$ (because all other (pre-)periodic sequences are contained in $\AD_\nu$). Iterating forward $n$ times, we obtain another triod $[\itn{t},\itn{u},\itn{v}]:=\Triod^{\circ n}[\tilde{\itn{t}},\tilde{\itn{u}},\tilde{\itn{v}}]$ with $\itn{t},\itn{u},\itn{v}\in\Omega^+(\nu)$ and $\itn{b}(\itn{t},\itn{u},\itn{v})=\nu$. We can assume w.l.o.g.\, that $\itn{t}=\sigma^{\circ j}(\nu)$, $\itn{u}=\sigma^{\circ k}(\nu)$, and $\itn{v}=\sigma^{\circ l}(\nu)$, where $j,k>0$ and $l\geq 0$. The triod $[\itn{t},\itn{u},\itn{v}]$ might be branched or linear. In both cases, by passing to an associated triod of external addresses, Lemma~\ref{lem:CombSeparatingRays} guarantees the existence of external addresses $\ad{s}',\ad{t},\ad{u}\in\mathcal{S}$ with the following properties:
	\begin{itemize}
		\item We have $\It(\ad{s}'~\vert~\ad{s})=\nu$, $\It(\ad{t}~\vert~\ad{s})=\sigma^{\circ j}(\nu)=\itn{t}$ and $\It(\ad{u}~\vert~\ad{s})=\sigma^{\circ k}(\nu)=\itn{u}$.
		\item We have $\ad{t}\in(\ad{s},\ad{s}')$ and $\ad{u}\in(\ad{s}',\ad{s})$.
	\end{itemize}
	By Proposition~\ref{pro:LandingBehavior}, the dynamic rays $g_{\ad{s}}$ and $g_{\ad{s}'}$ land at $0$, while $g_{\ad{t}}$ lands at  $p=E_\lambda^{\circ j}(0)$ and $g_{\ad{u}}$ lands at $q=E_\lambda^{\circ k}(0)$. By Lemma~\ref{lem:verticalLexicographicOrder}, the relations $\ad{t}\in(\ad{s},\ad{s}')$ and $\ad{u}\in(\ad{s}',\ad{s})$ imply that $p$ and $q$ are contained in different connected components of $\complex\setminus(\Tr(g_{\ad{s}})\cup\Tr(g_{\ad{s}'})\cup\{0\})$. But then, for one of these points, say $p$, the itineraries of the preimages of $p$ have different initial entries w.r.t.\ the dynamical partitions $\mathcal{D}$ induced by $g_{\ad{s}}$ and $\mathcal{D}'$ induced by $g_{\ad{s}'}$. At least one preimage of $p$ is itself a post-singular point and it has different itineraries w.r.t.\ the two dynamical partitions, contradicting Lemma~\ref{lem:DifferentPartitions}.

	Finally, let us prove $(3)$. We can assume that $[\itn{t},\itn{u},\itn{v}]$ is branched because otherwise the middle point of the triod is contained in $V_\nu$ anyway. Let $T:=(\ad{t}^i)_{i\in I}$, $U:=(\ad{u}^j)_{j\in J}$ and $V:=(\ad{v}^k)_{k\in K}$ be the sets of external addresses with itinerary $\itn{t}$, $\itn{u}$ and $\itn{v}$ respectively. By Lemma~\ref{lem:CombSeparatingRays} and Lemma~\ref{lem:Unlinked}, there are three external addresses $\ad{s}^1,\ad{s}^2,\ad{s}^3$ (w.l.o.g.\ indexed with respect to their cyclic order) with itinerary $\itn{b}(\itn{t},\itn{u},\itn{v})$ such that $T\subset(\ad{s}^1,\ad{s}^2)$, $U\subset(\ad{s}^2,\ad{s}^3)$ and $V\subset(\ad{s}^3,\ad{s}^1)$. If $\itn{t}\in\Omega^+(\star\nu)$, we set $\itn{t}':=\itn{t}$. Otherwise, $\itn{t}$ is the branch point of a triod $[\itn{t}^1,\itn{t}^2,\itn{t}^3]$ of post-singular points $\itn{t}^i\in\Omega^+(\star\nu)$ and by Lemma~\ref{lem:CombSeparatingRays} we can find an external address $\ad{t}\in(\ad{s}^1,\ad{s}^2)$ with itinerary $\itn{t}^j$ for a suitable value of $j\in\{1,2,3\}$. We set $\itn{t}':=\itn{t}^j$. Proceeding with $\itn{u}$ and $\itn{v}$ in the same way, we obtain a triod $[\itn{t}',\itn{u}',\itn{v}']$ of post-singular points $\itn{t}',\itn{u}',\itn{v}'\in\Omega^+(\star\nu)$, and it follows from Lemma~\ref{lem:CombSeparatingRays} (the if-direction), that $\itn{b}(\itn{t}',\itn{u}',\itn{v}')=\itn{b}(\itn{t},\itn{u},\itn{v})$. We have $\itn{b}(\itn{t}',\itn{u}',\itn{v}')\in V_\nu$ by definition of the formal vertex set.
\end{proof}
\begin{Def}[Separating sets]
If $v=v_{\itn{t}}$ is a pre-singular vertex, we set
\[
A_{v}:=\bigcup_{\It(\ad{t}~\vert~\ad{s})=\itn{t}}\Tr(g_{\ad{t}})\cup U_{\itn{t}}.
\]
If $v$ is not a pre-singular vertex, we set
\[
A_{v}:=\bigcup_{\It(\ad{t}~\vert~\ad{s})=\itn{t}}\Tr(g_{\ad{t}})\cup\{v\}.
\]
\label{def:RaysPlusLandingPoint}
\end{Def}

\begin{figure}[ht]
	\includegraphics[width=\textwidth]{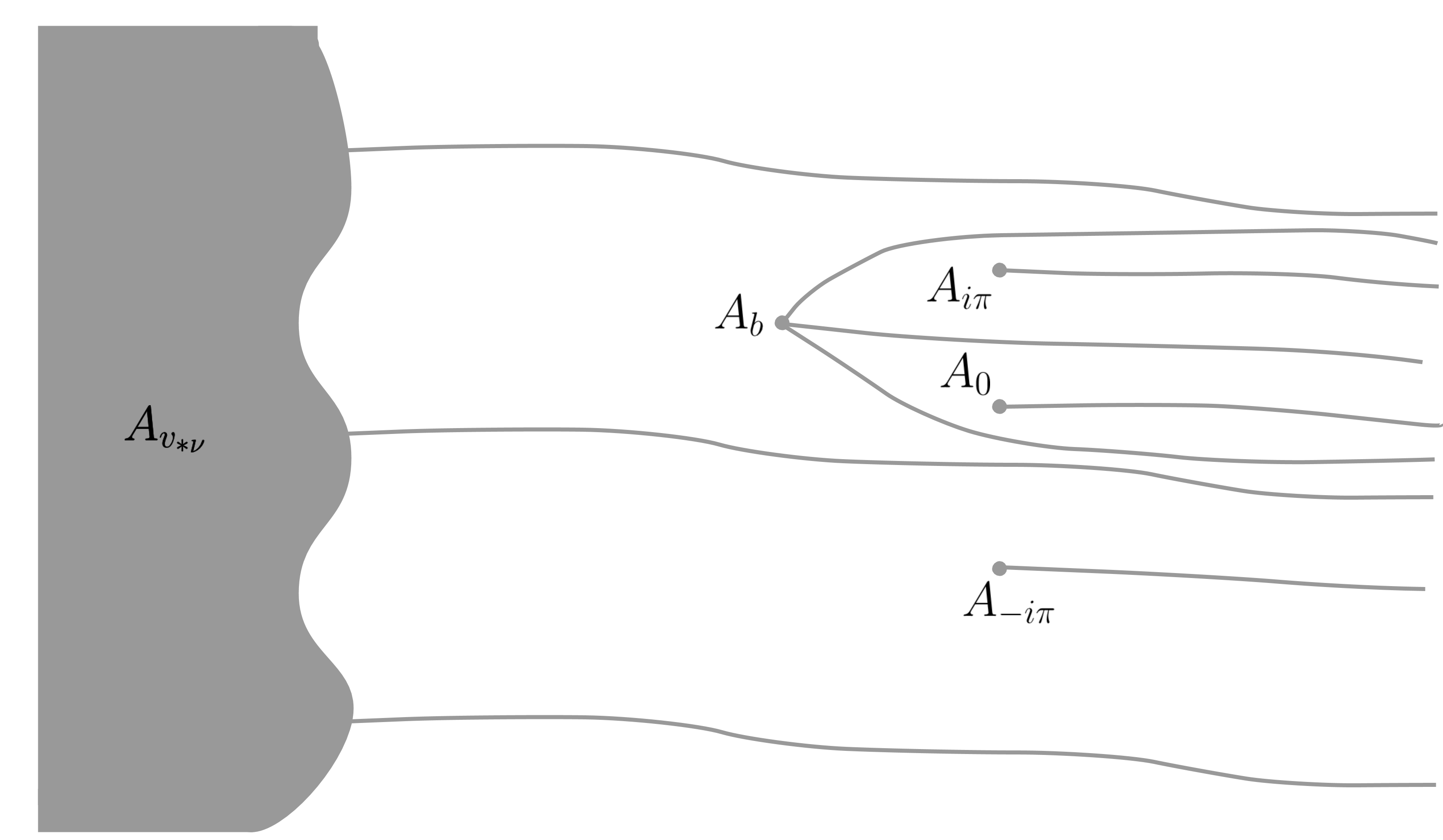}
	\caption{Sketch of the separating sets $A_v$ associated to the embedded vertices for the exponential map $E_{i\pi}$}
	\label{fig:SeparatingRays}
\end{figure}
It is clear (by the properties (1) and (2) imposed on $U$) that if $v,w\in V$ are different vertices then $A_v\cap A_w=\emptyset$. See Figure~\ref{fig:SeparatingRays} for a sketch of the separating sets $A_v$ for a particularly simple psf exponential map.

\begin{Lem}[Separating dynamic rays]
Let $[p,q,r]$ be a triod of vertex points. The triod $[p,q,r]$ is linear with $p$ in the middle if and only if $q$ and $r$ lie in different connected components of $\complex\setminus A_p$, and it is branched with branch point $b=b[p,q,r]\in V$ if and only if $p$, $q$ and $r$ lie in different connected components of $\complex\setminus A_b$.
\label{lem:SeparatingRays}\end{Lem}

\begin{proof}
This is just a matter of translating Lemma~\ref{lem:CombSeparatingRays} and Lemma~\ref{lem:VertexSet} into the language of dynamic rays and their landing points using Lemma~\ref{lem:verticalLexicographicOrder} and Proposition~\ref{pro:LandingBehavior}.
\end{proof}

We continue our construction of Homotopy Hubbard Trees. Choose a point $p_0\in\partial U\setminus\Tr(g_{\ad{s}})$ and write $p_{\star\nu, j}$ for the unique preimage of $p$ (under $E_\lambda$) in the sector $D_j$ of the dynamical partition. We choose disjoint arcs $\gamma_{p_{\star\nu,j}}\colon [0,1]\to U_{\star\nu}$ from $v_{\star\nu}$ to $p_{\star\nu,j}$ such that $\gamma_{p_{\star\nu,j}}((0,1))\subset\text{int}(U_{\star\nu})$. These curves are there to normalize the ends of the edges of the yet to be constructed Homotopy Hubbard Tree that have $v_{\star\nu}$ as an endpoint. For a pre-singular vertex $v_{\itn{t}}$ with $\itn{t}=k_1\ldots k_n\star\nu$, we define the pulled-back path $\gamma_{p_{\itn{t},j}}:=\log_{\ad{s},k_1}\ldots \log_{\ad{s},k_n}(\gamma_{\star\nu,j})$ with endpoints $v_{\itn{t}}$ and $p_{\itn{t},j}$. The points $p_{\itn{t},j}$ are just auxiliary points, that will lie on the interior of an edge at the end of the construction.

We call an arc $\gamma\colon [0,1]\to\complex$ between two different points $v,w\in V$ \emph{allowable} if the following conditions are satisfied:
\begin{itemize}
	\item If $r\in V$ is not a pre-singular vertex, then either $\Tr(\gamma)\cap A_{r}=\emptyset$ or $\Tr(\gamma)\cap A_{r}=\{r\}$.
	\item If $r=v_{\itn{t}}\in V$ is a pre-singular vertex, then either $\Tr(\gamma)\cap A_{r}=\emptyset$, or $\Tr(\gamma)\cap A_{r}=\Tr(\gamma_{p_{\itn{t},i}})$ for some $i\in\Z$ if $r$ is one of the endpoints of $\gamma$ or $\Tr(\gamma)\cap A_{r}=\Tr(\gamma_{p_{\itn{t},i}})\cup \Tr(\gamma_{p_{\itn{t},j}})$ for indices $i\neq j$ if $r$ lies in the interior of the path.
\end{itemize}
An allowable arc $\gamma$ satisfying $\gamma^{-1}(V)=\{0,1\}$ is called an \emph{edge}. We call two vertices \emph{incident} if they can be connected by an edge. Note that every allowable arc is a concatenation of finitely many edges; conversely, concatenating edges of a finite acyclic path always yields an allowable arc.
\begin{Lem}[Existence and uniqueness of allowable arcs]
	For any pair of distinct vertices $w,w'\in V$ there exists an allowable arc $\gamma\colon [0,1]\to\complex$ with $\gamma(0)=w$ and $\gamma(1)=w'$. If $\tilde{\gamma}\colon [0,1]\to\complex$ is another allowable arc from $w$ to $w'$, then $\gamma$ and $\tilde{\gamma}$ are homotopic rel $\bigcup_{v\in V} A_{v}$.
\label{lem:existenceAllowableCurve}\end{Lem}
\begin{proof}
Let $D$ be a connected component of $\complex\setminus\bigcup_{v\in V} A_{v}$. We start by investigating the topology of $D$ and $\partial D$. It is well known that if $A\subset\CC$ is closed and connected, every connected component of $\CC\setminus A$ is simply connected. As $\cl_{\CC}(\bigcup_{v\in V} A_{v})$ is closed and connected, $D$ is simply connected. Let $B$ be a connected component of $\partial D$. We have $B\subset\partial A_v$ for some $v\in V$. If $v=v_{\itn{t}}$ is a pre-singular vertex, then there exists a continuous bijection $\Gamma\colon (-\infty,\infty)\to B$ such that $\lim_{t\to -\infty}\Gamma(t)=\lim_{t\to +\infty}\Gamma(t)=\infty$ (in $\CC$) and $\Gamma(0)=p_{\itn{t},j}$ for some $j\in\Z$. We call $p_{\itn{t},j}$ the \emph{distinguished boundary point} of $D$ on $A_v$. If $v$ is not pre-singular, and there is only one dynamic ray $g$ landing at $v$, then we have $B=\partial D\cap A_v =\Tr(g)\cup\{v\}$. If there are at least two dynamic rays landing at $v$, we have $B=\partial D\cap A_v\Tr(g)\cup\Tr(g')\cup\{v\}$ for some rays $g,g'$ landing at $v$. In both cases, we call $v$ the \emph{distinguished boundary point} of $D$ on $A_v$.

We claim that $\partial D$ has at most two connected components. Assume to the contrary that there are vertices $v_1,v_2,v_3\in V$ such that $\partial D\cap A_{v_i}\neq\emptyset$ for all $i$.  Then, we can connect any two distinct $v_i$ and $v_j$ by a curve $\gamma_{(i,j)}$ such that $\Tr(\gamma)\cap A_v=\emptyset$ for all $v\in V\setminus\{v_i,v_j\}$.  Consider the triod $[v_1,v_2,v_3]$. If it is linear with middle point $v_i$, then $A_{v_{i-1}}$ and $A_{v_{i+1}}$ are separated by $A_{v_i}$ by Lemma~\ref{lem:SeparatingRays}, and if $[v_1,v_2,v_3]$ is branched with branch point $b$, then each pair $A_{v_i}$ and $A_{v_j}$ is separated by $A_b$ by Lemma~\ref{lem:SeparatingRays}. In both cases, we get a contradiction. Therefore, there are only three possibilities for the topology of $(D,\partial D)$: there exists a homeomorphism $\varphi\colon(D,\partial D)\to(\Omega_i,\partial\Omega_i)$, i.e., a homeomorphism $\varphi\colon\overline{D}\to\overline{\Omega_i}$ satisfying $\varphi(\partial D)=\partial\Omega_i$, for exactly one of three uniformizing domains $\Omega_i\subset\complex$, where
\begin{enumerate}
\item $\Omega_1:=\{z\in\complex\colon\Re(z)>0\}$ is the right half plane,
\item $\Omega_2:=\{z\in\complex\colon\Re(z)>0\}\setminus[1,+\infty)$ is the slit right half plane,
\item and $\Omega_3:=\{z\in\complex\colon 0<\Re(z)<1\}$ is a vertical strip.
\end{enumerate}

Let us now prove existence and uniqueness of allowable arcs. First, consider the special case that there exists a connected component $D$ of $\complex\setminus\bigcup_{v\in V} A_{v}$ such that $\partial D\cap A_w\neq\emptyset$ and $\partial D\cap A_{w'}\neq\emptyset$. Let $p_w,p_w'\in\partial D$ be the distinguished boundary points. The domain $D$ is of type (2) or (3), so there exists an arc $\delta\colon [0,1]\to\overline{D}$ from $p_w$ to $p_w'$ such that $\delta((0,1))\in D$, and $\delta$ is unique up to homotopy rel $\bigcup_{v\in V} A_{v}$. We obtain an edge $\gamma$ connecting $w$ and $w'$ after concatenating $\delta$ with $\gamma_{p_w}$ or $\gamma_{p_{w'}}$ if $p_w\neq w$ or $p_{w'}\neq w'$.

Finally, consider the general case, and set $v_0:=w$. Let $E(v_0)$ be the component of $\complex\setminus A_{v_0}$ containing $w'$, and let $D_1$ be the unique component of $\complex\setminus\bigcup_{v\in V} A_{v}$ such that $\partial D_1\cap A_{v_0}\neq\emptyset$ and $D_1\subset E(v_0)$. The component $D_1$ is of type (2) or (3) since $E(v_0)$ contains at least the vertex $w'$, so there exists $v_1\in V$ such that $\partial D_1\cap A_{v_1}\neq\emptyset$. If $v_1=w'$, we are reduced to the special case. Else, $D_1$ is of type (3), and there is a unique component $E(v_1)\subset E(v_0)$ of $\complex\setminus A_{v_1}$ containing $w'$. We continue this argument inductively, and as the number of vertices contained in $E(v_i)$ decreases in each step, we obtain a finite sequence $w=v_0, v_1, \ldots, v_n=w'$ of vertices such that $A_{v_{i-1}}$ and $A_{v_{i}}$ both intersect $\partial D_i$. The concatenation $\gamma:=\gamma_0\cdot\gamma_1\cdot\ldots \cdot\gamma_{n-1}$, where $\gamma_i$ is an edge from $v_i$ to $v_{i+1}$, is an allowable path from $w$ to $w'$, and it is unique up to homotopy rel $\bigcup_{v\in V} A_{v}$ because every part $\gamma_i$ is.
\end{proof}

\begin{Thm}[Existence of Homotopy Hubbard Trees]
For every post-singularly finite exponential map $E_\lambda$ there exists a Homotopy Hubbard Tree.
\label{thm:ExistenceProof}\end{Thm}
\begin{proof}
For every pair $w,w'\in V$ of incident vertices, choose an edge $e_{(w,w')}\colon [0,1]\to\complex$ connecting $w$ and $w'$. We choose the unique edge $e_{(0,w_0)}\colon [0,1]\to\complex$ connecting $0$ to another vertex $w_0$ in such a way that it intersects $\partial U$ only at the point $p_0\in\partial U\setminus\Tr(g_{\ad{s}})$ (see the paragraph before Lemma~\ref{lem:existenceAllowableCurve} for the definition of $p_0$). Define
\[
H:=\bigcup_{w,w'\in V~\text{incident}}\Tr(e_{(w,w')}).
\]
As every edge is contained in some component of $\complex\setminus\bigcup_{v\in V}A_v$ (except possibly for its normalized ends), and every such component contains at most one edge by Lemma~\ref{lem:existenceAllowableCurve}, different edges are disjoint (except possibly for their endpoints). Clearly, they are not homotopic rel
$\bigcup_{v\in V}A_v$. The embedded graph $H$ is connected by Lemma~\ref{lem:existenceAllowableCurve}, and the existence of a cycle of edges would contradict uniqueness in Lemma~\ref{lem:existenceAllowableCurve}. Hence, $H$ is an embedded tree. Moreover, every endpoint of $H$ is a post-singular point: every vertex $v\in V\setminus P(E_\lambda)$ is the middle point of a triod of the form $[p_1,v,p_2]$ with $p_1,p_2\in P(E_\lambda)$ and the in-tree connection of $p_1$ and $p_2$ has to contain $v$ by Lemma~\ref{lem:SeparatingRays} and the definition of allowable arcs.

Next, we want to see that $H':=[P(E_\lambda)]_{\widehat{E_\lambda}^{-1}(H)}$ is homotopic to $H$ rel $P(E_\lambda)$. By Lemma~\ref{lem:HomotopingMinusInfty}, we can find a homotopy $I_0$ between $H'$ and an embedded tree $H''\subset\complex$ rel $\CT\setminus\widehat{E_{\lambda}}^{-1}(U)$, and by our choice of the edge $e_{(0,w_0)}$, we can make sure that $H''\cap\widehat{E_{\lambda}}^{-1}(U)=H\cap\widehat{E_{\lambda}}^{-1}(U)=\bigcup_{i\in\{1,\ldots, k\}}\Tr(\gamma_{p_{\star\nu, j_i}})$ for certain $j_i\in\Z$.

For a vertex $v\neq v_{\star\nu}$, the map $E_\lambda\colon A_v\to A_{E_\lambda(v)}$ is a homeomorphism, and if $v=v_{\itn{t}}$ is pre-singular, we have $E_\lambda(\bigcup_j\Tr(\gamma_{\itn{t},j}))=\bigcup_j\Tr(\gamma_{\sigma(\itn{t}),j})$. If $v$ is not pre-singular, we have $H\cap A_{E_\lambda(v)}=\{E_\lambda(v)\}$, and therefore $\widehat{E_\lambda}^{-1}(H)\cap A_v=\{v\}$. If instead $v=v_{\itn{t}}\neq v_{\star\nu}$ is pre-singular, we have $H\cap A_{E_\lambda(v)}\subset\bigcup_j\Tr(\gamma_{\sigma(\itn{t}),j})$, and therefore $\widehat{E_\lambda}^{-1}(H)\cap A_v\subset\bigcup_j\Tr(\gamma_{\itn{t},j})$. Hence, for every pair of distinct post-singular points $p,q\in P(E_\lambda)$, the in-tree connection $[p,q]_{H''}$ is an allowable arc.

It follows that $V\subset H''$, because every $v\in V\setminus (P(E_\lambda)\cup\{v_{\star\nu}\})$ is the branch point of a post-singular triod. 
Moreover, the trees $H$ and $H''$ have the same graph structure on $V$. To see this, observe that $H$ is defined to contain an edge for every pair of incident vertices. Hence, vertices that are connected by an edge in $H''$ are also connected by an edge in $H$. Conversely, the tree $H''$ must contain all edges in $H$ since otherwise it would be disconnected.
By the uniqueness of allowable arcs, there exists a homotopy $I_1$ between $H''$ and $H$ rel $\bigcup_{v\in V}A_v$. We obtain a homotopy $I$ between $H$ and $H'$ rel $P(E_\lambda)$ by concatenating $I_0$ and $I_1$.

It remains to show that the induced self-map $f:=\widehat{E_\lambda}\circ I^1\colon H\to H$ is expansive (see the paragraph before Definition~\ref{def:HomotopyHubbardTree} for the definition of expansivity). By construction, we have $I^1(v_{\star\nu})=-\infty$, and we deduce that $V_f=V$ for the set of marked points of the self-map $f$. Let $p,q\in V_f=V$ be two different marked points of $f$. By construction, two vertices $v_{\itn{t}},v_{\itn{u}}\in V\setminus\{v_{\star\nu}\}$ are contained in different branches of $H$ at $v_{\star\nu}$ if and only if the first entries of $\itn{t}$ and $\itn{u}$ are different. As different vertices have different itineraries, there exists a smallest $n\geq 0$ such that the itineraries of $f^{\circ n}(p)$ and $f^{\circ n}(q)$ have different initial entries. Hence, we have $v_{\star\nu}\in[f^{\circ n}(p),f^{\circ n}(q)]\subset f^{\circ n}[p,q]$.
\end{proof}

\section{Classification of post-singularly finite exponential maps}
\label{sec:classification}

Post-singularly finite exponential maps have been classified in terms of the external addresses of the dynamic rays landing at the singular value in \cite{LSV}, and we have stated their main result in Theorem~\ref{thm:SpiderClassification}. As an application of our construction, we give another combinatorial classification of this class of maps in terms of abstract Hubbard Trees. Post-singularly finite polynomials have been classified combinatorially by Poirier in \cite{P} in terms of so called \emph{abstract Hubbard Trees}. These are graph-theoretic trees equipped with a self-map and certain extra information, but without an embedding into $\complex$. The exponential dynamical trees from Definition~\ref{def:ExpDynamicalTree} come quite close to what we mean by an abstract exponential Hubbard Tree. We have to add a cyclic order on the branches of the tree at each marked point in order to specify a homotopy type of embeddings of the tree into the complex plane.
\begin{Def}[Abstract exponential Hubbard Tree]
An \emph{abstract exponential Hubbard Tree} $(H,f,B_{k_i},\angle_v)$ is an exponential dynamical tree $(H,f,B_{k_i})$ together with a cyclic order $\angle_v$ on the set of branches of $H$ at $v$ for each marked point $v\in V_f\setminus\{v_T\}$ such that $f$ preserves the cyclic order. Two abstract Hubbard Trees $(H,f,B_{k_i},\angle_v)$ and $(\tilde{H},\tilde{f},\tilde{B}_{k_i}, \angle_{\tilde{v}})$ are called \emph{equivalent} if they are equivalent as exponential dynamical trees and in addition the homeomorphism $\varphi\colon H\to\tilde{H}$ from Definition~\ref{def:ExpDynamicalTree} can be chosen to preserve the cyclic order at marked points.
\end{Def}
Let $H$ be a Homotopy Hubbard Tree for the post-singularly finite exponential map $E_\lambda$. In Section~\ref{sec:triodAlgorithm} we have seen that for any choice of induced self-map $f\colon H\to H$ the triple $(H,f,B_{k_i})$ is an exponential dynamical tree. We have also seen that if $\tilde{H}$ is another Homotopy Hubbard Tree for $E_\lambda$ with choice of induced self-map $\tilde{f}\colon \tilde{H}\to\tilde{H}$, the exponential dynamical trees $(H,f,B_{k_i})$ and $(\tilde{H},\tilde{f},\tilde{B}_{k_i})$ are equivalent. As $H$ is embedded into the complex plane, the branches of $H$ at a marked point $v\in V_f\setminus\{v_T\}$ come equipped with a natural cyclic order $\angle_v$ (of course this also holds for $v_T$, but the cyclic order of edges at $v_T$ is already contained in the sector information). Also, the trees $H$ and $\tilde{H}$ are homotopic in $\complex$ rel $P(E_\lambda)$, and such a homotopy preserves the cyclic order of branches at marked points. It follows that there is a well-defined map
\[
\mathcal{F}\colon \{\text{psf exponential maps}\}\to\{\text{abstract exponential Hubbard Trees}\},
\]
where we consider abstract exponential Hubbard Trees up to equivalence. In this section, we prove that $\mathcal{F}$ is a bijection.

The major tool that has made classification results of this kind possible for certain classes of rational functions is Thurston's topological characterization of rational maps. An analogous result for exponential maps has been established in \cite{HSS}. In order to state this result precisely, we need to introduce some terminology.

\begin{conv}
In the following, $\mathbb{S}^2$ denotes an oriented topological $2$-sphere with two distinguished points $0$ and $\infty$. All homeomorphisms and coverings will be understood to be orientation-preserving. We also write $\R^2=\mathbb{S}^2\setminus\{\infty\}$.
\end{conv}
\begin{Def}[Topological exponential maps]
A covering map $g\colon \mathbb{R}^2\to\mathbb{R}^2\setminus\{0\}$ is called a \emph{topological exponential map}. It is called \emph{post-singularly finite} if the orbit of $0$ is finite, hence preperiodic. The \emph{post-singular set} is $P(g):=\bigcup_{n\geq 0}g^{\circ n}(0)$.
\label{def:TopExpMap}
\end{Def}
\begin{Def}[Thurston equivalence]
Two post-singularly finite topological exponential maps $f$ and $g$ with post-singular sets $P(f)$ and $P(g)$ are called \emph{Thurston equivalent} if there are two homeomorphisms $\varphi_0,\varphi_1\colon \mathbb{R}^2\to\mathbb{R}^2$ satisfying $\left.\varphi_0\right|_{P(f)}=\left.\varphi_1\right|_{P(f)}$ and $P(g)=\varphi_0(P(f))=\varphi_1(P(f))$ such that the diagram
\[ \begin{tikzcd}
(\mathbb{R}^2,P(f)) \arrow{r}{\varphi_1} \arrow[swap]{d}{f} & (\mathbb{R}^2, P(g)) \arrow{d}{g} \\%
(\mathbb{R}^2,P(f)) \arrow{r}{\varphi_0}& (\mathbb{R}^2, P(g))
\end{tikzcd}
\]
commutes, and $\varphi_0$ is homotopic to $\varphi_1$ relative to $P(f)$.
\label{def:ThurstonEquiv}
\end{Def}
The main result of \cite{HSS} is a criterion for a topological exponential map to be Thurston equivalent to a (necessarily post-singularly finite) holomorphic exponential map. There exists an equivalent holomorphic map if and only if the following obstruction does not occur; see below.
\begin{Def}[Essential curves and Levy cycles]
Let $g$ be a post-singularly finite topological exponential map. A simple closed curve $\gamma\subset\mathbb{S}^2\setminus(P(g)\cup\{\infty\})$ is called \emph{essential} if both connected components of $\mathbb{S}^2\setminus\gamma$ contain at least two points of $P(g)\cup\{\infty\}$. A \emph{Levy cycle} of $g$ is a finite sequence of disjoint essential simple closed curves $\gamma_0$, $\gamma_1$, \ldots, $\gamma_{m-1}$, $\gamma_m=\gamma_0$ such that for $i=0,1,\ldots,m-1$, some component $\gamma_i'$ of $g^{-1}(\gamma_{i+1})$ is homotopic to $\gamma_i$ in $\R^2$ relative to $P(g)$ and $g\colon \gamma_i'\to\gamma_{i+1}$ is a homeomorphism.

Let $U_i'$ be the bounded component of $\mathbb{R}^2\setminus\gamma_i'$ and let $U_i$ be the bounded component of $\mathbb{R}^2\setminus\gamma_i$. If all restrictions $g\colon \overline{U_{i'}}\to\overline{U_{i+1}}$ are homeomorphisms, then the Levy cycle is called \emph{degenerate}.
\label{def:LevyCycle}\end{Def}
It is easy to see that for topological exponential maps every Levy cycle is degenerate and that Levy cycles are preserved under Thurston equivalence. Furthermore, a routine hyperbolic contraction argument shows that a post-singularly finite exponential map does not have a Levy cycle. The following theorem is the main result of \cite{HSS}.
\begin{Thm}[Topological characterization of exponential maps]
A post-singularly finite topological exponential map is Thurston equivalent to a post-singularly finite holomorphic exponential map if and only if it does not admit a degenerate Levy cycle. The holomorphic exponential map is unique up to conjugation with an affine map.
\label{thm:TopChar}\end{Thm}
\begin{remark}
In our parametrization $\lambda\mapsto\lambda\exp(z)$ of the exponential parameter space, with $\lambda\in\complex\setminus\{0\}$, no two distinct $E_\lambda$ and $E_{\tilde{\lambda}}$ are affinely conjugate. By Theorem~\ref{thm:TopChar}, two post-singularly finite maps $E_\lambda$ and $E_{\tilde{\lambda}}$ are never Thurston equivalent for different parameters $\lambda,\tilde{\lambda}\in\complex\setminus\{0\}$.
\end{remark}

We show that distinct psf exponential maps always yield abstract Hubbard Trees that are not equivalent. For the proof, we need a result from \cite[Corollary~6.6]{BFH} about extensions of maps between embedded graphs.

\begin{Lem}[Extension of graph homeomorphisms]
  Let $\Gamma_1,\Gamma_2\subset\mathbb{S}^2$ be connected embedded graphs, and let $f\colon \Gamma_1\to\Gamma_2$ be a homeomorphism. Then $f$ extends to an orientation-preserving homeomorphism $\hat{f}\colon \mathbb{S}^2\to\mathbb{S}^2$ if and only if $f$ preserves the cyclic order of the edges at all branch points of $\Gamma_1$.\qed
\label{lem:GraphExtension}\end{Lem}

\begin{Thm}[Different maps have non-equivalent trees]
Let $E_\lambda$ and $E_{\tilde{\lambda}}$ be two post-singularly finite exponential maps. Let $H$ be a Homotopy Hubbard Tree for $E_\lambda$, and let $\tilde{H}$ be a Homotopy Hubbard Tree for $E_{\tilde{\lambda}}$. If $H$ and $\tilde{H}$ yield the same abstract Hubbard Tree, we have $\lambda=\tilde{\lambda}$.
\label{thm:F_Injective}\end{Thm}
\begin{proof}
Choose dynamic rays $g$ (resp.\ $\tilde{g}$) of $E_\lambda$ (resp.\ $E_{\tilde{\lambda}}$) landing at the singular value. By Proposition~\ref{pro:AccessesContainRays}, we can assume w.l.o.g.\ that $H$ (resp.\ $\tilde{H}$) does not intersect $\Tr(g)$ (resp.\ $\Tr(\tilde{g}))$.
By hypothesis, there exists a homeomorphism $\varphi\colon H\to\tilde{H}$ that restricts to a conjugation between $E_\lambda$ and $E_{\tilde{\lambda}}$ on $P(E_\lambda)$, and preserves the cyclic order at each marked point of $H$.
By Lemma~\ref{lem:GraphExtension}, $\varphi$ can be extended to a homeomorphism $\varphi_0\colon \complex\to\complex$, and we can choose $\varphi_0$ to satisfy $\varphi_0(\Tr(g))=\Tr(\tilde{g})$. In analogy to Lemma~\ref{lem:Lifting}, one shows that there exists a lift $\varphi_1\colon \complex\to\complex$ of $\varphi_0$ satisfying $\varphi_0\circ E_\lambda=E_{\tilde{\lambda}}\circ\varphi_1$ and $\varphi_1(0)=0$.
We want to see that $\varphi_0$ is homotopic to $\varphi_1$ rel $P(E_\lambda)$. This implies that $E_\lambda$ and $E_{\tilde{\lambda}}$ are Thurston equivalent, so we have $\lambda=\tilde{\lambda}$ by Theorem~\ref{thm:TopChar} and the remark following it.

Let us first see that $\varphi_1$ restricts to a conjugation between $E_\lambda$ and $E_{\tilde{\lambda}}$ on $P(E_\lambda)$, just as $\varphi_0$ does. Let $\mathcal{D}$ be the dynamical partition for $E_\lambda$ w.r.t.\ $g$, and let $\mathcal{\tilde{D}}$ be the dynamical partition for $E_{\tilde{\lambda}}$ w.r.t.\ $\tilde{g}$. It follows from $\varphi_0(\Tr(g))=\Tr(\tilde{g})$ and $\varphi_1(0)=0$ that the lifted map $\varphi_1$ maps the partition sector $D_k$ homeomorphically onto the sector $\tilde{D}_k$ of the same index. Hence, for all $n\geq 0$, $\varphi_1$ sends the unique preimage of $E_\lambda^{\circ n}(0)$ in the sector $D_k$ to the unique preimage of $E_{\tilde{\lambda}}^{\circ n}(0)$ in $\tilde{D}_k$. As $E_\lambda$ and $E_{\tilde{\lambda}}$ have the same kneading sequence (remember that the abstract Hubbard Tree in particular contains the sector information for the post-singular points), it follows that $\varphi_1(E_\lambda^{\circ (n-1)}(0))=E_{\tilde{\lambda}}^{\circ (n-1)}(0)$ for all $n>0$.

Set $H':=[P(E_\lambda)]_{\widehat{E_\lambda}^{-1}(H)}$ and $\tilde{H}':=[P(E_{\tilde{\lambda}})]_{\widehat{E_{\tilde{\lambda}}}^{-1}(H)}$. The map $\varphi_1$
extends to a map $\widehat{\varphi}_1\colon\CT\to\complex_{\tilde{T}}$ by setting $\widehat{\varphi}_1(-\infty):=-\infty$.
As $\varphi_1(P(E_\lambda))=P(E_{\tilde{\lambda}})$, we have $\widehat{\varphi}_1(H')=\tilde{H}'$. By Lemma~\ref{lem:HomotopingMinusInfty}, there exists an embedded tree $H''\subset\complex$ such that $H''$ is homotopic to $H'$ rel $P(E_\lambda)$. Any homotopy between $H'$ and $H''$ rel $P(E_\lambda)$ is pushed forward by $\varphi_1$ to a homotopy between $\tilde{H}'$ and an embedded tree $\tilde{H}'':=\varphi_1(H'')\subset\complex$ rel $P(E_{\tilde{\lambda}})$. By Lemma~\ref{lem:HomotopyImpliesAmbIsotopy}, $H''$ is ambient isotopic to $H$ rel $P(E_\lambda)$, and $\tilde{H}''$ is ambient isotopic to $\tilde{H}$ rel $P(E_{\tilde{\lambda}})$. Stated differently, there exist homeomorphisms $\Psi_0,\Psi_1\colon \complex\to\complex$ such that $\Psi_0$ is isotopic to $\text{id}$ rel $P(E_\lambda)$, $\Psi_1$ is isotopic to $\text{id}$ rel $P(E_{\tilde{\lambda}})$, and $\varphi_1':=\Psi_1\circ\varphi_1\circ\Psi_0$ satisfies $\varphi_1'(H)=\tilde{H}$. The composition $\varphi_1'$ is isotopic to $\varphi_1$ rel $P(E_\lambda)$, and it is equal to $\varphi_0$ on the set of marked points of $H$. By \cite[Corollary~6.3]{BFH}, there exists a homeomorphism $\Psi\colon \complex\to\complex$ isotopic to $\text{id}$ relative to the set of marked points of $H$ (and in particular relative to $P(E_\lambda)$) such that $\varphi_1'\circ\Psi=\varphi_0$. But $\varphi_1'\circ\Psi$ is isotopic to $\varphi_1$ rel $P(E_\lambda)$, so $E_\lambda$ and $E_{\tilde{\lambda}}$ are Thurston equivalent.
\end{proof}

Next, we explain how to obtain a topological exponential map from an abstract exponential Hubbard Tree. For convenience, we are going to embed abstract Hubbard Trees into the complex plane $\complex$ (or suitable extensions of it). Note, however, that the complex structure of $\complex$ does not play any role; it only simplifies the construction. At first, we define a suitable embedding of the abstract tree into the complex plane. Afterwards, we show that the self-map of the embedded tree extends to a map on the plane, and that this map is a topological exponential map. The constructions are inspired by and similar to \cite[Chapter~5]{LSV}, so we skim some of them and refer to \cite{LSV} for more details.

\begin{figure}[ht]
\includegraphics[width=\textwidth]{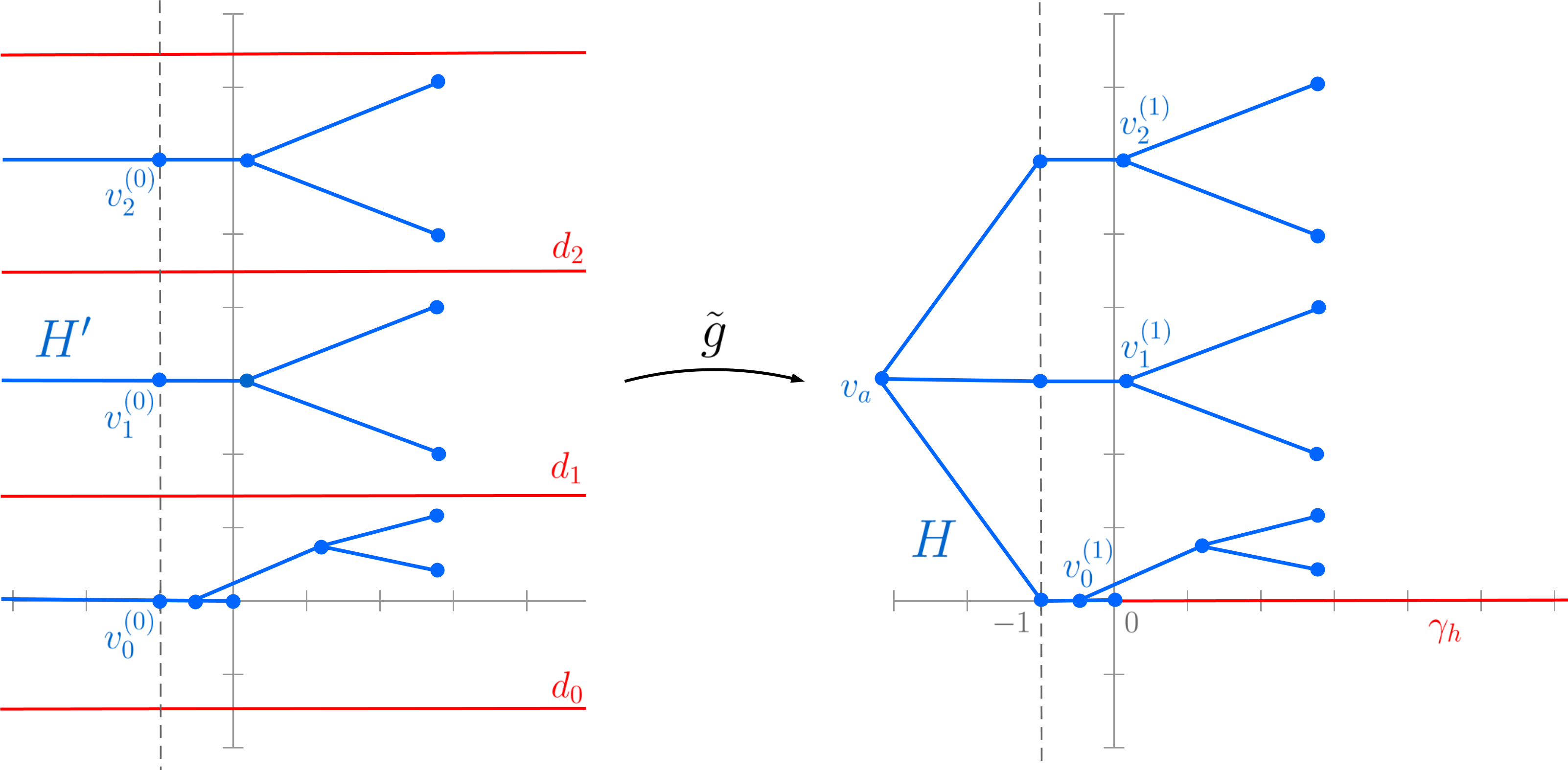}
\caption{Sketch of the embedded trees $H:=\iota(\mathtt{H})$ and $H':=\iota{(\mathtt{H})}$. The two embedded trees agree on $\{z\in\complex~\vert~\Re(z)\geq -1\}$.}
\label{fig:GraphMap}\end{figure}

At first, we specify a general mapping layout for all of our topological exponential maps which is independent of the abstract Hubbard Tree that we want to realize. See Figure~\ref{fig:GraphMap} for a sketch of the upcoming construction. Let $\gamma_h\colon (0,\infty)\to\complex$, $\gamma_h(t):=t$ be a parametrization of the horizontal line of positive reals. For each $k\in\Z$, let $d_k\colon (-\infty,+\infty)\to\complex$, $d_k(t):= t+(2k-1)\pi i$ be a parametrization of the straight horizontal line at constant imaginary part $(2k-1)\pi$. We define a continuous map $g_0\colon \bigcup\Tr(d_k)\to\Tr(\gamma_h)$ by $g_0(d_k(t)):=\gamma_h(\exp(t))$. Denote by $D_j$ the connected component of $\complex\setminus\bigcup\Tr(d_k)$ bounded by $\Tr(d_j)$ and $\Tr(d_{j+1})$.

Let $\mathtt{(H,f,B_{k_i},\angle_v)}$ (we use a different font here to distinguish the abstract tree from the embedded tree) be an abstract Hubbard Tree. The branches $\mathtt{B_{k_0}, \ldots, B_{k_n}}$ of $\mathtt{H}$ at $\mathtt{v_T}$ are indexed by distinct integers $\mathtt{k_i\in\Z}$, where $\mathtt{k_0}=0$ (compare Definition~\ref{def:ExpDynamicalTree}). Denote by $\mathtt{v_{k_i}^{(1)},\ldots, v_{k_i}^{(j_{k_i})}\in V_f\cap B_{k_i}}$ the marked points of $\mathtt{(H,f,B_{k_i},\angle_v)}$ contained in the branch $\mathtt{B_{k_i}}$ where $\mathtt{v_{k_i}^{(1)}}$ is the unique marked point in $\mathtt{B_{k_i}}$ incident to $\mathtt{v_T}$. Choose auxiliary points $\mathtt{v_{k_i}^{(0)}\in(v_T,v_{k_i}^{(1)})}$. We define an embedding $\iota\colon \mathtt{H}\to\complex$ in the following way:
\begin{itemize}
	\item For $i\in\{0,\ldots, n\}$, set $\iota(\mathtt{v_{k_i}^{(0)}}):= -1 + 2k_i\pi i=:v_{k_i}^{(0)}$.
	\item Choose an arbitrary point $v_T\in\{z\in\complex~\vert~\Re(z)<-1\}$, and choose pairwise disjoint arcs $\gamma_i\colon [0,1]\to\complex$ connecting $v_T$ to $v_{k_i}^{(0)}$ such that $\gamma_i((0,1))\subset\{z\in\complex~\vert~\Re(z)<-1\}$. Let $\iota\colon [\mathtt{v_T},\mathtt{v_{k_i}^{(0)}}]\to\Tr(\gamma_i)$ to be a homeomorphism satisfying $\iota(\mathtt{v_T})=v_T$ and $\iota(\mathtt{v_{k_i}^{(0)}})=v_{k_i}^{(0)}$.
	\item For $i\in\{1,\ldots, n\}$, choose distinct points $v_{k_i}^{(1)},\ldots, v_{k_i}^{(j_{k_i})}\in D_{k_i}$ satisfying $\Re(v_{k_i}^{(j)})>-1$. Define $\iota\colon {\mathtt{B_{k_i}\setminus[v_T,v_{k_i}^{(0)}]}}\to D_{k_i}\cap\{z\in\complex~\vert~\Re(z)>-1\}$ to be an embedding satisfying $\iota(\mathtt{v_{k_i}^{(j)}})=v_{k_i}^{(j)}$ such that the cyclic order of branches at $v_{k_i}^{(j)}$ coincides with the cyclic order of the corresponding branches of the abstract tree at $\mathtt{v_{k_i}^{(j)}}$. For the embedding of $\mathtt{B_0}$, we require $\iota(\mathtt{f(v_T)})=0$ and $\iota(\mathtt{B_0})\cap\Tr(\gamma_h)=\emptyset$ in addition. Note that every oriented topological tree is planar, so there always exists such an embedding.
\end{itemize}

The image $H:=\iota(\mathtt{H})\subset\complex$ is an embedded tree. Next, we define a slightly different embedding $\iota'$ of $\mathtt{H}$ whose image will be the subset $H'$ of the preimage tree of $H$ spanned by $P(g)$ under the yet to be constructed topological exponential map $g$.

Let $-\infty$ and $+\infty$ be abstract points not contained in $\complex$ and let $C_{\infty}:=\complex\cup\{-\infty,+\infty\}$ denote the extension of the complex plane by these abstract points.
We turn $\complex_{\infty}$ into a topological space by defining $U_n:=\{z\in\complex~\vert~\Re(z)<-n\}\cup\{-\infty\}$ to be a neighborhood basis of $-\infty$ and
$V_n:=\{z\in\complex~\vert~\Re(z)>n\}\cup\{+\infty\}$ to be a neighborhood basis of $+\infty$. For every $i\in\{0,\ldots, n\}$ let $\tilde{\gamma}_i\colon [0,1]\to\complex_{\infty}$ be an arc connecting $-\infty$ to $v_{k_i}^{(0)}$ and satisfying $\gamma_i((0,1))\subset D_{k_i}\cap\{z\in\complex~\vert~\Re(z)< -1\}$. Let $\left.\iota'\right|_{[\mathtt{v_T},\mathtt{v_{k_i}^{(0)}}]}\colon [\mathtt{v_T},\mathtt{v_{k_i}^{(0)}}]\to\Tr(\tilde{\gamma}_i)$ be a homeomorphism where $\iota'(\mathtt{v_T})=-\infty$ and $\iota'(\mathtt{v_{k_i}^{(0)}})=v_{k_i}^{(0)}$. On the complement of the $[\mathtt{v_T},\mathtt{v_{k_i}^{(0)}}]$ define $\iota'$ to be equal to $\iota$ and set $H':=\iota'(\mathtt{H})\subset\complex_{\infty}$.

The next step is to define a map from $H'$ to $H$. From now on, let $d_i\colon [-\infty,+\infty]\to\complex_{\infty}$ denote the extended curve, where $d_i(-\infty)=-\infty$ and $d_i(+\infty)=+\infty$. We also extend $\gamma_h$ to a curve $\gamma_h\colon [0,+\infty]\to\CC$, where $\gamma_h(0)=0$ and $\gamma_h(+\infty)=\infty$. Then we define an extended map $g_0\colon \bigcup\Tr(d_i)\to\Tr(\gamma)$ by setting $g_0(-\infty)=0$ and $g_0(+\infty)=\infty$. We define a graph map
\[
\tilde{g}\colon H'\cup(\bigcup\Tr(d_i))\to H\cup\Tr(\gamma_h),~\tilde{g}(z):=
\begin{cases}
g_0(z) & \text{if}~z\in\bigcup\Tr(d_i),\\
\iota\circ \mathtt{f}\circ(\iota')^{-1}(z) & \text{if}~z\in H'.
\end{cases}
\]

Our next goal is to show that the graph map $\tilde{g}$ can be extended to a post-singularly finite topological exponential map. Let us see how to use Lemma~\ref{lem:GraphExtension} in our setting. There exist continuous maps $\varphi_1\colon \overline{\D}\to\text{cl}_{\complex_{\infty}}(D_i)$ and $\varphi_2\colon \overline{\D}\to\text{cl}_{\CC}(\complex\setminus\Tr(\gamma_h))$ that restrict to orientation-preserving homeomorphisms from $\D$ onto $D_i$ and $\complex\setminus\Tr(\gamma_h)$ respectively. The map $\varphi_1$ is a homeomorphism, while every interior point of $\Tr(\gamma_h)$ has two preimages under $\varphi_2$. We choose $\varphi_1$ to satisfy $\varphi_1(1)=+\infty$, $\varphi_1(-1)=-\infty$, and $\varphi_2$ to satisfy $\varphi_2(-1)=0$, $\varphi_2(1)=\infty$. Then, there exist inverse branches $\psi_u\colon \Tr(\gamma_h)\to\partial\D\cap\{z\in\complex~\vert~\Im(z)\geq 0\}$ and $\psi_l\colon \Tr(\gamma_h)\to\partial\D\cap\{z\in\complex~\vert~\Im(z)\leq 0\}$ of $\varphi_2$. We define a homeomorphism $f\colon \partial\D\to\partial\D$ via
\[
f(z):=\begin{cases} \psi_u(g_0(\varphi_1(z))) & \text{if} \Im(z)\geq 0, \\ \psi_l(g_0(\varphi_1(z))) & \text{if} \Im(z)\leq 0. \end{cases}
\]
We view $\overline{\D}\subset\CC$ as a subset of the Riemann sphere. Setting $\Gamma_1:=\varphi_1^{-1}(\partial D_i\cup (H'\cap D_i))$, $\Gamma_2:=\varphi_2^{-1}(\tilde{g}(H'\cap D_i)\cup\Tr(\gamma_h))$, the $\Gamma_i\subset\CC$ are connected embedded graphs. We extend $f$ to a homeomorphism $f\colon \Gamma_1\to\Gamma_2$ by setting $f(z)=\varphi_2^{-1}(\tilde{g}(\varphi_1(z)))$ for $z\in\Gamma_1\setminus\partial\D$. Then $f$ preserves the circular order at branch points of $\Gamma_1$ because $\tilde{g}$ preserves the circular order at branch points. By Lemma~\ref{lem:GraphExtension}, there exists an orientation-preserving homeomorphism $\hat{f}\colon \CC\to\CC$ extending $f$ and satisfying $\hat{f}(\overline{\D})=\overline{\D}$. We define $g(z):=\varphi_2\circ\hat{f}\circ\varphi_1^{-1}$ for $z\in D_i$.

Pasting these extensions together, we obtain a globally defined map $g\colon\complex\to\complex\setminus\{0\}$. As it is a covering map, it is a post-singularly finite topological exponential map. To see that $g$ is Thurston equivalent to a holomorphic map, we have to show that $g$ does not admit a Levy cycle. We can assume w.l.o.g.\ that $\gamma_h$ is preperiodic (as a set) under the dynamics of $g$, and the forward images of $\gamma_h$ are pairwise disjoint (except possibly for their landing points). By construction, the boundary of the partition $(D_j)_{j\in\Z}$ consists of the preimages of $\gamma_h$, and by expansivity of the abstract Hubbard Tree, different post-singular points have different itineraries w.r.t.\ $(D_j)_{j\in\Z}$.

\begin{Lem}\label{lem:LevyCycle}
The map $g$ does not admit a Levy cycle.
\end{Lem}
\begin{proof}
Assume to the contrary that $g$ admits a Levy cycle $\Lambda=\{\delta_0,\delta_1,\ldots, \delta_k=\delta_0\}$. Set $\gamma_0:=\gamma_h$ and denote by $\gamma_i:=g^{\circ i}(\gamma_0)$ the forward iterates of $\gamma_0$. By construction $\gamma_0$ is preperiodic as a set. For curves $\delta_j$ and $\gamma_i$ let
\[
\vert[\delta_j]\cap[\gamma_i]\vert:=\min_{\delta'~\text{isotopic to}~\delta_j~\linebreak\text{rel}~P(g)\cup\{\infty\}}\vert\Tr(\delta')\cap\Tr(\gamma_i)\vert
\]
denote the minimal intersection number of $\delta_j$ and $\gamma_i$. It follows from the periodicity of the curves $\gamma_i$ and $\delta_j$ that $\vert[\delta_j]\cap[\gamma_0]\vert=0$ for all $j\in\{0,\ldots, k-1\}$. See \cite[Lemma~8.7]{BFH} for a proof of this fact. Therefore, the curves $\delta_j$ of the Levy cycle do not intersect the partition boundary (up to homotopy), so all post-singular points surrounded by the same curve $\delta_j$ have equal itineraries. This contradicts the fact that different post-singular points have different itineraries w.r.t.\ $(D_j)_{j\in\Z}$. See \cite[Lemma~5.5]{LSV} for a detailed proof.
\end{proof}

By Theorem~\ref{thm:TopChar} and Lemma~\ref{lem:LevyCycle}, the topological exponential map $g$ is Thurston equivalent to a holomorphic exponential map $E_{\lambda_g}$. We are in the position to close the classification cycle.

\begin{Thm}[Classification of post-singularly finite exponential maps by abstract Hubbard Trees]
The map
\[
\mathcal{F}\colon \{\text{psf exponential maps}\}\to\{\text{abstract exponential Hubbard Trees}\}
\]
is a bijection between post-singularly finite exponential maps (up to affine conjugation) and abstract exponential Hubbard Trees (up to equivalence).
\end{Thm}
\begin{proof}
We have already shown that $\mathcal{F}$ is well-defined. We choose a representative for every equivalence class of abstract exponential Hubbard Trees and perform the construction above to obtain a post-singularly finite exponential map associated with it. This sets up a map
\[
\mathcal{G}\colon \{\text{abstract exponential Hubbard Trees}\}\to\{\text{psf exponential maps}\}.
\]
Let $\mathtt{(H,f,B_{k_i},\angle_v)}$ be the chosen representative for a given equivalence class of abstract Hubbard Trees, and let $E_{\lambda_g}=\mathcal{G}([\mathtt{(H,f,B_{k_i},\angle_v)}])$ be the corresponding holomorphic exponential map obtained by the construction above. We want to see that $\mathcal{F}(E_{\lambda_g})=[\mathtt{(H,f,B_{k_i},\angle_v)}]$, i.e., $\mathcal{F}\circ\mathcal{G}=\text{id}$. As we have already seen in Theorem~\ref{thm:F_Injective} that $\mathcal{F}$ is injective, this shows that $\mathcal{F}$ is indeed a bijection.

By the definition of Thurston equivalence, there exist two homeomorphisms $\varphi_0,\varphi_1\colon \complex\to\complex$ such that the diagram

\[ \begin{tikzcd}
(\complex,P(g)) \arrow{r}{\varphi_1} \arrow[swap]{d}{g} & (\complex,P(E_{\lambda_g})) \arrow{d}{E_{\lambda_g}} \\%
(\complex, P(g)) \arrow{r}{\varphi_0}& (\complex,P(E_{\lambda_g}))
\end{tikzcd}
\]
commutes and $\varphi_0$ and $\varphi_1$ are homotopic rel $P(g)$. We want to see that $H_0:=\varphi_0(H)$ is a Homotopy Hubbard Tree for $E_{\lambda_g}$, and that it yields the abstract exponential Hubbard Tree that we have started with. Obviously, $H_0$ is a finite embedded tree spanned by $P(E_{\lambda_g})$. We want to see that $H_0$ is invariant up to homotopy rel $P(E_{\lambda_g})$. Set $H_1:=\varphi_1(H)$.
By Proposition~\ref{pro:homotopyImpliesIsotopy}, $\varphi_0$ and $\varphi_1$ are isotopic rel $P(g)$, hence the embedded trees $H_0$ and $H_1$ are ambient isotopic. More precisely, the map $\varphi:=\varphi_1\circ\varphi_0^{-1}\colon \complex\to\complex$ is isotopic to $\text{id}$ rel $P(E_{\lambda_g})$ such that $\varphi(H_0)=H_1$. We want to see that $H_0'=[P(E_{\lambda_g})]_{\widehat{E_{\lambda_g}}^{-1}(H_0)}$ is homotopic to $H_0$ in $\CT$ rel $P(E_{\lambda_g})$.

By construction, there exists a homotopy $I\colon H\times[0,1]\to\complex\cup\{-\infty\}$ between $H$ and $H'$ rel $P(g)$ satisfying $I^1=\iota'\circ\iota^{-1}$. The homeomorphism $\varphi_1$ extends to a homeomorphism $\widehat{\varphi_1}\colon \complex\cup\{-\infty\}\to\CT$ by setting $\widehat{\varphi_1}(-\infty):=-\infty$. By commutativity of the Thurston diagram, we have $H_0'=\widehat{\varphi_1}(H')$. The homotopy $I$ is pushed forward via $\widehat{\varphi_1}$ to a homotopy $\tilde{I}:=\widehat{\varphi_1}\circ I$ between $H_1$ and $H_0'$ rel $P(E_{\lambda_g})$. We have shown above that $H_1$ is homotopic to $H_0$ rel $P(E_{\lambda_g})$, so $H_0$ is homotopic to $H_0'$ rel $P(E_{\lambda_g})$.

The identification of $H_0$ with $H_0'$ yielded by the constructed homotopy is given by $\psi:=\widehat{\varphi_1}\circ\iota'\circ\iota^{-1}\circ\varphi_0^{-1}\colon H_0\to H_0'$. It remains to show that the induced self-map $f:=E_{\lambda_g}\circ\tilde{\psi}\colon H_0\to H_0$ is expansive. Looking carefully at the definition of the involved maps, we see that $f=(\varphi_0\circ\iota)\circ \mathtt{f}\circ(\varphi_0\circ\iota)^{-1}$, i.e., the self-map of $H_0$ is conjugate to the self-map $\mathtt{f}$ of the abstract Hubbard Tree $\mathtt{H}$. Hence, $f$ is expansive by the expansivity of $\mathtt{f}$.

This also shows that the abstract Hubbard Trees $(H_0,f,B_{k_i},\angle_v)$ and $\mathtt{(H,f,B_{k_i},\angle_v)}$ are equivalent via $\varphi_0\circ\iota$ except for the fact that the sector information of $\mathtt{H}$ and $H_0$ is consistent. For a post-singular point $\mathtt{p\in B_{k_i}}$, we have $\iota(\mathtt{p})\in D_{k_i}$, where $(D_i)_{i\in\Z}$ is the partition for $g$ defined above. As $\Tr(\gamma_h)\cap H=\{0\}$, we have $H_0\cap\Tr(\gamma)=\{0\}$, where $\gamma:=\varphi_0\circ\gamma_h$. By Proposition~\ref{pro:AccessesContainRays} and Lemma~\ref{lem:HomotopyImpliesAmbIsotopy}, the arc $\gamma$ is ambient isotopic rel $P(E_{\lambda_g})$ to a dynamic ray $g_{\ad{s}}$ landing at the singular value. Let $(\tilde{D}_i)_{i\in\Z}$ be the dynamical partition for $E_{\lambda_g}$ w.r.t.\ $g_{\ad{s}}$. As $\varphi_0$ restricts to a conjugation between $P(g)$ and $P(E_{\lambda_g})$, and as the ambient isotopy between $\gamma$ and $g_{\ad{s}}$ does not move post-singular points, we have $\varphi_0\circ\iota(\mathtt{p})\in\tilde{D}_{k_i}$. Hence, the sector information of $\mathtt{H}$ and $H_0$ is indeed consistent.
\end{proof}
This concludes the classification of psf exponential maps in terms of abstract exponential Hubbard Trees.

\section{Outlook}
\label{sec:outlook}

The key tools for our construction of Homotopy Hubbard Trees are dynamic rays, their landing properties, and the combinatorial encoding of both. It is natural to ask whether more general maps also have Homotopy Hubbard Trees. One putative obstacle is that the escaping set of an arbitrary post-singularly finite entire function need not consist of dynamic rays. Recently, the concept of \emph{dreadlocks} has been introduced in \cite{BR} as a generalization of dynamic rays. For every post-singularly finite entire function, the escaping set naturally decomposes into dreadlocks, and the collection of dreadlocks can be encoded in terms of external addresses. Every repelling periodic point is the landing point of at least one and at most finitely many dreadlocks, all of which are periodic of the same period. In short, dreadlocks possess all the important properties needed to study the branching of the Julia set in terms of the escaping set. Using these findings, the existence and uniqueness of Homotopy Hubbard Trees for all post-singularly finite transcendental entire functions has been established based on a study of branch points of the Julia set in \cite{Pf}. These results will be published in a series of forthcoming articles, the first of which is \cite{PPS}.

Another interesting question that has not been answered in \cite{Pf} is which post-singularly finite transcendental entire function have a Hubbard Tree in the classical sense, i.e., a compact embedded tree that is forward invariant as a set, not just forward invariant up to homotopy. By Theorem~\ref{thm:denseEscapingPoints}, psf exponential maps do not have a Hubbard Tree essentially because of the existence of an asymptotic value. We expect this to be the only obstacle to the existence of Hubbard Trees.

\let\oldaddcontentsline\addcontentsline
\renewcommand{\addcontentsline}[3]{}

\let\addcontentsline\oldaddcontentsline

\end{document}